\documentclass[11pt]{article}
\usepackage{a4}
\usepackage{santheorem}
\usepackage{amssymb,latexsym}
\usepackage[mathscr]{eucal}
\usepackage{citesort}
\usepackage{url}
\usepackage{deleq}
\ifx\macrosloaded\relax\endinput\else\let\macrosloaded\relax\fi

\newcommand{\U}{{\mcU}}
\newcommand{\mcU}{{\cal U}}
\newcommand{\mcS}{{\mycal I}}
\newcommand{\mcA}{{\mycal A}}
\newcommand{\mcO}{{\mycal O}}
\newcommand{\thm}[1]{{Theorem~\ref{#1}}}
\newcommand{\lem}[1]{{Lemma~\ref{#1}}}
\newcommand{\mca}{{\mcA_\infty^\delta}}

\newcommand{\bel}[1]{\begin{equation}\label{#1}}
\newcommand{\beal}[1]{\begin{eqnarray}\label{#1}}

{\catcode `\@=11 \global\let\AddToReset=\@addtoreset}
\AddToReset{equation}{section}
\renewcommand{\theequation}{\thesection.\arabic{equation}}

\newcommand{\be}{\begin{equation}}
\newcommand{\bea}{\begin{eqnarray}}
\newcommand{\eea}{\end{eqnarray}}
\newcommand{\beaa}{\begin{eqnarray*}}
\newcommand{\eeaa}{\end{eqnarray*}}
\newcommand{\bseq}{\begin{subeq}}
\newcommand{\eseq}{\end{subeq}}

\def \rectangle#1#2{\hbox{\vrule\vbox to #2
{\hrule\hbox to
#1{\hfil}\vfil\hrule}\vrule}}\def\square{\,\,\rectangle{7pt}{7
pt}\,\,}
\newcommand{\edd}{\end{document}}

\newcommand{\nn}{\nonumber}

\def\frac#1#2{{{#1}\over{#2}}}

\newcommand{\mxtauz }{(M_{x_1})}
\newcommand{\mxtau }{(M_{x_1-2\tau})}

\newcommand{\Mxonexzerot}{} 
\newcommand{\Mxonexzeros}{} 
\newcommand{\OmxT}{\Omega_{x_1,T}}
\newcommand{\BB}{{\mycal B}}
\newcommand{\cO}{{\cal O}}
\newcommand{\UU}{{\cal U}}
\newcommand{\cUxx}{\UU_{x_2,x_1}}
\newcommand{\hyp}{{\mycal S}}

\newcommand{\mxt}{(M_{x_1-2t})}
\newcommand{\mxs}{(M_{x_1-2s})}

     \newcommand{\setZ}{{\mathord{\mathbb Z}}}
     \newcommand{\lsemantics}{\mathopen{\lbrack\mkern-3mu\lbrack}}
     \newcommand{\rsemantics}{\mathclose{\rbrack\mkern-3mu\rbrack}}
\newcommand{\mcCzk}{{\mcC}^0_k}
\newcommand{\mcCaz}{{\mcC}^\alpha_0}
\newcommand{\Z}{\setZ}

\newcommand{\qed}{\hfill $\Box$\bigskip}
\newcommand{\proof}{\noindent {\sc Proof:\ }}
\newcommand{\remark}{\noindent {\bf Remark:\ }}
\newcommand{\remarks}{\noindent {\bf Remarks:\ }}
\newcommand{\eeq}{\end{equation}}
\newcommand{\ee}{\end{equation}}
\newcommand{\beqa}{\begin{eqnarray}}
\newcommand{\beqar}{\begin{deqarr}}
\newcommand{\eeqa}{\end{eqnarray}}
\newcommand{\eeqar}{\end{deqarr}}
\newcommand{\beqan}{\begin{eqnarray*}}
\newcommand{\eeqan}{\end{eqnarray*}}
\newcommand{\ba}{\begin{array}}
\newcommand{\ea}{\end{array}}

\newcommand{\decal}{{\mycal D}}
\newcommand{\mcC}{{\mycal C}}
\newcommand{\cH}{{\mycal H}}
\newcommand{\cB}{{\mycal B}}
\newcommand{\cG}{{\mycal G}}
\newcommand{\mcCak}{{\mcC}^\alpha_k}

\newcommand{\backmg}{h} 
\newcommand{\Id}{\mbox{\rm Id}} 
\newcommand{\const}{\mbox{\rm const}} 

\newcommand{\hide}[1]{}
\newcommand{\ovlocO}{\,\,\overline{\!\!\locO}}

\newcommand{\HH}{{\mycal H}}
\newcommand{\Hk}[1]{{\mycal H}_k^{#1}}
\newcommand{\Hak}{{\mycal H}_{k} ^{\alpha}}
\newcommand{\Haok}{{\mycal H}_{k} ^{-1/2}}
\newcommand{\Sigx}{\hyp_{x_2,x_1}}
\newcommand{\sigxx}[1]{\Sigma_{x_2,x_1,#1}}
\newcommand{\sigx}[1]{\Sigma_{x_1,#1}}
\newcommand{\demi}{{1\over 2}}
\newcommand{\emoins}{E_{-}^{\mu}}
\newcommand{\eplus}{E_{+}^\mu}
\newcommand{\nmu}{\nabla_\mu}
\newcommand{\pmu}{\partial_\mu}
\newcommand{\px}{\partial_x}

\newcommand{\xadu}{x^{-2\alpha-1+2\beta_1}}

\newcommand{\err}{\zeta} 
\newcommand{\prho}{p}
\newcommand{\dx}{\,dx}

\newtheorem{Theorem} {Theorem} [section]
\newtheorem{Corollary} [Theorem] {Corollary}
\newtheorem{Lemma} [Theorem] {Lemma}
\newtheorem{Proposition} [Theorem] {Proposition}

\DeclareFontFamily{OT1}{rsfs}{}
\DeclareFontShape{OT1}{rsfs}{m}{n}{ <-7> rsfs5 <7-10> rsfs7 <10->
rsfs10}{} \DeclareMathAlphabet{\mycal}{OT1}{rsfs}{m}{n}
\def\scri{{\mycal I}}%
\def\scrip{\scri^{+}}%
\newcommand{\proofend}{\qed}
\newcounter{mnotecount}[section]

\newcommand{\R}{\mathbb R}
\newcommand{\N}{\mathbb N}
\newcommand{\bM}{\overline{\! M}}

\newcommand{\backg}{b}
\newcommand{\pM}{\partial M}
\newcommand{\eq}[1]{(\ref{#1})}

\newcommand{\rspone}{\rangle_{\!{}_{1}}}
\newcommand{\rsptwo}{\rangle_{\!{}_{2}}}

\newcommand{\dv}{\;dv}
\newcommand{\Eqsone}[1]{Equations~\eq{#1}}
\newcommand{\Eq}[1]{Equation~\eq{#1}}
\newcommand{\Eqs}[2]{Equations~\eq{#1}-\eq{#2}}
\newcommand{\dea}{{D}_{e_{A}}}
\newcommand{\loc}{{\mbox{{\rm\scriptsize loc}}}}

\newcommand{\boa}{{\mycal B} _0 ^\alpha}
\newcommand{\cob}{\mcC _0 ^\beta}

\newcommand{\GG}{{\mycal G}}
\newcommand{\Gk}[1]{\GG_k ^{#1}}
\newcommand{\Gak}{\Gk \alpha}
\newcommand{\gok}{\Gk 0}
\newcommand{\Gbk}{\Gk \beta}
\newcommand{\tM}{\,\,\,\,\,\widetilde{\!\!\!\!\!\mycal M}}
\newcommand{\dmu}{d\mu} 
\newcommand{\volu}{\,d \nu} 
\newcommand{\dnu}{\volu} 

\newcommand{\ooi}{\Omega_i} 

\newcommand{\cM}{\mycal M}
\newcommand{\cN}{\mycal N}

\newcommand{\cAp}{{\mycal A}_{\mbox{\scriptsize phg}}}

\newcommand{\cApM}{\cAp(M)}
\newcommand{\stsg}{{\mathfrak g}}
\newcommand{\tf}{\widetilde f}
\newcommand{\hf}{\widehat f}
\newcommand{\hphi}{\widehat \phi}

\newcommand{\complementaire}{\complement}

\newcommand{\locO} {{\mycal O}}

\newcommand{\beqd}{\begin{deqarr}}
\newcommand{\eeqd}{\end{deqarr}}

\begin{document}
\title{Solutions of wave equations in the radiation regime}
\author{Piotr T.\ Chru\'sciel\thanks{ Supported in part by the Polish
    Research Council grant KBN 2 P03B 130 16, and by an
    A.~von~Humboldt fellowship. E-mail:
    \texttt{chrusciel@univ-tours.fr}; URL
    {\protect\url{http://www.phys.univ-tours.fr/}$\sim$\protect\url{piotr}}.}
    \\ O. Lengard\thanks{E-mail:
    \texttt{olengard@infonie.fr}} \\ D\'epartement de
    Math\'ematiques \\ Facult\'e des Sciences\\ Parc de Grandmont\\
    F-37200 Tours, France}

\maketitle

\vfill\begin{abstract}
  We study the ``hyperboloidal Cauchy problem'' for linear and
  semi-linear wave equations on Minkowski space-time, with initial
  data in weighted Sobolev spaces allowing singular behaviour at
  the boundary, or with polyhomogeneous initial data. Specifically,
  we consider nonlinear symmetric hyperbolic systems of a form which
  includes scalar fields with a $\lambda\phi^p$ nonlinearity, as well
  as wave maps, with initial data given on a hyperboloid; several of the results proved
  apply to general space-times admitting conformal completions at null infinity, as well to
  a large class of equations with a similar non-linearity structure.
  We prove existence of solutions with controlled
  asymptotic behaviour, and asymptotic expansions for solutions when the initial data
  have such expansions. In particular we prove that polyhomogeneous
  initial data (satisfying compatibility conditions) lead to solutions
  which are polyhomogeneous at the conformal boundary $\scrip$ of the
  Minkowski space-time.
\end{abstract}
\eject
\tableofcontents

\section{Introduction}
Bondi \emph{et al.} \cite{BBM} together with Sachs \cite{Sachs}
and Penrose \cite{penrose:scri}, building upon the pioneering work
of Trautman \cite{T,Tlectures}, have proposed in the sixties a set
of boundary conditions appropriate for the gravitational field in
the radiation regime. A somewhat simplified way of introducing the
Bondi-Penrose (BP) conditions is to assume existence of
``asymptotically Minkowskian coordinates'' $(x^\mu)=(t,x,y,z)$ in
which the space-time metric $\stsg  $ takes the form
\begin{equation}
  \label{eq:1}
  {\stsg }_{\mu\nu}-\eta_{\mu\nu}=
  \frac{\stackrel{\phantom{x}_1}{h}_{\mu\nu}(t-r,\theta,\varphi)}{r} +
  \frac{\stackrel{\phantom{x}_2}{h}_{\mu\nu}(t-r,\theta,\varphi)}{r^2} + \ldots\;,
\end{equation}
where $\eta_{\mu\nu}$ is the Minkowski metric diag$(-1,1,1,1)$,
$u$ stands for $t-r$, with $r,\theta,\varphi$ being the standard
spherical coordinates on $\R^3$. The expansion above has to hold
at, say, fixed $u$, with $r$ tending to infinity. Existence of
classes of solutions of the vacuum Einstein equations satisfying
the asymptotic conditions \eq{eq:1} follows from the work in
\cite{friedrich:cauchy} together with
\cite{AndChDiss,CorvinoSchoen,Corvino,ACF}.  As of today it
remains an open problem how general, within the class of radiating
solutions of vacuum Einstein equations, are those solutions which
display the behaviour \eq{eq:1}.  Indeed, the results in
\cite{AndChDiss,AC,ChMS,PRLetter,ACF,ChMS} suggest
strongly\footnote{\emph{Cf.}  \cite{Kroon1} and references therein
for
  some further related results.} that a more appropriate setup for such
gravitational fields is that of \emph{polyhomogeneous} asymptotic
expansions:
\begin{equation}
  \label{eq:2}
  {\stsg }_{\mu\nu}-\eta_{\mu\nu}\in \cAp\;.
\end{equation}
In the context of expansions in terms of a radial coordinate $r$
tending to infinity, the space of polyhomogeneous functions is
defined as the set of smooth functions which have an asymptotic
expansion of the form
\begin{equation}
  \label{eq:3}
  f \sim \sum_{i=0}^\infty \sum_{j=0}^{N_i} f_{ij}(u,\theta,\varphi) \frac{\ln
  ^j r}{r^{n_i}}\;,
\end{equation}
for some sequences $n_i,N_i$, with $n_i\nearrow\infty$. Here the
symbol $\sim$ stands for ``being asymptotic to'': if the
right-hand-side is truncated at some finite $i$, the remainder
term falls off appropriately faster. Further, the functions
$f_{ij}$ are supposed to be smooth, and the asymptotic expansions
should be preserved under differentiation.\footnote{ The choice of
the sequences $n_i,N_i$ is not arbitrary, and is dictated by the
equations at hand. For example, the analysis of $3+1$ dimensional
Einstein equations in \cite{ChMS} suggests that consistent
expansions can be obtained with $n_i=i$. On the other hand,
Theorem~\ref{Twavemap} below gives actually $n_i=i/2$ for
wave-maps on $2+1$ dimensional Minkowski space-time.  We note that
the $2+1$ dimensional wave map equation is related to the vacuum
Einstein equations with cylindrical symmetry (\emph{cf., e.g.},
\cite{BChM,CT93,CT293}).}

The suggestion, that the expansions \eq{eq:2} are better suited
for describing the gravitational field in the radiation regime
than \eq{eq:1}, arises from the fact that {\em generic} -- in a
well defined sense -- initial data constructed in
\cite{AndChDiss,AC,ChMS,PRLetter,ACF,ChMS} are polyhomogeneous.
This leads naturally to the question, whether polyhomogeneity of
initial data is preserved under evolution dictated by wave
equations. In this paper we answer in the affirmative this
question for semi-linear wave equations, and for the wave map
equation, on Minkowski space-time. We develop a functional
framework appropriate for the analysis of such questions.
We prove local in time existence of solutions for classes of
equations that include the semi-linear wave equations and the wave
map equation on Minkowski space-time, with conormal and with
polyhomogeneous initial data. We show that polyhomogeneity is
preserved under evolution when appropriate (necessary) corner
conditions are satisfied by the initial data. We note that the
initial data considered here are  more singular than allowed in
the existing related results \cite{Bony,MelroseRitter,Joshi}. We
are planning to analyse the corresponding problems for the vacuum
Einstein equation in a forthcoming publication, see also
\cite{OLthese}.

Our main results are the existence and polyhomogeneity of
solutions with appropriate polyhomogeneous initial data for the
nonlinear scalar wave equation, and for the wave map equation. We
achieve this in a few steps. First, we prove local existence of
solutions of these equations in weighted Sobolev spaces, {\em
cf.}\/ Theorems~\ref{T2} and \ref{T2w}. The next step is to obtain
estimates on the time derivatives, {\em cf.}\/ Theorems~\ref{T2t},
\ref{T2wt} and \ref{2dT2wt}. Those estimates are uniform in time
in a neighbourhood of the initial data surface if the initial data
satisfy compatibility conditions. Somewhat surprisingly, we show
that all initial data in weighted Sobolev spaces, not necessarily
satisfying the compatibility conditions, evolve in such a way that
the compatibility conditions will hold on all later time slices;
see Corollary~\ref{C1t} and Theorems~\ref{T2wt} and \ref{2dT2wt}.
Finally, in Theorems~\ref{T2phg} and \ref{Twavemap} we prove
polyhomogeneity of the solutions with polyhomogeneous initial
data; this requires a hierarchy of compatibility conditions. We
hope to be able to show in a near future that polyhomogeneity of
solutions can be established, for polyhomogeneous initial data,
with a finite number of compatibility conditions.

The restriction to Minkowski space-time in Theorem~\ref{Twavemap}
is not necessary, and is only made for simplicity of presentation
of the results; the same remark applies to Theorem~\ref{T2}.
Similarly the choice of the initial data hypersurface as the
standard unit hyperboloid is not necessary.

This work is organised as follows: First, the reader is referred
to  Appendix~\ref{S2} for definitions, notations, and the
functional spaces involved; we also develop  calculus in those
spaces there. In Section~\ref{ss1} we briefly recall Penrose's
conformal completions, as they provide the link between the
asymptotic behavior of fields and the local analysis carried on in
this work. In Section~\ref{S3} we consider linear equations. There
the key elements of our analysis are: a) Proposition~\ref{PL.1}
and its variations, which give {\em a priori} estimates in
weighted Sobolev spaces; b) the mechanism for proving
polyhomogeneity, provided in the proof of Theorem~\ref{Tlemme1}.
The transition from the linear weighted Sobolev estimates to their
nonlinear counterparts is done in Sections~\ref{sslwe} and
\ref{Swave}. This has already been outlined above, and requires a
considerable amount of work. In Appendix~\ref{SODEsws} we prove
several auxiliary results on ODE's, some of which are fairly
straightforward; as those results are  used in the body of the
paper in various, sometimes involved, iterative arguments, it
seemed convenient to have precise statements at hand.

Some of the results proved here have been announced in
\cite{ChLNantes}.

\newcommand{\phg}{\mbox{\scriptsize \rm phg}}
\section{Conformal completions}\label{ss1}
The aim of this section is to set-up the framework necessary for
our considerations; the results here are well known to
relativists, but perhaps less so to the PDE community. In any case
they are needed to establish notation. Consider, thus, an $n+1$
dimensional space-time $({\mycal M},\stsg)$ and let
\begin{equation}
  \label{C.1}
\tilde{\stsg} = \Omega^{2}\stsg \;.
\end{equation}
Let $\Box _h$ denote the wave operator associated with a
Lorentzian metric $h$,
$$\Box_h f= {1\over\sqrt{|\det h_{\rho\sigma}|}} \pmu (\sqrt{|\det
  h_{\alpha\beta}|}h^{\mu\nu} \partial_\nu f).$$
We recall that the scalar curvature $R =R(\stsg)$ of $\stsg $ is
related to the corresponding scalar curvature $\tilde{R} =
\tilde{R}(\tilde{\stsg })$ of $\tilde{\stsg }$ by the formula
\begin{equation}
\tilde{R}\Omega^2 = R -2n \left\{ {1\over
\Omega}\Box_{\stsg}\Omega +{n-3 \over
  2} {|\nabla\Omega|^2_\stsg  \over \Omega^2}\right\}\;.
\label{C.2}
\end{equation}
It then follows from (\ref{C.2}) that we have the identity
\begin{equation}
\Box_{\tilde{\stsg }} (\Omega^{-{n-1\over 2}}f)=
\Omega^{-{n+3\over 2}}\left( \Box_\stsg  f +{n-1 \over
4n}(\tilde{R}\Omega^2 -R)f\right)\;. \label{C.3}
\end{equation}
It has been observed by Penrose~\cite{penrose:scri} that the
Minkowski space-time $({\mycal M},\eta)$ can be conformally
completed to a space-time with boundary $(\tM,\tilde{\eta})$,
$\tilde{\eta}=\Omega^{-2} \eta$ on ${\mycal M}$, by adding to
${\mycal M}$ two null hypersurfaces, usually denoted by $\scrip$
and $\scri^-$, which can be thought of as end points ($\scrip$)
and initial points ($\scri^-$) of inextendible null geodesics
\cite{NewmanAF,WaldBook,penrose:scri}. We will only be interested
in ``the future null infinity'' $\scrip$; an explicit construction
(of a subset of $\scrip$) which is convenient for our purposes
proceeds as follows: for $(x^0) ^2 <\displaystyle \sum_{i}(x^i
)^2$ we define \be y^\mu = {x^\mu\over x^\alpha x_\alpha} \;
;\label{C.4} \ee in the coordinate system $\{ y^\mu\} $ the
Minkowski metric $\eta \equiv -(dx^0)^2 +(dx^1)^2 +(dx^2)^2
+(dx^3)^2
= \eta_{\alpha\beta}dx^\alpha dx^\beta$ takes the form \beqa &\eta
= \displaystyle{1\over \Omega^2} \eta_{\alpha\beta} dy^\alpha
dy^\beta \;,\label{C.5}\qquad
\Omega =  \eta_{\alpha\beta}y^\alpha y^\beta \;.&
\eeqa We note that under (\ref{C.4}) the exterior of the light
cone $C_0^{x^\mu} \equiv\{ \eta_{\alpha\beta} x^\alpha x^\beta
=0\} $ emanating from the origin of the $x^\mu$-coordinates is
mapped to the exterior of the light cone $C_0^{y^\mu} = \{
\eta_{\alpha\beta}y^\alpha y^\beta =0 \} $ emanating from the
origin of the $y^\mu$-coordinates. The conformal completion is
obtained by adding $C_0^{y^\mu}$ to ${\mycal M}$,
$$\tM = {\mycal M} \cup (C_0^{y^\mu}\setminus \{0\} )\;,$$
with the obvious differential structure arising from the
coordinate system $y^\mu$.  We shall use the symbol $\scri$ to
denote $C_0^{y^\mu}\setminus \{0\}$, and $\scrip$ to denote
$C_0^{y^\mu}\setminus \{0\}\cap \{ y^0>0\} $. As already
mentioned, $\scri$ so defined is actually a subset of the usual
$\scri$, but this will be irrelevant for our purposes.

We note that (\ref{C.4}) is singular at the light cone
$C_0^{x^\mu}$. This is again irrelevant from our point of view
because we are only interested in the behaviour of the solutions
near $\scrip$, and finite speed of propagation allows us, for that
purpose, to disregard what happens near $C_0^{x^\mu}$.

The above procedure can be adapted for several metrics of
interest, such as the Schwarzschild, Kerr, or Robinson-Trautman
metrics, to similarly yield conformal completions of space-time by
the addition of null hypersurfaces $\scrip$. This observation was
at the origin of Penrose's proposal to describe systems which are
asymptotically flat in lightlike directions through the use of
conformal completions.

It is noteworthy that the conformal technique allows one to reduce
global-in-time existence problems to local ones; this has been
exploited by various authors
\cite{ChBPisa,ChristodoulouCPAM,ChBGu,ChBdeSitter,ChBglwa,ChBNou}
for wave equations on a fixed background space-time. Further,
Friedrich \cite{Friedrich,FriedrichSchmidt,HelmutJDG} has used
this approach to obtain a global existence result for Einstein
equations to the future of a ``hyperboloidal'' Cauchy surface,
with ``small'' smoothly conformally compactifiable initial data,
\emph{cf.} also \cite{FriedrichdS,Friedrich:aDS,Friedrich:Pune}.

On a more modest level, the identity \eq{C.3} can be used as a
starting point for the analysis of the asymptotic behaviour of
solutions of the scalar wave equation near $\scrip$, as it reduces
the problem to a study of solutions near a null hypersurface. This
is the approach used in this paper.  There are associated
identities for fields of any spin \cite{penrose:scri}, which
provide a convenient framework for similar questions for those
fields.

\section{A  class of linear symmetric   hyperbolic systems} \label{S3}

In this section we  define a
  class of linear symmetric hyperbolic first order systems on a set of
  the form $ M_{x_0}\times I$, where $I$ is an interval corresponding
  to the time variable, which will be denoted by $\tau$, and
  we derive our key energy inequality
  in \emph{weighted}
 Sobolev spaces. (We note that
  in some of our further applications the vector $\partial/\partial
  \tau$ will be lightlike, and not timelike as is usually the case.
  It should be pointed out that in our conventions the time variable
  is  the last coordinate, allowing $x$ to be the first variable,
  consistently with  the conventions of the preceding sections.)
 We start by introducing some notation for the sets within the
 ``space-time'' $M_{x_0}\times I$, which will be relevant in what
 follows\footnote{The motivation for the factors of $2$, and the general form
  of the sets   considered, arises as follows: The set $\partial
  M\times I$ should be thought of as a smooth null hypersurface in
  space-time; \emph{e.g.}, in Minkowski
  space-time with Minkowskian coordinates $y^\mu$, it can be the
  intersection of the half-space $\{y^0\geq\frac{1}2\} $ with the   light cone
  emanating from   the origin $y^{\mu}=0$ . Then $\tau$ is the
  Minkowski time, perhaps shifted by a constant, say
  $\tau=y^0-\frac{1}{2}$. The coordinate $x$ is  a coordinate
  which vanishes on $\partial M\times I$, in the current example \emph{e.g.}
  $x=\sqrt{\sum (y^i)^2}-y^0$. Finally, in such a Minkowskian setup,
  the hypersurfaces $x=x_1-2\tau$,
  which determine  one of the boundaries of the $\Sigma$'s and
  $\Omega$'s defined in \eq{L.7al}, correspond to the
  converging light cones $y^0+\sqrt{\sum
    (y^i)^2}=\const$. The restrictions $2(x_2+t)<x_1\leq x_0$  (in the
  definition of $\Sigma_{x_2,x_1,t}$)  and $2(x_2+T)<x_1$ (in the
  definition of $\Omega_{x_2,x_1,T}$) are not   necessary, and are
  only made    for simplicity of discussion.}:
 \begin{deqarr} \label{L.7.0}& t\geq 0,
\quad 2(x_2+t)<x_1\leq x_0,\qquad \Sigma_{x_2,x_1,t} = \{
\tau=t, \ x_2< x < x_1 -2t\}, &\nn \\ &&\\
& T>0,\quad 2(x_2+T)<x_1\leq x_0,\qquad \Omega_{x_2,x_1,T}
=\bigcup_{0<
  \tau< T} \Sigma_{x_2,x_1,\tau},  & \\
&  0\leq 2t<x_1\leq x_0\;,\qquad \sigx t = \{ \tau=t, 0< x<
x_1-2t\}, &\\
& 0<2T<x_1, \qquad \Omega_{x_1,T} = \bigcup_{0< t< T} \sigx t \;.&
\label{L.7}\arrlabel{L.7al}\end{deqarr} There is a natural
identification between $\Sigma_{x_2,x_1,t}$ and $M_{x_2,x_1-2t}$,
similarly between $\sigx t$ and $M_{x_1-2t}$,
and we shall freely make use of such identifications throughout.
We shall write $\|f(t)\|_{\Hak}$ for $ \|f|_{\sigxx t}\|_{\Hak
(\sigxx t
  )}$, or for $ \|f|_{\sigx t}\|_{\Hak (\sigx t )}$, \emph{etc.}; the
distinction should be clear from the context.

We shall be interested in symmetric hyperbolic first order systems
which in local coordinates take the form
\begin{equation} \left[A^\mu (z^\nu)\partial_\mu + A(z)\right]
f=F,\label{L.1}
\end{equation}
where $z^\nu = (y^i,\tau)$, with the following properties:

$\mcC 1$) 
$f$ and $F$ are sections of a bundle which is a direct sum of two
$N_1$ and $N_2$ dimensional Riemannian bundles; we will write
\begin{equation}
f=\left(\begin{array}{c}\varphi \\ \psi \end{array}\right)
\;,\qquad F=\left(\begin{array}{c} a \\ b \end{array}\right) \; .
\label{L.2} \end{equation} In local coordinates $\varphi$ and $a $
are thus $\R^{N_1}$ valued, while $\psi$ and $b $ are $\R^{N_2}$
valued. The respective scalar products will be denoted by
$\langle\cdot,\cdot\rspone$ and $\langle\cdot,\cdot\rsptwo$. We
shall use the generic symbol $\nabla$ to denote\footnote{In some
situations
  \eq{eq:covder} might fail to hold, and some undifferentiated
  supplementary terms will occur at the right-hand-side of
  \eq{eq:covder}. We note that our results will not be affected by the
  occurrence of such terms, provided those terms satisfy bounds as in
  \eq{L.6.2}.} a covariant derivative compatible with those scalar
products, \emph{e.g.}, if $X$ is a vector field on
$\Omega_{x_0,T}$, then
\begin{equation}
  \label{eq:covder}
  X(\langle\phi,\psi\rspone)=\langle\nabla_X\phi,\psi\rspone +
\langle\phi,\nabla_X\psi\rspone\;,
\end{equation}
similarly for $\langle\cdot,\cdot\rsptwo$. $\nabla$ will also be
assumed to be compatible with every other structure at hand
whenever useful in the context, \emph{e.g.} a Riemannian metric on
$M$, \emph{etc.}

$\mcC2$) The left hand side of (\ref{L.1}) can be written as
\begin{equation}
\left(\begin{array}{ll} E^\mu_- \nabla _\mu \varphi& + L\psi \\ -
L^\dagger \varphi &+ E^\mu_+ \nabla_\mu \psi \end{array}\right)+
\left(\begin{array}{cc}B_{11} & B_{12} \\ B_{21} &B_{22}
\end{array}\right) \left(\begin{array}{c}\varphi\\\ \psi
\end{array}\right)\;,
\label{L.3} \end{equation} where $L$ is a first order differential
operator. Here $L^\dagger$ denotes the formal adjoint of $L$, in
the sense that if $\Omega = M$, or $M_{x_1}$, or $M_{x_2,x_1}$,
and if $\varphi,\psi$ are in $C_1 ({\overline{\Omega}})$, then
\begin{equation}\int_{\Omega} \langle \varphi, L \psi \rspone \;\dmu=
\int_{\Omega} \langle L^\dagger \varphi,\psi
\rsptwo\;\dmu\;,\label{L.4}
\end{equation} where $\dmu$ is a measure on $M$ which will, we hope,
be obvious from the context. By density Equation (\ref{L.4}) will
still hold with $\Omega=M_{x_2,x_1}$ for all $\alpha,\beta\in\R$,
all $\varphi \in {{\HH }^\alpha_1}(M_{x_2,x_1})$ and all $ \psi
\in {{\HH }^\beta_1}(M_{x_2,x_1})$. Equation (\ref{L.4}) forces
$L$ not to contain any $\tau$- or $x$- derivatives, where the
letter $x$ denotes a coordinate as defined in Section~\ref{S2},
thus
\begin{equation}
\label{eq:L.4.1} L=\ell^A(x,v,\tau)\partial_A + \ell(x,v,\tau)\;.
\end{equation}
 It follows that the principal part of the system \eq{L.3} is of
the form \be \label{sheq}\left(\begin{array}{ll} E^\mu_-
\partial_\mu & \ell^A\partial_A \\ (\ell^A)^t\partial_A & E^\mu_+
\partial_\mu
\end{array}\right)\;, \ee where $A^t$ denotes the transpose of a matrix
$A$. \Eq{sheq} explicitly shows that \eq{L.3} is symmetric
hyperbolic when the $E^\mu_{\pm}$'s are symmetric with
$E^\tau_{\pm}$ positive definite; the notions of ``symmetric
hyperbolic'' and ``symmetrizable hyperbolic'' are identified
throughout this work.

The hypotheses above will be assumed throughout this section.
\subsection{Estimates on the space derivatives of the solutions, $\alpha<-1/2$}
 Let us
pass now to the description of the hypotheses needed to derive
weighted energy estimates for space derivatives of $f$. To obtain
such estimates, we shall require the existence of a constant $C_1$
such that the (matrix-valued) coefficients $\ell^A$ and $\ell$
satisfy, in the relevant range of $\tau$'s,
\begin{equation}
\label{eq:L.4.2}
\|\ell(\tau)\|_{\gok\mxtau } +\sum_A\|\ell^A(\tau)\|_{\gok\mxtau } \leq C_1 
\;.
\end{equation}
Similarly writing
\begin{equation}
\label{eq:L.4.3} L^\dagger=\ell^{\dagger A}(x,v,\tau)\partial_A +
\ell^\dagger(x,v,\tau)\;,
\end{equation}
we require
\begin{equation}
\label{eq:L.4.4}
\|\ell^{\dagger}(\tau)\|_{\gok\mxtau } +\sum_A\|\ell^{\dagger A}(\tau)\|_{\gok\mxtau } \leq C_1 
\;.
\end{equation}

$\mcC3$) The matrices $E^\mu_{\pm}$ are symmetric and  satisfy
\begin{equation}
E^\mu_\pm n_\mu \geq \varepsilon \Id\;,\qquad E^\mu_+ \partial_\mu
x \leq - \varepsilon \Id\;,\qquad |E^\mu_- \partial_\mu x| \leq
C_1x\;, \label{L.5} \end{equation} for some $\varepsilon > 0$.
Here $n_\mu$ denotes the field of future directed (\emph{i.e.},
$\backg(d\tau,n)>0$) $\backg$-unit normals to the surfaces
$\{\tau=\const\}$, where $\backg$ is an auxiliary Riemannian
metric $\backmg$ on $M$. (Later on we will mainly be interested in
the case of $E^\mu_+$s of the form $E^\mu_\pm = e^\mu_\pm \otimes
\Id$, for some vector fields $e^\mu_\pm$.) For simplicity we shall
also assume
\begin{equation}
  \label{eq:simpl}
  \partial_i E^{\tau}_\pm = 0\;;
\end{equation}
this is by no means necessary, but is sufficient for the purposes
of this paper. We will further assume\footnote{\label{densfoot}We
use a convention in which the covariant derivatives $\nmu
E^\mu_{\pm}$ include terms associated with the vector density
character of $X^\mu$ defined by \eq{L.11}; in particular this
should be taken into account when verifying that the estimates
\eq{L.6.0}-\eq{L.8.n} hold.} that the $E^\mu_ -$'s satisfy a bound
of the form:
\begin{eqnarray}
\| E^A_- (\tau)\|_{{\cG}^0_k\mxtau }+ \| E^x_-
(\tau)\|_{{\cG}^1_{k}\mxtau }+\|\nabla_\mu  E^\mu_ -(\tau)
\|_{L^\infty\mxtau } & \leq & C_1\;. \nn\\
&&\label{L.6.0}\end{eqnarray}
As far as the $ E^\mu_+$'s are concerned, we allow singular
behaviour which should, however, be somewhat less singular than
$1/x$; to control that, we require existence of a function
$\err:\R^+\to\R^+$, satisfying $\lim_{x\to0}\err(x)=0 $, such that
for $0<x\leq x_1-2\tau$ we have
\begin{equation}
\| E^A_ + (\tau)\|_{{\cG}^{-1}_k(M_{x})}+\|  E^x_ +
(\tau)\|_{{\cG}^{0}_{k}(M_{x})}+ \|x\nabla_\mu E^\mu_+(\tau)
\|_{L^\infty(M_{x})}\leq  \err(x)\;.\label{L.8.n}\end{equation}
When the operators $E^\mu_\pm \nabla_\mu$ are written out
explicitly as
\begin{equation}
  \label{L.6.1}
  E^\mu_\pm \nabla_\mu = E^\mu_\pm \partial_\mu + B_\pm\;,
\end{equation}
we require that
\newcommand{\mx}{(M_x)}
\begin{equation}
  \label{L.6.2}
  \|B_-(\tau)\|_{\gok\mxtau } \leq C_1\;, \quad
\|B_+(\tau)\|_{{\cG}^{-1}_{k}\mx}\leq  \err(x)\;, \quad
0<x<x_1-2\tau\;.
\end{equation}

$\mcC4$) The matrices $B_{ab}$, $a,b=1,2$, satisfy the bounds
\begin{eqnarray}
& \|B_{12}(\tau)\|_{{\cG}^{-1/2}_k\mxtau}
+\|B_{21}(\tau)\|_{{\cG}^{-1/2}_k\mxtau}+\|B_{11}(\tau)\|_{{\cG}^0_k\mxtau
}
 \leq
C_1\;, & \nonumber \\ & \|B_{22}(\tau)\|_{{\cG}^{-1}_k\mx}\leq
\err(x)\;, & \label{L.6}\end{eqnarray} this last equation holding
again for $0<x<x_1-2\tau$.

Our final hypothesis concerns the ``acausal'' nature of the
boundary of $\Omega_{x_2,x_1,T}$:

$\mcC5$) $\partial \Omega_{x_2,x_1,T}$ is ``non-timelike'', in the
sense that for any covector $n_\mu$,
which is positive on outwards-pointing vectors and vanishes on
vectors tangent to $\partial \Omega_{x_2,x_1,T}$  we have, on
$\partial \Omega_{x_2,x_1,T}\cap \{\tau
> 0\} $,
\begin{equation}
  \label{eq:L.5.1}
  E^\mu_\pm n_\mu \geq 0\;.\end{equation}
(We note that \eq{L.5} already guarantees that \eq{eq:L.5.1} holds
on $\partial \Omega_{x_2,x_1,T}\cap \{\tau =T \mbox{\ or\ }
x=0\}$.)

 The essential point of the above hypotheses is
that the boundary $\{x=0\}$ is characteristic for \Eq{L.1}, with
the dimension of the relevant kernel being constant over the
boundary.\footnote{We are grateful to H.~Friedrich for useful
discussions concerning this point.} Weighted estimates,  in the
spirit of Proposition~\ref{PL.1} below, near such characteristic
boundaries hold for general symmetric hyperbolic systems, this
will be discussed elsewhere.

Weighted energy inequalities in $\Hak$ spaces with arbitrary
values of $k$ may be proved under various hypotheses on the
coefficients which appear in \eq{L.1}.  We note one such result
for systems satisfying $\mcC1)$-$\mcC5)$, which lies in line with
our remaining investigations. The restriction $\alpha\le-1/2$
seems to be inherent to the problem at hand. We consider first the
case $\alpha<-1/2$; the case $\alpha=-1/2$ is handled by the same
methods, under somewhat more restrictive conditions on the
coefficients, in Section~\ref{Sonehalf}.

\begin{Proposition} \label{PL.1}Suppose that
  $\alpha<-\demi$, $k>{n\over 2}+1$, $ k\in \N$, and set either
  $f(t)=f|_{\sigx t}$, $0<x_1\leq x_0$, $0\leq t \leq t_{\max}\equiv
  x_1/2$, or $f(t)=f|_{\sigxx t}$, $0<2x_2< x_1\leq x_0$, $0\leq t <
  t_{\max} \equiv x_1-2x_2$.  Under
  the  hypotheses $\mcC1)$-$\mcC5)$,
  there exists a constant $C_2$ depending upon $x_1$, $C_1$, $n$, $N$,
  $k$ and $\alpha$, as well as upon the ``error function''
  $\err$ and the boundary manifold $\partial M$, such that for all
$f$ satisfying~\eq{L.1} for which $f(0)\in
H_{k}^{loc} $ and for all $
0<t\leq t_{\max}$ we have
\begin{eqnarray}
\|f(t)\|^2_{\Hak\mxt} & \leq &
C_2e^{C_2t}\left(\|f(0)\|^2_{\Hak(M_{x_1})} +
  \int_0^t
  e^{C_2(t-s)}\left(\|a(s)\|^2_{\Hak\mxs}\right. \right.
\nonumber \\ & &  
+\left.\|b(s)\|^2_{\Hk {\alpha-1/2}\mxs}\Big)
  ds\right)\;,
 \label{L.10}\end{eqnarray}
with an identical estimate with $M_{x_1-*}$ replaced by
$M_{x_2,x_1-*}$.
\end{Proposition}

\remark The condition $k>n/2+1$ is needed to obtain
$C_1$--weighted control of the solution; there are no restrictions
on $k$ if we have at our disposal an  \emph{a  priori\/} $C_1$
weighted bound for $f$, and if the coefficients in the equation
are suitably regular. In such a case, for $k\leq n/2+1$, the
inequality \eq{L.10} should be modified by adding a term
$\|f(s)\|^2_{{{\cB}^\alpha_{1}(M_{x_1-2s})}}$ under the integral
appearing in  \eq{L.10}.

\medskip

\proof We start by proving \eq{L.10} on sets $M_{x_2,x_1-t}$; in
that case we are mainly interested to obtain uniform control for
small values of $x_2$, with eventually $x_2$ tending to zero;
without the uniformity in $x_2$ the estimate would of course be
standard. Keeping this in mind, let $X^\mu$ be the
``energy-momentum vector density'', \be X^\mu = \sum_{0\leq
  |\beta|\leq k}\xadu\{ \langle \decal^\beta\varphi,\emoins
\decal^\beta\varphi\rspone +\langle \decal^\beta\psi,\eplus
\decal^\beta\psi\rsptwo\}. \label{L.11} \ee Suppose, first, that
$f(0)\in H_{k+1}^{loc} $; standard results \cite[Vol.~III]{Taylor}
show that  $f(t)\in H_{k+1}^{loc} $, and we then
have$^{\mbox{\scriptsize \ref{densfoot}}}$ \be \nmu X^\mu = N_1
+D_1 + D_2 + E_1 + E_2 + E_3\;, \quad \label{L.12} \ee where
\begin{eqnarray}
N_1&=& \sum_{0\leq|\beta|\leq k}  {(2\beta_1-2\alpha-1)}x^{-2\alpha-2+2\beta_1} \langle \decal^\beta \psi,(\eplus\pmu x)\decal^\beta \psi\rsptwo\;,\nonumber\\
D_1&=& 2\sum_{0\leq|\beta|\leq k}\xadu \langle \decal^\beta \varphi,\emoins\nmu\decal^\beta \varphi\rspone\;,\nonumber\\
D_2&=& 2\sum_{0\leq|\beta|\leq k}\xadu \langle \decal^\beta \psi,\eplus\nmu\decal^\beta \psi\rsptwo\;,\nonumber\\
E_1&=& \sum_{0\leq|\beta|\leq k} (2\beta_1-2\alpha-1)\xadu
\langle\decal^\beta \varphi\;,{(\emoins\pmu x)\over x} \decal^\beta \varphi\rspone\;,\nonumber\\
E_2&=& \sum_{0\leq|\beta|\leq k}\xadu \langle\decal^\beta \varphi,(\nmu\emoins)\decal^\beta \varphi\rspone\;,\nonumber\\
E_3&=& \sum_{0\leq|\beta|\leq k} \xadu \langle\decal^\beta
\psi,(\nmu\eplus)\decal^\beta \psi\rsptwo\; .  \label{L.13}
\end{eqnarray}
Since $2\alpha+1 <0$, from (\ref{L.5}) one 
finds that
\be
\int_{\sigxx s}N_1 \dx\volu \leq -|2\alpha+1|\varepsilon
\|\psi\|^2_{\Hk {\alpha+\demi}} \label{L.14} \ee which is strictly
negative except if $\psi$ is identically zero, and can be used to
control some of the error terms which occur at the right hand side
of (\ref{L.12}). (Here we have used the form \eq{S2.0x} of
$\|\psi\|^2_{\Hk
  {\alpha+\demi}}$.) For example, to control $E_3$ we take any $x_3$
satisfying $2x_2\leq x_3\leq x_1-2s$ (we will make a more precise
choice of $x_3$ later), and we write
\begin{eqnarray*}
\int_{\sigxx s}E_3 \dx\volu & = & E_{3,1} +E_{3,2}\;,\\
E_{3,1} &\equiv& \int_{\sigxx s \bigcap \{ x\geq x_3\}} E_3 \dx\volu \;,\\
E_{3,2} &\equiv& \int_{\sigxx s \bigcap \{ x\leq x_3\}} E_3
\dx\volu\;.
\end{eqnarray*}
By (\ref{L.8.n}), $E_{3,2}$ can be estimated as follows:
\begin{eqnarray*}
|E_{3,2}|&\leq& \sum_{0\leq \beta\leq k} \int_{\sigxx s \bigcap \{ x\leq x_3\}} \err(x)x^{-2\alpha-2+2\beta_1} |\decal^\beta \psi|^2 \dx\volu\\
&\leq& {(2\alpha+1)\varepsilon \over 10} \|\psi\|^2_{\Hk
{\alpha+\demi}}\;,
\end{eqnarray*}
if $x_3$ is chosen small enough.  Once this choice has been done,
we can clearly estimate $E_{3,1}$ as
$$E_{3,1}\leq C\|\psi\|^2_{\Hak} \;,$$
with some constant which is determined by $x_3$.  The integrals of
the error terms $E_1$ and $E_2$ are estimated in the obvious way,
{\em cf.\/} (\ref{L.5}) and (\ref{L.6.0}):
$$\int_{\sigxx s} (E_1+E_2) \dx\volu \leq C\|\varphi(s)\|^2_{\Hak}\;.$$  To
control the terms $D_1$ and $D_2$ we use the evolution equations
(\ref{L.3}): \begin{eqnarray}
\emoins\nmu\decal^\beta \varphi &=& \decal^\beta  (\emoins\nmu\varphi )+ [\emoins \nmu, \decal^\beta  ]\varphi \nonumber \\
 &=& -\decal^\beta (L\psi+B_{11}\varphi+B_{12}\psi -a) + [\emoins\nmu, \decal^\beta ]\varphi \nonumber \\
&=& -L\decal^\beta  \psi +\decal^\beta  a +E_4^\beta\;,
\label{L.15}
\\
E_4^\beta&=& -[\decal^\beta , L]\psi + [\emoins\nmu, \decal^\beta
]\varphi - \decal^\beta (B_{11}\varphi+B_{12}\psi)\;,\nonumber
\\
\eplus\nmu\decal^\beta \psi & = & L^\dagger \decal^\beta \varphi + \decal^\beta  b +E_5^\beta, \label{L.16} \\
E_5^\beta & = & [\decal^\beta , L^\dagger]\varphi +[\eplus \nmu,
\decal^\beta ]\psi -\decal^\beta
(B_{21}\varphi+B_{22}\psi)\;.\nonumber
\end{eqnarray}
Integrating $D_1+D_2$ over $\sigxx s$, one finds that the terms
containing $L\decal^\beta \psi$ and $-L^\dagger\decal^\beta
\varphi$ in (\ref{L.15}) and (\ref{L.16}) cancel out; the terms
containing $\decal^\beta  a$ and $\decal^\beta  b$ are estimated
as (here the somewhat arbitrarily chosen factor $10$ can be
replaced by any other larger number if desired)
\begin{eqnarray*}
\lefteqn{2\sum_{0\leq |\beta|\leq k}\displaystyle\int_{\sigxx s}
\xadu \left(\langle \decal^\beta  \varphi, \decal^\beta  a\rspone
+\langle \decal^\beta  \psi,\decal^\beta  b\rsptwo \right)\dx\volu
} &&
\\
&&\leq \|\varphi\|^2_{\Hak} +\|a\|^2_{\Hak} +
{(2\alpha+1)\varepsilon
  \over 10} \|\psi\|^2_{\Hk {\alpha+\demi}}
+ {10 \over (2\alpha+1)\varepsilon} \|b\|^2_{\Hk
{\alpha-\demi}}\;.
\end{eqnarray*}
 The  terms containing the commutators $[\decal^\beta , L]\psi$ and $[\decal^\beta ,
 L^\dagger]\varphi$,
can be estimated\footnote{This step requires  weighted $L^\infty$
control of $\phi$ and $\psi$, and weighted $W^{1,\infty}$ control
of the coefficients in the equation. The hypothesis $k>n/2+1$ is
not needed if such {\em a priori\/} bounds are known.} using the
weighted commutator inequality (\ref{Mo2}), while the $B_{11}$,
$B_{12}$, \emph{etc.}, terms can be estimated using (\ref{Mo1}),
by an expression of the form
\begin{equation}
  \label{eq:gest}
  CC_1\left(\|\psi\|^2_{{\HH }^\alpha_k} + \|\varphi\|^2_{{\HH }^\alpha_k }+
{(2\alpha+1)\varepsilon \over 10} \|\psi\|^2_{\Hk
  {\alpha+\demi}}\right)\;.
\end{equation}   To estimate the commutator terms
arising from $E_{\pm }^\mu $, we calculate, {\em e.g.}
\begin{eqnarray*}
x^k [E_{\pm}^\mu \pmu, \partial_x^k]\chi &=& \sum_{i=1}^k ({}^i_k)
x^i
(\partial_ x ^i E_{\pm}^\mu)x^{k-i} \px^{k-i}\pmu \chi \\
 &=& E_6 + E_7\;,\\
E_6 &=& \sum_{i=1, \mu \neq x}^k  ({}^i _k )x^i (\px^i
E_{\pm}^\mu) x^{k-i}\px^{k-i} \pmu \chi\;. \end{eqnarray*} The
terms arising from $E_6$ are estimated in a straightforward way as
in \eq{eq:gest} using \eq{Mo1.1}.  The dangerous term $E_7$ can be
written as
\begin{eqnarray*}
 E_7 \equiv \sum_{i=1}^k  ({}^i _k )x^i (\px^i E_{\pm}^x)
x^{k-i} \partial_x^{k-i+1} \chi= \sum_{i=1}^k  ({}^i _k )x^{i-1}
(\px^{i-1} \px E_{\pm}^x) x^{k-i+1} \partial_x^{k-i+1} \chi\;,
\end{eqnarray*} and can thus again be estimated as in \eq{eq:gest}
provided that $\partial_x E^x_-\in{\cG}^0_{k-1}$, that
$\partial_xE^x_+\in{\cG}^{-1}_{k-1}$, and that \eq{L.8.n}  holds.
Other terms in the $E^\mu_\pm$ commutators are handled in a
similar way.

Summarizing, we have derived
\begin{equation}
  \label{eq:gest2}
  \left| \int_{\Sigma_{x_2,x_1,s}} \nmu X^\mu d^{n}\mu\right| \leq
  CC_1\left(\|a(s)\|^2_{{\HH }^\alpha_k} +
\|b(s)\|^2_{{\HH }^{\alpha-\demi}_k} +\|\psi(s)\|^2_{{\HH
}^\alpha_k} + \|\varphi(s)\|^2_{{\HH }^\alpha_k }\right)\;,
\end{equation}
where $d^n\mu$ stands for $dx d\nu$, or for any measure uniformly
equivalent to $dx d\nu$. Stokes theorem,
\begin{eqnarray*}
\int_{\Omega_{x_2,x_1,t}} \nmu X^\mu \,d^{n}\mu \,d\tau&=&
\int_{\partial\Omega_{x_1,x_2,t}} X^\mu dS_\mu
\;, \label{L.18}
\end{eqnarray*}
and our hypotheses on the geometry of the problem lead to
\begin{eqnarray*}
  \|f(t)\|^2_{\Hak} \leq C\left( \|f(0)\|^2_{\Hak} + C_1\int_0^t
    \left(\|a(s)\|^2_{{\HH }^\alpha_k} +
\|b(s)\|^2_{{\HH }^{\alpha-1/2}_k} +\|f(s)\|^2_{{\HH }^\alpha_k}
\right)ds \right)\;.
\end{eqnarray*}
Gronwall's lemma establishes (\ref{L.10}) on the family of
hypersurfaces $\displaystyle \sigxx t$ for $f(t)\in H_{k+1}^{loc}
$. If $f(t)\in H_{k}^{loc} $, we  approximate $f(0)$ by a sequence
of functions $f_n(0)$, with $f_n(0)\in H_{k+1}^{loc} $ converging
to $f(0)$ in ${\HH }^\alpha_{k}(\sigxx t)$, and we solve
Equation~\eq{L.1} with initial data $f_n(0)$. The inequality
(\ref{L.10}) applied to the functions $f_n(t)-f_m(t)$ shows that
$f_n(t)$ is Cauchy in $\Hak$; passing to the limit $n\to\infty$
the desired result for $f$'s such that $f(0)\in {\HH
}^\alpha_{k}(\sigxx t) $ easily follows.

Since all the constants above are $x_2$ independent,  an
elementary argument using the monotone convergence theorem shows
that the inequality (\ref{L.10}) for the $\displaystyle \sigx t$'s
follows from the one for the $\displaystyle \sigxx t$'s by passing
to the limit $x_2\to0$.  \qed

\subsection{Estimates on the space derivatives of the solutions, $\alpha=-1/2$}
\label{Sonehalf}  When $\alpha=-1/2$ we do not have the
$\beta_1=0$ negative terms in $N_1$ at our disposal in \Eq{L.13},
so that we cannot allow coefficients as singular as in the
previous section. To handle that case we keep all the structure
and regularity conditions already made, with the following
supplementary restrictions: \Eq{L.8.n} is replaced by
\begin{equation}
\| E^A_ + (\tau)\|_{{\cG}^{0}_k(M_{x})}+\|  E^x_ +
(\tau)\|_{{\cG}^{1}_{k}(M_{x})}+ \|\nabla_\mu E^\mu_+(\tau)
\|_{L^\infty(M_{x})}\leq  C_1\;.\label{L.8.n1}\end{equation}
Instead of \eq{L.6.2} we require that
\begin{equation}
  \label{L.6.2n}
  \|B_\pm(\tau)\|_{\gok\mxtau }  \leq C_1\;,
\end{equation}
while condition \eq{L.6} becomes
\begin{eqnarray}&\|B_{ab}(\tau)\|_{{\cG}^0_k\mxtau }
 \leq
C_1\;. & \label{L.6n}\end{eqnarray}  We then obtain:
\begin{Proposition} \label{PL.1n}Suppose that
 $k>{n\over 2}+1$, $ k\in \N$, and set either
  $f(t)=f|_{\sigx t}$, $0<x_1\leq x_0$, $0\leq t \leq t_{\max}\equiv
  x_1/2$, or $f(t)=f|_{\sigxx t}$, $0<2x_2< x_1\leq x_0$, $0\leq t <
  t_{\max} \equiv x_1-2x_2$.  Under
  the  hypotheses $\mcC1)$-$\mcC5)$ together with \eq{L.8.n1}-\eq{L.6n}
  there exists a constant $C_2$ depending upon $x_1$, $C_1$, $n$, $N$,
  $k$ and the boundary manifold $\partial M$, such that for all
$f$ satisfying $f(0)\in
H_{k}^{loc} $ and for all $
0<t\leq t_{\max}$ we have
\begin{eqnarray}
\|f(t)\|^2_{\Haok\mxt} & \leq &
C_2e^{C_2t}\left(\|f(0)\|^2_{\Haok(M_{x_1})} +
  \int_0^t
  e^{C_2(t-s)}
  \Big(\|a(s)\|^2_{\Haok\mxs}\right.
\nonumber \\ & &  
+\left.\|b(s)\|^2_{\Hk {-1/2}\mxs}\Big)
  ds\right)\;,
 \label{L.10n}\end{eqnarray}
with an identical estimate with $M_{x_1-*}$ replaced by
$M_{x_2,x_1-*}$.
\end{Proposition}

\proof The proof is essentially identical, but simpler, to that of
Proposition~\ref{PL.1}. We simply note that the key inequality
\eq{eq:gest2} gets replaced by
\begin{equation}
  \label{eq:gest2n}
  \left| \int_{\Sigma_{x_2,x_1,s}} \nmu X^\mu d^{n}\mu\right| \leq
  CC_1\left(\|a(s)\|^2_{{\HH }^{-1/2}_k} +
\|b(s)\|^2_{{\HH }^{{-1/2}}_k} +\|\psi(s)\|^2_{{\HH }^{-1/2}_k} +
\|\varphi(s)\|^2_{{\HH }^{-1/2}_k }\right)\;.
\end{equation}\qed

\subsection{Estimates on the time derivatives of the solutions}
\label{ssltest}

The hypotheses assumed in the previous section ensure that we can
algebraically solve Equation \eq{L.1} for $\partial_\tau f$, and
then recursively obtain formulae for $\partial^i_\tau f$. Under
the hypotheses of Proposition~\ref{PL.1}, it is then
straightforward to obtain estimates on the norms
\newcommand{\hki}{{\HH}^\alpha_{k-i}(\Sigma_{x_1-2\tau})}
\newcommand{\hkiz}{{\HH}^\alpha_{k-i}(\Sigma_{x_1})}
$$ \|((x\partial_\tau)^i f)(\tau)\|_{\hki}\;, \qquad 0
\leq i\leq k\;,$$ provided suitable weighted conditions are
imposed on the $\tau$ derivatives of the coefficients of
Equation~\eq{L.1}. However, we would like to obtain derivative
estimates without the $x$ factors, uniformly in $\tau$. Clearly a
necessary condition for the existence of such  estimates is that
\begin{equation}
  \label{eq:t1}
  \|(\partial_\tau^i f)(0)\|_{\hkiz}< \infty\;, \qquad 0
\leq i\leq k\;.
\end{equation}
It turns out that \eq{eq:t1} does not need to hold  for arbitrary
initial
 data $f(0)\in \HH ^\alpha_k$, and the requirement that it does lead to the
 \emph{$j$-th order compatibility
   conditions}:
 by definition, these are the conditions on $f(0)$ which ensure that
 Equation~\eq{eq:t1} holds for $0\leq i \leq j$. Since, for sufficiently differentiable
 solutions of Equation~\eq{L.1}, all the
 derivatives $\partial^i_\tau f(0)$ can be explicitly written as an
 $i$-th order differential operator acting on $f(0)$, the compatibility
 conditions are conditions on the behaviour of the initial data $f(0)$
 near the ``corner'' $x=0$; we shall therefore sometimes refer to them
 as ``{\em corner conditions}''. We note that there can be corner conditions
 in weighted Sobolev spaces, or in weighted H{\"o}lder spaces; in this
 section we will be mainly interested in the latter, defined by
 Equation~\eq{s2} below.

The following example is instructive in this context: For $0\leq t
< y$ let $g$ be a solution of the $1+1$ dimensional wave equation
\begin{equation}
  \label{eq:2dwe}
  \left({\partial^2 \over \partial t ^2 } - {\partial^2 \over \partial
      y ^2 }\right) g = 0\;,
\end{equation}
with initial condition
$$g\Big|_{t=0} = 2C y^{\alpha +1}\;, \qquad {\partial  g\over \partial t
  }\Big|_{t=0} = 2(\alpha+1) y^{\alpha } \;,$$
for some constants $C,\alpha \in \R$. From Equation~\eq{eq:2dwe}
we can obtain a system of the form \eq{L.3} by introducing $\tau
=t$, $x=y-t$, $\varphi= (g, (\partial_\tau-2\partial_x)g)$, $\psi=
\partial_\tau g$, and setting $L=0$, $E^\mu_-\partial_\mu =
\partial_\tau \otimes \mathrm{id}$, $E^\mu_+\partial_\mu =
(\partial_\tau -2\partial_x)$, so that we have
\begin{eqnarray*}
  \partial_\tau \left(\begin{array}{c} g \\
      (\partial_\tau-2\partial_x)g \end{array} \right)  -
      \left(\begin{array}{c} \psi \\
      0 \end{array} \right) & = &
  \left(\begin{array}{c}
0 \\ 0
  \end{array}\right)\;,
\\
(\partial_\tau-2\partial_x) \psi  & = & 0\;.
\end{eqnarray*}
The solution is
\begin{eqnarray*}
 g & =&  (C+1) (y+t)^{\alpha+1} +
(C-1)(y-t)^{\alpha+1}
\\  & = & (C+1) (2\tau+x)^{\alpha+1} +
(C-1) x^{\alpha+1}\;.
\end{eqnarray*}
It follows that for each $0\leq \tau\leq 1$, $k\in\N$, and
$\beta<\min\{0,\alpha+1\}$, we have $g(\tau,\cdot)\in
\HH^\beta_k((0,10]) $, consistently with Proposition~\ref{PL.1}.
Somewhat surprisingly, for $\tau>0$ and for all $i\in\N$ the
functions $\partial^i_\tau g(\tau,\cdot)$ are smooth in $x$ up to
$x=0$. However, the $L^\infty$ bound for $\partial^i_\tau
g(\tau,\cdot)$ blows up as $\tau$ tends to zero except in the case
\begin{equation}
  \label{eq:2dwe1}
  C= -1\;.
\end{equation}
Condition \eq{eq:2dwe1} is precisely the corner condition needed
for $\partial_\tau g(0,\cdot)$ to be better behaved than
$\partial_x g(0,\cdot)$ at $x=0$. In the example under
consideration the fulfillment of the first order corner condition
guarantees already that all the $\tau$ derivatives of $g$ will be
well behaved, but we do not expect this to be true in general.

Let us pass to a derivation of the desired estimates. We shall use
a method which works directly in weighted H\"older spaces,
avoiding the use of weighted Sobolev spaces; the price one pays is
the need to consider systems somewhat less general than \eq{L.3},
but still general enough for our purposes. More precisely, in this
section we restrict our attention to systems of the form
\begin{deqarr}
  \partial_\tau\varphi + B_{11}\varphi + B_{12}\psi&=& L_{11}\varphi +
  L_{12}\psi + a
  \label{s1a}\;,\\ e_+\psi+
  B_{21}\varphi + B_{22}\psi &=&   L_{21}\varphi +
  L_{22}\psi + b \label{s1b}\;,\arrlabel{s1}
\end{deqarr}
with
$$e_+\psi\equiv (\partial_\tau-2\partial_x)\psi\;.$$
We assume that the $L_{ab}$'s, $a,b=1,2$ are first order
differential operators of the form
\begin{equation}
  \label{s2}
  L_{ab} = L_{ab}^A\partial_A + x L_{ab}^\tau \partial _\tau + x L_{ab}^x
  \partial _x\;,
\end{equation}
with bounded coefficients $L_{ab}^\mu$; no symmetry hypotheses are
made. Clearly the intersection of systems of equations satisfying
\eq{s1} with those of the form \eq{L.3} is non-empty. (As we will
see in Sections~\ref{sslwe} and \ref{Swave} below, non-linear wave
equations on Minkowski space-time can be written in the form
\eq{s1}.) In particular Proposition~\ref{PL.1} provides a large
class of solutions of \eq{s1} such that
$$(\varphi,\psi)(\tau)\in \Hak\mxtau\subset{\mcC}_{\ell} ^  \alpha
\mxtau $$ for $\ell<k-n/2$. We shall therefore assume that a
solution $f=(\varphi,\psi)$ satisfying $f(\tau)\in {\mcC}_{\ell} ^
\alpha \mxtau$ is given, and study its $\tau$-differentiability
properties. For the purposes of the proof below it is convenient
to introduce auxiliary spaces ${\mcC}_{\ell|p} ^ \alpha
   (\Omega)$ defined, for $p\leq \ell$, as the space of functions $f$
   in $C_\ell (\Omega)$ such that the norm
$$\|f\|_{{\mcC}_{\ell|p} ^  \alpha
   (\Omega)}\equiv \sup_{\Omega}\sum _{\begin{array}{c} 0\leq
     i+j+k+|\gamma|\leq \ell
\cr 0\leq k \leq p\end{array}} x^{-\alpha}|(x\partial_x)^i
(x\partial_\tau)^j \decal _v^\gamma \partial_\tau^k f| $$ is
finite. Obviously, ${\mcC}_{\ell|\ell} ^ \alpha={\mcC}_{\ell} ^
\alpha$. Similarly one defines ${\mcC}_{\ell|p} ^ {\alpha,\beta}
   (\Omega)$ using the norm $$\|f\|_{{\mcC}_{\ell|p} ^  {\alpha,\beta}
   (\Omega)}\equiv \sup_{\Omega}\sum _{\begin{array}{c} 0\leq
     i+j+k+|\gamma|\leq \ell
\cr 0\leq k \leq p\end{array}} (1+|\ln
x|)^{-\beta}x^{-\alpha}|(x\partial_x)^i (x\partial_\tau)^j \decal
_v^\gamma \partial_\tau^k f| \;.$$ Clearly ${\mcC}_{\ell|p} ^
{\alpha}   (\Omega)={\mcC}_{\ell|p} ^ {\alpha,0}    (\Omega)$. We
shall write ${\mycal C}^{\alpha,\beta}_\ell$ for ${\mycal
C}^{\alpha,\beta}_{\ell|\ell}$.
\begin{Proposition}
  \label{Ptd} Let $ \alpha\leq 0$, $\ell \in \N$, write $\Omega$ for
  $\OmxT$ (with $\OmxT$ as in \eq{L.7}), and suppose that $L_{ab}^\mu, B_{ab}\in
  \mcC^0_{\ell}(\Omega)$, $a\in \mcC^\alpha_{\ell-1}(\Omega)$, $b\in
  \mcC^{\alpha-1}_{\ell-1}(\Omega)$. Consider $f\equiv (\varphi,\psi)$
  --- a solution of \eq{s1} satisfying $$\forall
  \tau\in[0,T]\qquad f(\tau)\in {\mcC}_{\ell} ^ \alpha \mxtau\;.$$  Then:
\begin{enumerate}
\item For all $\epsilon>0$ we have
  $$ (\varphi,\psi)\in {\mcC}_{\lfloor\ell/2\rfloor} ^ {\alpha,\beta}
  \left(\Omega\cap\{x+2\tau > \epsilon\}\right)\;.$$ Further,
  for any $\tau>0$ the compatibility conditions of order
  $p=\lfloor\ell/2\rfloor$ (the integer part of $\ell/2$) are satisfied by
  $(\varphi(\tau),\psi(\tau))$:
  \begin{equation}
    \label{s2.}
    \forall \;1\leq i\leq p  \qquad \partial^i_\tau \varphi(\tau),
    \partial^i_\tau \psi(\tau) \in {\mcC}_{\ell-i} ^ {\alpha,\beta}
\mxtauz\;,
  \end{equation}
 Here $\beta=\lfloor\ell/2\rfloor$ if
  $\alpha =0$, and $\beta=0$ otherwise.
\item If there exists $1\leq p\leq \ell/2$, $p\in\N$, such that \Eq{s2.}
holds with $\beta=0$ at $\tau=0$, then
\begin{equation}
    \label{s3}
    (\varphi, \psi) \in{\mcC}_{\ell-p|p} ^{\alpha,\beta}
   (\Omega) \subset {\mcC}_{p} ^{\alpha,\beta}    (\Omega)\;,
\end{equation}
with $\beta =p $  if
  $\alpha =0$, and $\beta=0$ otherwise.
 \end{enumerate}
\end{Proposition}
\remark The method of proof here gives a  number of well
controlled time derivatives smaller by a factor 2 than the number
of space ones. This is, however, irrelevant, when $\ell=\infty$,
which is the main point of interest in this work. We note that
energy estimates as in the proof of Theorem~\ref{T2t} below
provide an alternative, more complicated way of establishing a
stronger statement, with more controlled time derivatives for
large $\ell$'s. In the linear case considered here the function
$F$ occuring there vanishes, so that all the complications arising
from the non-linearities disappear, and somewhat stronger results
can be obtained using the methods there.

\medskip

\proof By rearranging terms and redefining the $L_{ab}$'s, the
$B_{ab}$'s, and the source functions $a$ and $b$  we may without
loss of generality assume that
$$L^\tau_{ab}\equiv 0\;.$$
One can  rewrite Equations~\eq{s1} as $x\partial_\tau
(\varphi,\psi)=$ a partial differential operator linear in
$x\partial_x$ and $\partial_v$; by iteration this immediately
yields $(\varphi,\psi)\in \mcC^\alpha_{\ell|0}$. \Eq{s1a} shows
then that $\partial_\tau\varphi\in \mcC^\alpha_{\ell-1|0}$, hence
$\varphi\in \mcC^\alpha_{\ell|1}$. On the other hand, \Eq{s1b}
gives
 $e_+(\psi)\in \mcC^\alpha_{\ell-1|0}+ \mcC^{\alpha-1}_{\ell-1}$,
hence   $\partial_\tau e_+(\psi)\in \mcC^{\alpha-1}_{\ell-2|0}$.
In order to extract further information out of \Eq{s1b} we use the
straightforward identity
\begin{equation}
  \label{s6}
 \psi(x,v^A,\tau)=\psi(x+2\tau,v^A,0) + \int_{x/2}^{\tau+x/2}
e_+(\psi)(2v,v^A,\tau-v+x/2) \dv\;.
\end{equation}
(We note that for each $\epsilon>0$ the first term above is
uniformly $C_\ell$ on the set $\Omega\cap\{x+2\tau\geq
\epsilon\}\cap \{x\leq x_0\}$.) Differentiating \Eq{s6} one
obtains
$$ \partial_\tau\psi(x,v^A,\tau)=\partial_\tau\psi(x+2\tau,v^A,0) +
 \int_{x/2}^{\tau+x/2} \partial_\tau e_+(\psi)(2v,v^A,\tau-v+x/2) \dv\;;
$$
since  $\alpha\leq0$ and $\partial_\tau e_+(\psi)\in
\mcC^{\alpha-1}_{\ell-2|0}$, straightforward estimations show that
$\partial_\tau \psi\in\mcC^\alpha_{\ell-2|0}$, hence $
\psi\in\mcC^\alpha_{\ell-1|1}$ if $\alpha\neq 0$, while
$\psi\in\mcC^{0,1}_{\ell-1|1}$ when $\alpha= 0$.

Let $\beta_r=0$ if $\alpha\neq 0$ and $\beta_r=r$ when $\alpha=0$,
and suppose that $\varphi\in\mcC^{\alpha,\beta_r}_{\ell+1-r|r}$
and $ \psi\in\mcC^{\alpha,\beta_r}_{\ell-r|r} $ for some $1\leq r
\leq (\ell-1)/2$; we have already shown this to hold for $r=1$.
Equation~\eq{s1a} gives $$\partial_\tau\varphi\in
\mcC^{\alpha,\beta_r}_{\ell-r-1|r}\quad \Longrightarrow \quad
\varphi\in\mcC^{\alpha,\beta_r}_{\ell-r|r+1}\;.$$ It then follows
from \Eq{s1b} that
$$e_+(\psi)\in \mcC^{\alpha,\beta_r}_{\ell-r-1|r}\quad \Longrightarrow \quad
\partial_\tau^{r+1}e_+(\psi)\in\mcC^{\alpha-1,\beta_r}_{\ell-2r-2|0}\;.$$
Differentiating $r+1$ times \Eq{s6} with respect to $\tau$ we
obtain
$$ \partial_\tau^{r+1}\psi(x,v^A,\tau)=\partial^{r+1}_\tau\psi(x+2\tau,v^A,0) +
 \int_{x/2}^{\tau+x/2} \partial_\tau^{r+1} e_+(\psi)(2v,v^A,\tau-v+x/2) \dv\;,
$$
which gives $\partial_\tau^{r+1}
\psi\in\mcC^{\alpha,\beta_r}_{\ell-2r-2|0}$, hence $\psi\in
\mcC^{\alpha,\beta_r}_{\ell-r-1|r+1}$, and the induction is
completed. \qed

\subsection{Polyhomogeneous  solutions }
\label{sslphg}

We now wish to show that solutions with polyhomogeneous initial
data will be polyhomogeneous.
 Let $\Omega_{x_0,T}$ be defined by \Eq{L.7}; we shall denote by
$\mcA_k^{\delta}(\Omega_{x_0,T})$ the space of functions $f$
defined on $\Omega_{x_0,T}$ which can be written in the form
$$\sum_{i=0}^k \sum_{j=0}^{N_i}
x^{i\delta}\ln^j x \, f_{ij} + f_{\alpha+k\delta+\epsilon}\;,$$
for some $\epsilon>0$, some functions $f_{ij}\in
C_{\infty}(\overline{\Omega_{x_0,T}})$, and some   sequence
$(N_i)$ of non-negative integers. We also require that
$f_{\alpha+k\delta+\epsilon}\in
\mcC_{\infty}^{\alpha+k\delta+\epsilon}(\Omega_{x_0,T})$. We set
$$ \mcA_\infty^{\delta}:=\cap_{k\in\N}\mcA_k^{\delta}\;.$$
The following properties are useful in what follows:
\begin{itemize}
\item If $0<x_1<x_0 -T/2$, then a function $f\in
C_\infty(\Omega_{x_0,T})$ is in $\mcA_k^\delta(\Omega_{x_0T})$ if
and only if for any coordinate patch $\mcO$ of $\pM$ we have $f\in
\mcA_k^\delta({\mcU}_{x_1})$, where ${\mcU}_{x_1} = ]0,x_1[\times
\mcO \times [0,T]$, and if $f\in
C_\infty(\overline{\Omega_{\textrm{\scriptsize int}}})$, where
$\Omega_{\textrm{\scriptsize int}} = \Omega_{x_0,T} \cap \{x\geq
x_1\} $.

\item For all $\epsilon >0$ we have $ \mcC_\infty^{\beta+p
\delta+\epsilon}\subset  x^{\beta}\mcA_p^\delta $; in particular $
\mcC_\infty^\epsilon \subset \mcA_0^\delta$;
\item
It  does not hold that $\mcA_k^\delta \subset \mcC_\infty^0$,
however, for all $\epsilon
>0$ we have  $ \mcA_k^\delta \subset
\mcC_\infty^{-\epsilon}$. More precisely, if  $f\in
\mcA_k^{\delta}$,
 then there exists $p\in\N$ such that  $(1+|\ln x|^2)^{-p/2}f \in
\mcC_\infty^0$.
\item As before we assume that $1/\delta \in \N$, which implies $x
\mcA_k^\delta \subset \mcA_{k+1/\delta}^\delta\subset
\mcA_{k+1}^\delta\;. $
\item $\mcA_k^{\delta}$ is stable under multiplication: if $f,g\in
\mcA_k^{\delta} $, then $fg \in \mcA_k^{\delta} $.
\item $\mcA_k^{\delta}$ is stable under differentiation with respect
to $\tau$ and to $v$, as well as under $x\partial_x$: if $f\in
\mcA_k^{\delta} $, then $\partial_\tau f \;, X_i \cdot f $ ($i\geq
2$), $x\partial_xf \in \mcA_k^{\delta} $, with the vector fields
$X_i$ defined in Section~\ref{S2}, {\em cf.\/}~\Eq{champscoord}.

\end{itemize}

In this section we will consider  systems of the  form
\begin{deqarr}
  \partial_\tau\varphi + B_{11}\varphi + B_{12}\psi&=& L_{11}\varphi +
  L_{12}\psi + a
  \label{equa1}\;,\\ \partial_x\psi+ B_{21}\varphi + B_{22}\psi &=&
  L_{21}\varphi +
  L_{22}\psi + b \label{equa2}\;,\arrlabel{equa}
\end{deqarr}
with the $L_{ij}$'s, $i,j=1,2$ of the form
 \be L_{ij}=
L_{ij}^{A}\partial_A + L_{ij}^\tau
\partial_\tau + xL_{ij}^x\partial_x\,,
\label{equa1+}\ee with
 \be
  L_{11}^{\mu}\in x^\delta\mcA_{k-1}^{\delta}\;, \quad L_{21}^{\mu}\;,
  L_{12}^{\mu}\;, L_{22}^{\mu}\in \mcA_k^\delta \;.
\label{HLdelta} \ee No symmetry hypotheses  are made on the
matrices $L^\mu_{ij}$. Conditions~\eq{equa1}-\eq{HLdelta} are
easily seen to be compatible with those made elsewhere in this
paper, {\em cf.,
  e.g.,\/} the proof of Corollary~\ref{Cpolylineaire} below. The reader
is warned, however, that the operators $L_{ij}$ here do {\em
not\/} coincide with those in \eq{s1}: to bring \eq{s1} into the
form \eq{equa} one needs to multiply \Eq{s1b} by $-1/2$, transfer
the operator $\partial_\tau$ from the left- to the right-hand-side
of \eq{s1}, and appropriately redefine the $L_{2j}$'s.

We start with the following result, which assumes that the
solutions have both space and time derivatives controlled, in the
sense of  weighted Sobolev spaces; recall that this hypothesis can
be justified for equations satisfying moreover the hypotheses of
the previous sections:
\begin{Theorem}\label{Tlemme1}
  Let $\beta,\, \beta'\in \R$,
$k\in \N\cup\{\infty\}$, 
and let $(\varphi,\psi)$ be  a solution of
  (\ref{equa}) in $
  \mcC_{\infty}^{\beta'}(\Omega_{x_0,T}) $. Suppose that \eq{HLdelta}
  holds, and that
  \begin{deqarr}&
    B_{11}\in \left(\mcA_k^\delta \cap L^\infty\right)
    (\Omega_{x_0,T})\;,\qquad
B_{12}, B_{22}, B_{21}\in \mcA_k^\delta (\Omega_{x_0,T}) \;,&
\label{H3first}
\\ &a,b \in x^\beta\mcA_k^\delta (\Omega_{x_0,T})\;,
\qquad\varphi(0)\in x^\beta\mcA_k^\delta(M_{x_0})\;.&
\label{H4first}
  \end{deqarr}
 Then $$\varphi\in
\left(x^\beta\mcA_k^\delta+\mcA_k^\delta\right)(\Omega_{x_0,T})\;,\quad
\psi\in  \left(x^{\beta+1}\mcA_{k}^\delta
+x\mcA_{k}^\delta\right)(\Omega_{x_0,T})+
C_\infty(\overline{\Omega_{x_0,T}})\;.$$ If one further assumes
$$L_{12}^\mu,B_{12},a,\varphi(0)\in L^\infty(\Omega_{x_0,T})\;,$$
then it also holds that
$$\varphi\in
\left(x^\beta\mcA_k^\delta+\mcA_k^\delta\cap
L^\infty\right)(\Omega_{x_0,T})\;.$$
\end{Theorem}
\proof It is convenient to decompose $B_{11}$ in the obvious way
as
 $$B_{11} = B_{11}^0+ B_{11}^\delta\;,$$ with $B_{11}^\delta\in
x^\delta\mcA_{k-1}^\delta$ and $B_{11}^0\in C_\infty$. We rewrite
(\ref{equa}) as
\begin{deqarr}
\partial_\tau \varphi + B_{11}^0 \varphi & = & c_1\;, \label{equa1.1}
\\
\partial_x \psi  & = & c_2\;, \label{equa1.2}
\end{deqarr}
where
\begin{deqarr}
c_1&:=& L_{11}\varphi + L_{12}\psi + a -B_{12}\psi -
B_{11}^\delta\varphi \;, \label{equa2.1}
\\c_2 &:=&
L_{21}\varphi + L_{22}\psi + b-B_{21}\varphi -B_{22}\psi\;,
\label{equa2.2}
\end{deqarr}
In what follows we let $\epsilon>0$ be a positive constant, which
can be made as small as desired, and which may change from line to
line. We note that $c_2$ is in $ \mcC_\infty
^{\beta'-\epsilon}+x^\beta\mcA_k^\delta$, and integration in $x$
of \eq{equa1.2}, together with Propositions~\ref{propositionx1}
and \ref{propositionx}, gives
 $$\psi=\psi_{0}+\psi_{\beta'+1-\epsilon}+\psi_{\mbox{\rm\scriptsize phg}}\;,$$
where
$$\psi_{0}(\cdot) =\cases{
         \lim_{x\to 0}\psi(x,\cdot)\;, & {if} $
                        \beta'+1-\epsilon >0$,
\cr  0 \;,& otherwise,}\qquad $$ with $$\psi_{0} \in
C_\infty(\overline{\Omega_{x_0,T}})\;,\quad
\psi_{\beta'+1-\epsilon}\in
\mcC_\infty^{\beta'+1-\epsilon}(\Omega_{x_0,T})\;,\quad\psi_{\mbox{\rm\scriptsize
phg}}\in x^{\beta+1}\mcA_k^\delta(\Omega_{x_0,T})\;,$$ hence
$$\psi\in C_\infty +
\mcC_\infty^{\beta'+1-\epsilon}+x^{\beta+1}\mcA_k^\delta\;.$$
Since $L_{11}\varphi \in \mcC_\infty ^{\beta'+\delta-\epsilon}$
($\partial_x\varphi \in \mcC_{\infty}^{\beta'-1}$ and
$xL_{11}^x\in x\mcA_k^\delta\cap \mcC^\delta_0\subset
\mcC^\delta_\infty $; similarly for the other derivatives), we
find that $$c_1 \in \mcA_k^\delta+x^\beta \mcA_k^\delta+
\mcC_\infty ^{\beta'+\delta-\epsilon}\;.$$ We can then apply
Proposition~\ref{proprietetau} to (\ref{equa1.1}) to conclude that
\be\label{indp}\varphi \in \mcA_k^\delta+x^\beta
\mcA_k^\delta+\mcC_\infty ^{\beta'+p\delta-\epsilon}\;,\ee with
$p=1$. Coming back to $c_2$ we find now that  $c_2\in
\mcA_k^\delta+x^\beta \mcA_k^\delta+\mcC_\infty
^{\beta'+p\delta-\epsilon}$, and by
 Proposition
\ref{propositionx} we obtain  \be\label{indp+1}\psi \in C_\infty+
x\mcA_k^\delta+x^{\beta+1}\mcA_k^\delta+
\mcC_\infty^{\beta'+p\delta+1-\epsilon}\;,\ee still with $p=1$. To
conclude, we proceed by induction; let  $\beta'+p\delta\leq
\beta+k$ and suppose that \Eqs{indp}{indp+1} hold; it follows that
$c_1\in \mcA_k^\delta+ x^{\beta}\mcA_k^\delta+
\mcC_\infty^{\beta'+(p+1)\delta-\epsilon}\;.$  Applying
Proposition~\ref{proprietetau} to (\ref{equa1.1}) gives \eq{indp}
with $p$ replaced by $p+1$. It follows that $c_2\in \mcA_k^\delta+
x^{\beta}\mcA_k^\delta+
\mcC_\infty^{\beta'+(p+1)\delta-\epsilon}$;
Proposition~\ref{propositionx} applied to (\ref{equa2}) gives
\eq{indp+1} with $p$ replaced by ${p+1}$, and the result is
established.  \qed

As a straightforward consequence of Theorem~\ref{Tlemme1} we
obtain:

\begin{Corollary}\label{Cpolylineaire} Let $\beta'\in\R$, let
$(\varphi,\psi)\in \mcC_\infty^{\beta'}(\Omega_{x_0,T})$ be a
solution of the system (\ref{L.3}), and suppose that
\begin{deqarr}&B_{ij}, E_{\pm}^\mu, B_\pm ,\ell,
  \ell^\dagger,\ell^A,(\ell^A)^\dagger\in
  \mcA_k^\delta(\Omega_{x_0,T})\;, &
\label{HP1} \\&
 E_-^\tau \ \mbox{and }\
  E_+^x \ \mbox{ --- invertible,  with  }\  (E_-^\tau)^{-1},
(E_+^x)^{-1}\in   \mcA_k^\delta(\Omega_{x_0,T}) \;,
&\label{HP2}\\&
  {(E_-^\tau)^{-1}E_-^x}   \in
  x\left(\mcA_k^\delta\cap\mcC^\delta_0\right)(\Omega_{x_0,T})\;,\quad
  (E_-^\tau)^{-1}E_-^A\in
  x^\delta\mcA_{k-1}^\delta
  (\Omega_{x_0,T})
\;, \nonumber & \\ && \label{HP3}\\ &(E_-^\tau)^{-1}(B_{11}+
B_-)\in L^\infty(\Omega_{x_0,T})\;.& \label{HP4}\arrlabel{HP}
\end{deqarr}
  If $$a,b \in x^\beta\mcA_k^\delta(\Omega_{x_0,T})\;,\quad \varphi(0)\in
  x^\beta\mcA_k^\delta(M_{x_0})\;,$$
  with $\beta\in \R$,
  then $$\varphi\in
\left(x^\beta\mcA_k^\delta+\mcA_k^\delta\right)(\Omega_{x_0,T})\;,\quad
\psi\in
\left(x^{\beta+1}\mcA_{k}^\delta+x\mcA_k^\delta\right)(\Omega_{x_0,T})+
C_\infty(\overline{\Omega_{x_0,T}})\;.$$
  In particular, if $k=\infty$ then the solution is polyhomogeneous.
\end{Corollary}
\proof: We write Equation~(\ref{L.3}) as
  \begin{eqnarray} \partial_\tau\varphi + (E_-^\tau)^{-1}\left\{ (B_{11}+
  B_-)\varphi +\ell \psi\right\}&=& (E_-^\tau)^{-1}
(E_-^i\partial_i\varphi - \ell^A\partial_A\psi + a
  )\label{equa1'}\nn\\ \partial_x\psi-(E_+^x)^{-1}\left\{\ell^\dagger
\varphi- ( B_{22}+B_+)\psi\right\}  &=&
  (E_-^\tau)^{-1}( (\ell^A)^\dagger\partial_A\varphi +
E_+^\tau\partial_\tau \psi +
  E_+^A\partial_A\psi+ b)\;,
\nn\\ &&\label{equa2'} \end{eqnarray} which is of
  the form  (\ref{equa}), and we note that the hypotheses made on
  the coefficients of \Eq{equa2'} imply those of Theorem~\ref{Tlemme1}.
  \proofend.

  An unsatisfactory feature of results such as Theorem~\ref{Tlemme1}
  is that uniform estimates both on space and time derivatives of the
  solutions are assumed.
  Recall that uniform control of time derivatives can be
  obtained only if corner conditions are satisfied, and the hypotheses
  of Theorem~\ref{Tlemme1} require an infinite number of those to be
  fulfilled. The same techniques can be used to obtain various
  expansions of solutions when a finite number of time derivatives
  are controlled only, but the statements turn to be out somewhat less
  elegant. We give an example of such results when $\delta=1$:
\begin{Theorem}\label{Tlemme2}
  Let $\beta\in \R$,
$k\in \N\cup\{\infty\}$, 
and let $(\varphi,\psi)$ be  a solution of
  (\ref{equa}) in $
  \mcC_{\ell}^{\beta}(\Omega_{x_0,T}) $ for some $\ell \geq 1$.
  If \Eqs{HLdelta}{HP}   hold with $\delta=1$, then for any
  $\lambda<1$ we have
  \begin{eqnarray}
    &
\displaystyle\varphi\in \left(x^\beta\mcA_k^1+\mcA_k^1
+\cap_{\ell-2j-2\geq0} \mcC^{\beta+j+\lambda} _{\ell-2j-2}
\right)(\Omega_{x_0,T})\;,& \nn \\& \displaystyle \psi\in
\left(x^{\beta+1}\mcA_{k}^1+x\mcA_{k}^1+\cap_{\ell-2j-1\geq0}
\mcC^{\beta+j+1+\lambda} _{\ell-2j-1}\right)(\Omega_{x_0,T})+
C_\infty(\overline{\Omega_{x_0,T}})\;. & \label{partexphg}
\end{eqnarray} If one further assumes
$$L_{12}^\mu,B_{12},a,\varphi(0)\in L^\infty(\Omega_{x_0,T})\;,$$
then it also holds that
$$\varphi\in
\left(x^\beta\mcA_k^1+\mcA_k^1\cap L^\infty +\cap_{\ell-2j-2\geq0}
\mcC^{\beta+j+\lambda} _{\ell-2j-2}\right)(\Omega_{x_0,T})\;.$$
\end{Theorem}

\proof The result is obtained through a repetition of the proof of
Theorem~\ref{Tlemme1}, keeping track of the differentiability of
the remainder terms. \qed

  We are ready now to prove polyhomogeneity of solutions of the Cauchy
  problem for \Eq{L.3}. We consider only the simplest case of
  equations satisfying the conditions \eq{HP+} below, considerably
  more general statements can be proved using similar methods. The
  differentiability hypotheses below are clearly satisfied by
  equations with smooth bounded coefficients; however, they also allow
  for a wide class of equations with polyhomogeneous coefficients.
We restrict ourselves to the case in which the corner conditions
are
  satisfied to arbitrary order; if not, one obtains expansions as in
\eq{partexphg}, with a remainder in which a finite number only of
time
  derivative are controlled; such results can be proved by identical
  arguments, compare the proof of Theorem~\ref{Tlemme2}. We hope
  to be able to show in a near future that the corner conditions
  are not needed, in which case one should obtain polyhomogeneous
  expansions in which uniformity is lost when the corner
  $\tau=x=0$ is approached; this will be discussed elsewhere.
\begin{Theorem}
  \label{Tphglin}
Consider a solution $(\varphi,\psi)\in C_\infty\times C_\infty$ of
the system (\ref{L.3}), suppose that in addition to \eq{L.5},
\eq{eq:simpl}, \eq{eq:L.5.1},
 and \eq{HP1} we have
\begin{deqarr}&B_{11}, B_-,E_{\pm}^\mu, \ell,
  \ell^\dagger\in L^\infty(\Omega_{x_0,T})\;, &
\label{hp2}\\
&E^\mu_-\Big|_{x=0}=\partial_\tau \otimes \mathrm{id}\;, \qquad
E^\mu_+\Big|_{x=0}=(\partial_\tau-2\partial_x) \otimes
\mathrm{id}\;,\label{hp3} & \\ &E^x_\pm-E^x_\pm\Big|_{x=0}\;,
E^\tau_\pm-E^\tau_\pm\Big|_{x=0} \in
x^{1+\delta}
\mcA_\infty^\delta
(\Omega_{x_0,T}) \;,& \label{hp1} \\ & \quad E_-^A\in
x\mcA_\infty^\delta(\Omega_{x_0,T})\;.& \label{hp4}\arrlabel{HP+}
\end{deqarr}
  If $$a,b \in x^\beta\mcA_k^\delta(\Omega_{x_0,T})\;,\quad \varphi(0)\in
  x^\beta\mcA_k^\delta(M_{x_0})\;,$$
  with $\beta\in \R$, and if the initial data satisfy {\em corner
  conditions to arbitrary order}, in the sense that
\begin{equation}
  \label{eq:cptphg}
  \forall \; i\in \N \qquad \partial_\tau^i \varphi(0),
  \;\partial_\tau^i \psi(0)\in
  \mcC^\lambda_\infty(M_{x_0})\;,
\end{equation}
for some ($i$-independent) $\lambda\in\R$,  then
$$\varphi\in
\left(x^\beta\mcA_k^\delta+\mcA_k^\delta\right)(\Omega_{x_0,T})\;,\quad
\psi\in
\left(x^{\beta+1}\mcA_{k}^\delta+x\mcA_k^\delta\right)(\Omega_{x_0,T})+
C_\infty(\overline{\Omega_{x_0,T}})\;.$$
  In particular, if $k=\infty$ then the solution is polyhomogeneous.
\end{Theorem}

\remark The class of initial data satisfying corner conditions to
arbitrary order is rather large; for example, if an initial data
set $(\varphi(0),\psi(0))$  satisfies them, and if $f,g$ are
arbitrary functions smooth up to boundary on the initial data
hypersurface, then $(\varphi(0)+f,\psi(0)+g)$ will also satisfy
those conditions. More generally, large classes of such initial
data can be constructed using a polyhomogeneous generalization of
the Borel summation lemma.

\medskip

 \proof The hypothesis \eq{eq:cptphg} with $i=0$ and
Proposition~\ref{PL.1} show that for all $\tau \in [0,T ]$ we have
\begin{eqnarray*} \varphi(\tau),
  \; \psi(\tau)\in
  \mcC^\lambda_\infty(M_{x_0/2})\;.
\end{eqnarray*}
Proposition~\ref{Ptd} shows then that the hypotheses of
Corollary~\ref{Cpolylineaire} are satisfied, and the result
follows. \qed

\section{The semi-linear scalar wave equation}
\label{sslwe} Let $f$ be a solution of the semi-linear wave
equation
\be
\Box_\stsg  f = H(x^\mu,f)\;,\label{SE.1} \ee here $\Box_\stsg $
is the d' Alembertian associated with $\stsg $. Set
\be
\tilde{f} = \Omega^{-{(n-1)\over 2}}f \;; \label{SE.2} \ee Letting
$\tilde{\stsg } = \Omega^{2} \stsg $ as in (\ref{C.1}), from
(\ref{C.3}) we obtain \be \Box_{\tilde{\stsg }}\tilde{f} =
{n-1\over 4n} (\tilde{R} - {R\over \Omega^2})\tilde{f} +
\Omega^{-{n+3\over 2}} H(x^\mu, \Omega^{n-1\over 2}\tilde{f})\;.
\label{SE.3} \ee Let $\stsg =\eta$ be the Minkowski metric; under
the conformal transformation (\ref{C.4}) one obtains from
(\ref{C.5}) that $\tilde{\stsg }$ is again the Minkowski metric,
and (\ref{SE.3}) becomes
\be \Box_\eta \tilde{f} = \Omega^{-{n+3\over 2}} H(x^\mu,\Omega^{n-1\over
  2} \tilde{f})\;.
\label{SE.4} \ee We shall assume that  the initial data for $f$
are given on a hypersurface $\Sigma\subset {\mycal M}$, which, in
a neighbourhood ${\mycal O}$ of $\scrip$ is given by the equation
\be \Sigma\cap {\mycal O} = \{ y^0 = \demi\} \;. \label{SE.5} \ee
This corresponds to a hyperboloid in ${\mycal M}$ given by the
equation $x^0 + 1 = \sqrt{1+\vec{x}^2}\;.  $ It is convenient to
introduce the following coordinate system $(x,v,\tau)$ in a
$\tM$-neighbourhood of $\scrip$: \beqa
\tau &=& y^0 -1/2\geq 0\;,\nonumber \\
x &=& (\sum (y^i)^2)^{\demi} -y^0 \geq 0  \;,\nonumber \\
y^i& = & (\sum (y^i)^2)^{\demi}n^i(v)\;, \label{SE.6} \eeqa
$n^i(v)\in S^{n-1}$, with $v=(v^A)$ denoting spherical coordinates
on $S^{n-1}$.  Equation~\eq{C.5} gives
\begin{equation}
  \label{SE.6.1}
   \Omega = x(2\tau +x+1) \approx x\;.
\end{equation}
 If we let $h$ denote the unit round metric
on $S^{n-1}$, we then have \be \eta = 2 dx d\tau + dx^2 +
(x+\tau+1/2)^2 h \;, \label{SE.7} \ee and  \beqa \Box_\eta
\tilde{f} &=& \frac{1}{(x+\tau+1/2)^{n-1}\sqrt{\det h}} \pmu
\left( (x+\tau+1/2)^{n-1}\sqrt{\det h}\;
\eta^{\mu\nu}\partial_\nu\tilde{f}\right)
\nonumber\\
&=& \{-\partial_\tau(\partial_\tau-2\partial_x) + {{n-1}\over
  x+\tau+1/2}\partial_x + {\triangle_h \over (x+\tau+1/2)^2}\} \tilde{f}\;,
\label{SE.8} \eeqa where $\triangle_h$ is the Laplace-Beltrami
operator of the metric $h$. We set \beqa & e_- = \partial_\tau \;,
\quad e_+ =
\partial_\tau - 2\partial_x\;,\quad e_A=\frac{1}{(x+\tau+1/2)
  }h_A\;,&\label{SE.9}
\\
&\phi_-=e_-(\tilde{f})\;, \quad \phi_+ = e_+ (\tilde{f})\;, &
\label{SE.10}\\ &\phi_{A}= \psi_A
=\frac{1}{(x+\tau+1/2)}h_A(\tilde{f})\;,&\label{SE.11} \eeqa where
$h_A$ denotes an $h$-orthonormal frame on $S^{n-1}$. We use the
symbol $D$ to denote the covariant derivative operator associated
to the metric $h$. (The usefulness of introducing two different
objects for ${h_A(\tilde{f})}/{(x+\tau+1/2)}$ will become clear
shortly.) Equation (\ref{SE.4}) implies the following set of
equations:
\begin{eqnarray} &
  \begin{array}{rrrcr}
 e_-(\phi_+) &  - \dea
\psi_A & - {n-1\over 2(x+\tau+1/2)}
\phi_+ &=& - {n-1\over 2(x+\tau+1/2)}\phi_- +a_+\;,\label{SE.12.1}\\
-e_A(\phi_+)   & + e_+ (\psi_A) & -\frac{1}{(x+\tau+1/2) }\psi_A
&=& b_A\;,\label{SE.13.1}\end{array} & \\ &  \begin{array}{rrrcr}
 e_- (\phi_A) &- e_A (\phi_-) &  + \frac{1}{(x+\tau+1/2)} \phi_A &=& a_A\;,
\label{SE.15.1}
\\
- \dea \phi_A & + e_+(\phi_-) &  + {n-1\over 2(x+\tau+1/2)}\phi_-
&=& {n-1\over2(x+\tau+1/2)}\phi_++ b_-\;,
\label{SE.14.1}\end{array}
& \\& e_-(\tilde{f}) = \phi_- \;, & \label{SE.17.1.0}\\&\quad e_+
(\tilde{f}) = \phi_+ \;,
 \label{SE.17.1} &
\eeqa with $a_A=b_A=0$ and
\begin{eqnarray}
&a_+\equiv b_-\equiv -G \equiv -\Omega^{-{n+3\over 2}}
H(x^\mu,\Omega^{n-1\over 2} \tilde{f}) \;.&\label{SE.18} \eeqa

\subsection{Existence of solutions, space derivatives estimates}
We note that the partial differential operator standing on the
left-hand-side of \eq{SE.12.1} is symmetric hyperbolic; the same
holds true for \eq{SE.15.1}, or for the joint system
\eq{SE.12.1}-\eq{SE.17.1}. Now, part of our technique consists in
deriving weighted energy estimates for symmetric hyperbolic
systems having the structure above, \emph{cf.\/} Section~\ref{S3}.
Each such system comes with his own estimates, so that for the
systems \eq{SE.12.1} and \eq{SE.15.1} we can obtain estimates with
different weights. This allows us to handle a reasonably wide
range of non-linearities, giving existence and blow-up control for
initial data in weighted Sobolev spaces (with conormal-type
blow-up at $\scrip$):
\begin{Theorem}\label{T2} Consider Equation (\ref{SE.1}) on
  ${\R}^{n,1}$ with initial data given on a hyperboloid $\hyp\supset
  \Sigma_{x_0,0}$ in Minkowski space-time, and satisfying
  \beqa\label{cpt1.0} \tf |_{\Sigma_{x_0,0}}
  \equiv\Omega^{-{n-1\over 2}}f|_{\Sigma_{x_0,0}} & \in&
  \HH_{k+1}^{\alpha} (\Sigma_{x_0,0})
\;,\\
  \partial_x (\Omega^{-{n-1\over 2}}f) |_{\Sigma_{x_0,0}} &\in &
  \mcC_0^{\alpha}( \Sigma_{x_0,0})\cap
  \HH_k^{\alpha-1/2}(\Sigma_{x_0,0})\;,\label{cpt1.1}\\
  \partial_\tau(\Omega^{-{n-1\over 2}}f) |_{\Sigma_{x_0,0}} &\in&
  \HH^\alpha_k
  (\Sigma_{x_0,0}) 
\label{cpt1} \;,  \eeqa  with some $k > {n\over2} +1$, $ -1<
\alpha <-1/2$. Suppose further that $H$ has a uniform zero of
order $\ell$ at $f=0,$ in the sense of (\ref{S2.71}), with \be
\ell\geq \cases{ 4\; , & $n=2$ , \cr
  3\; , & $n=3$ , \cr 2\; , & $n\geq 4$ .}
\label{condH} \ee Then:
\begin{enumerate}
\item There exists $0<\tau_+ \leq T \; (<x_0/2)$, depending only upon $x_0$ and
a bound on the norms of the initial data in the
  spaces appearing in \Eqs{cpt1.0}{cpt1}, and a solution $ f$ of
  Equation~\eq{SE.1}, defined on a set containing $\Omega_{x_{0},
  \tau_+}$, satisfying the given initial conditions, and satisfying
   $$\|\widetilde f\|_{L^\infty(\Omega_{x_0,\tau_+})}<\infty\;.$$
\item Further,  if $\tau_*$ is such that $f$
  exists on $\Omega_{x_0,\tau_*}$ and satisfies
  $\|\widetilde f\|_{L^\infty(\Omega_{x_0,\tau_*})}<\infty$, then for
$0 \leq \tau < \tau_*$
  we have $$\widetilde f |_{\Sigma_{x_0,\tau}} \in L^\infty
  (\Sigma_{x_0,\tau})\cap \HH_{k+1}^\alpha(\Sigma_{x_0,\tau})\; ,$$
$$\partial_\tau \widetilde f|_{\Sigma_{x_0,\tau}} \in \HH{^\alpha_k}
  (\Sigma_{x_0,\tau}) 
\;,\qquad \partial_x \widetilde f|_{\Sigma_{x_0,\tau}} \in
  \HH^{\alpha-{1\over 2}}_k (\Sigma_{x_0,\tau}) \cap
  \mcC_0^\alpha(\Sigma_{x_0,\tau})\;,$$
  with a $\tau$-independent bound on all the norms.
\end{enumerate}\end{Theorem}

\remarks

 1. Integration in $x$ of condition \eq{cpt1.1} implies  that
$\tf  \in L^\infty(\Sigma_{x_0,0})$.

 2. Some further information can be found in
Theorem~\ref{T2g} below.

3. If the inequality in \eq{condH} is not an equality for $n=2,3$
(no further restrictions for $n\ge 4$), then a proof similar, but
simpler, basing on Proposition~\ref{PL.1n} instead of \ref{PL.1},
leads to the same result with $\alpha=-1/2$.

\medskip

\proof As before, we write $\|f(\tau)\|_{\HH{^\alpha_k}}$ for
$\|f|_{\Sigma_{x_0,\tau}}\|_{\HH{^\alpha_k}(\Sigma_{x_0,\tau})}$,
\emph{etc.\/} Recall that the standard theory of hyperbolic
systems ({\em cf., e.g.,\/}~\cite[Chapter~16,
Vol.~III]{Taylor}$^{\mathit{\ref{Taylorref}}}$) shows that for any
$0<x_1\leq x_0$ there exists $T(x_1)>0$, satisfying $2x_1+T \leq
x_0$, and a solution $\tf $ of (\ref{SE.4}), defined on
$\Omega_{x_1,x_0,T}$, with initial data on $\Sigma_{x_1,x_0}$
obtained from those on $\Sigma_{x_0}$ by restriction. The idea of
the proof is to derive $x_1$-independent, weighted \emph{a priori}
estimates for the solution. These estimates will guarantee that
the existence time $T(x_1)$ does not shrink to zero as $x_1$ goes
to zero; they will also guarantee that the weighted Sobolev
regularity is preserved by evolution.  We start with the
following:
\begin{Lemma}
  \label{LT2} Under the hypotheses of Theorem~\ref{T2},
  consider on $\Omega_{x_1,x_0,T}$ the system \eq{SE.11}-\eq{SE.17.1},
  set
\begin{eqnarray}\nonumber
  E_\alpha(t)&=&\|\widetilde f(t)\|^2_{\HH^{\alpha}_{k}}+ \|\phi_-(t)\|^2
_{\HH^{\alpha}_{k}\Mxonexzerot} %
\\&&
+ \|\phi_+(t)\|^2_{\HH^{\alpha-\demi}_{k}\Mxonexzerot}
+\sum_{A}\|\phi_A (t)\|^2_{\HH^{\alpha}_{k}\Mxonexzerot}\;.
  \label{eq:LT2}
\end{eqnarray}
 Then there exists a $x_1$-independent constant $C$ such that
\beqa E_{\alpha}(t) & \leq & C\left\{ E_{\alpha}(0)e^{Ct}  +
\int_{0}^{t}e^{C(t-s)}S(s)ds\right\}\;, \label{N1.3.1.1}\;\eeqa
where
\beqa S(s)&\equiv& \left.\sum_A\|a_A(s)\|^2_{\HH^{\alpha}_{k}\Mxonexzeros}\right.
+\|a_+(s)\|^2_{\HH^{\alpha-1/2}_{k}\Mxonexzeros}  \nonumber \\
& & + \|b_-(s)\|^2_{\HH^{\alpha-1/2}_{k}\Mxonexzeros}
\left.
+ \sum_A\|b_A(s)\|^2_{\HH^{\alpha-1}_{k}\Mxonexzeros}
\right.
\;.
 \label{N1.3.1}
\eeqa
\end{Lemma}
\proof We wish, first, to apply Proposition~\ref{PL.1} to the
system consisting of Equation~(\ref{SE.14.1}) together with
$e_-(\tilde{f}) = \phi_- $; in order to do this we set
\[
\varphi=\left(\begin {array}{c}
\widetilde f \\
\phi_A\end {array}\right)
  \;,
\qquad \psi=\phi_-\;.
\]
We choose $E^\mu_\pm \pmu = e_\pm \otimes \Id$, we set \be
L\psi =\left(\begin{array} {c} 0 \\
    - e_A (\psi)\end{array}\right) \;,  \ee and we define
$$\widetilde E_\alpha(t) = \|\widetilde
f(t)\|^2_{\HH^\alpha_k}+ \|e_-(\widetilde f)(t)\|^2_{\HH^\alpha_k}
+ \sum_{A}\|e_A(\widetilde f)(t)\|^2_{\HH^\alpha_k}\;.
$$
The hypotheses ${\mcC}1-{\mcC}5$ of Proposition~\ref{PL.1} are
readily verified, and for any $\alpha < - \demi$ the
inequality~(\ref{L.10}) gives \beqa \widetilde E_\alpha(t) &\leq&
C \left\{\widetilde E_\alpha(0)e^{Ct} + \int_{0}^{t}e^{C(t-s)}
  \left(\sum_A\|a_A(s)\|^2_{\HH^{\alpha }} \right.\right.\nonumber\\
& & \left.+ \|\phi_+(s)\|^2_{\HH^{\alpha -1/2}} +
  \left.\|b_-(s)\|^2_{\HH^{\alpha-1/2}}\right)
  ds
\right\} \;.
 \label{N1.1}
\eeqa Next, we apply Proposition~\ref{PL.1} directly to
(\ref{SE.12.1}):
setting
\beqan
\hat E_{\alpha'}(t)&=&
\|e_+(\widetilde f)(t)\|^2_{\HH^{\alpha'}_{k}\Mxonexzerot}
+\sum_{A}\|e_A (\widetilde f)
(t)\|^2_{\HH^{\alpha'}_{k}\Mxonexzerot}\;, \eeqan for any $\alpha'
< -1/2$ it follows from (\ref{L.10}) that
\beqa \hat E_{\alpha'}(t) &\leq&
C \left\{\hat E_{\alpha'}(0)e^{Ct} + \int_{0}^{t}e^{C(t-s)}
  \left(\|a_+(s)\|^2_{\HH^{\alpha' }} \right.\right.\nonumber\\
& & \left.+ \|\phi_-(s)\|^2_{\HH^{\alpha' -1/2}} +
  \sum_A\left.\|b_A(s)\|^2_{\HH^{\alpha'-1/2}}\right)
  ds
\right\} \;.
 \label{N1.3}
\eeqa
%
We set $$E(t) = \widetilde E_\alpha(t)+ \hat
E_{\alpha-1/2}(t)\;.$$ It follows from (\ref{N1.1}) and
(\ref{N1.3}) with $\alpha'=\alpha -1/2$ that we have
\be
E(t)\leq C \left(E({0})e^{Ct} + \int_{0}^{t} e^{C(t-s)}(E(s) +
S(s))ds\right)\;, \label{N1.4.0} \ee with $S(s)$ as in
\eq{N1.3.1}. Equation~\eq{N1.3.1.1} with $E_\alpha$
replaced\footnote{The constant $C$ in Equation~\eq{N1.3.1.1} does
not
  necessarily coincide with that in \eq{N1.4.0}.} by $E$ follows now
from  Gronwall's Lemma. Since $E_\alpha$ is equivalent to $E$, our
claims follow. \qed

Returning to the proof of Theorem~\ref{T2}, Lemma~\ref{LT2}
applied to \eq{SE.12.1}-\eq{SE.17.1.0} gives (recall that $G$ was
defined in \eq{SE.18}) \be E_\alpha(t)\leq C
\left(E_\alpha({0})e^{Ct} + \int_{0}^{t}
e^{C(t-s)}\|G(s)\|^2_{\HH^{\alpha-1/2}_{k}}ds\right)\;.
\label{N1.4} \ee By hypothesis  the function $H$ appearing in
(\ref{SE.1}) has a uniform zero of order $\ell \geq 2,$ in the
sense of (\ref{S2.71}); we wish to use (\ref{S2.7.2}) to control
the term containing $G(s)$ in (\ref{N1.4}). This requires an
$L^\infty$ bound on $\widetilde f$, which will be obtained next.
As $k > n/2+1$, the Sobolev embedding (\ref{S2.4}) gives \be
\|e_-(\widetilde f)(s)\|^2_{\mcC^{\alpha}_{1}\Mxonexzeros}+ \|e_+
(\widetilde f) (s)\|^2_{\mcC^{\alpha-1/2}_{1}\Mxonexzeros} +
\|e_A(\widetilde f)(s)\|^2_{\mcC^\alpha_1\Mxonexzeros} \leq C
E_\alpha(s). \ee Now the conditions (\ref{condH}) on $n$ and
$\ell$ give
$$
|G(\tau)| \leq C \|\widetilde f(\tau)\|_{L^{\infty}
}^\ell\;
x^{\ell{n-1\over 2} -{n+3\over 2}}  \leq C \|\widetilde f(\tau)\|_{L^{\infty}
}^\ell\; x^{-1/2}  \;,$$
so that (recall that $\alpha<-\frac 12$)
\begin{equation}
  \label{eq:gd}
  \|G(\tau)\|_{\mcC_0^\alpha
} \leq
C\;\|\widetilde f(\tau)\|^\ell_{L^{\infty}
}\;.
\end{equation}
{}From (\ref{SE.12.1}) we have
\be
\partial_\tau\phi_+ - {n-1 \over 2(x+\tau+1/2)} \phi_+ = \dea \psi_A - {n-1\over 2(x+\tau+1/2)}\phi_- -G\;,
\label{equatauphi} \ee and \eq{eq:gd} together with
Proposition~\ref{Plemme} yield \beqa
\|\phi_+(t)\|_{\mcC_0^\alpha\Mxonexzerot}& \leq& C e^{Ct}\|\phi_
+(0)\|_{\mcC_0^\alpha\Mxonexzerot}\nn\\ &&
+C\int_0^te^{C(t-s)}(\|\dea
\psi_A (s)\|_{\mcC_0^\alpha\Mxonexzerot} 
+ \|\phi_- (s)\|_{\mcC_0^\alpha\Mxonexzerot} +
\|G(s)\|_{\mcC_0^\alpha\Mxonexzerot}) \;ds \nn\\
& \leq& C
e^{Ct}\|\phi_+(0)\|_{\mcC_0^\alpha\Mxonexzerot}+\int_0^te^{C(t-s)}
C(E_\alpha(s), \|\widetilde f(s)\|_{L^\infty}) ds \;,
\label{eq:gd1}\eeqa for some continuous function
$C(E_\alpha(\cdot), \|f(\cdot)\|_{L^\infty})$. Integration of
$\partial_x \widetilde f = {1\over
   2}(\phi_--\phi_+)$ over $[x,x_0-2\tau]$ gives
\beqa\nonumber & |\widetilde f(\tau,x)|\leq |\widetilde
f(\tau,x_0-2\tau)| + \displaystyle\frac 12
\|(\phi_--\phi_+)(\tau)\|_{\mcC^\alpha_0}\int_x^{x_0-2\tau}
 s^{\alpha} ds \;.&
\eeqa For any $0\leq \tau\leq \tau_*<x_0/2$ the
$f(\tau,x_0-2\tau)$ term is estimated by a multiple of the initial
energy in a standard way, which leads to the estimate (recall that
$\alpha>-1$) \beqa \|\widetilde f(\tau)\|_{L^\infty\Mxonexzerot}
&\leq & CE_\alpha(\tau)+ C
e^{C\tau}\|\phi_+(0)\|_{\mcC_0^\alpha\Mxonexzerot} \nonumber \\ &&
+ \displaystyle\int_0^\tau e^{C(\tau-s)}C(E_\alpha(s),
\|\widetilde f(s)\|_{L^\infty}) ds \;. \label{N1.5} \eeqa Next,
$$\|G(s)\|_{\HH^{\alpha-1/2}_{k}\Mxonexzeros} \leq C
\|H(s,\cdot,x^{{n-1\over 2}} \widetilde f)\|_{\HH^{\alpha-1/2+
{n+3\over 2}} \Mxonexzeros}\;,$$ and our hypothesis that $H$ has a
uniform zero of order $\ell$ together with (\ref{S2.7.2}) gives
$$\|G(s)\|_{\HH^{\alpha-1/2}_{k}\Mxonexzeros} \leq C\left(\|\widetilde
f(s)\|_{L^\infty\Mxonexzeros}\right) \|\widetilde
f\|_{\HH_k^{\alpha+
  {n+2\over 2}  - l {n-1\over 2}}\Mxonexzeros}\;.$$
In view of (\ref{N1.5}) this can be estimated by a function of
$E_\alpha(s)$ and of $\|\widetilde f(s)\|_{L^\infty}$,
 \beqa
\|G(s)\|^2_{\HH^{\alpha-1/2}_{k}\Mxonexzeros} &\leq&
C\left(\|\widetilde f(s)\|_{L^\infty\Mxonexzeros}\right)
\|\widetilde f(s)\|^2_{\HH^{\alpha}_{k}\Mxonexzeros} \nn
\\ &\leq & C\left(\|\widetilde
f(s)\|_{L^\infty\Mxonexzeros}\right) E_\alpha(s)\;, \label{N1.5.1}
\eeqa provided that \be l \geq  {n+2\over n-1} \label{N1.6} \ee
(which coincides again with \eq{condH}).
 If (\ref{N1.6})
holds, from (\ref{N1.4}) and \eq{N1.5} we obtain
\begin{eqnarray}
  \|\widetilde f(\tau)\|_{L^\infty}+ E_\alpha(\tau) &\leq&
  Ce^{C\tau}\left(E_\alpha({0}) + \|\partial_x\widetilde
f(0)\|_{\mcC^\alpha_0} + \|\partial_\tau\widetilde f (0)\|_{\HH
^\alpha_{k}}\right) \nonumber \\ &&
 + \int_{0}^{\tau}
\Phi\left(\tau,s,\|\widetilde f(s)\|_{L^\infty},
E_\alpha(s)\right) ds\;,
  \label{eq:glbd}\end{eqnarray} for some constant $C$, and for a
  function $\Phi$ which is bounded on bounded sets. It then easily
  follows that there exists a time $\tau_+$ and a constant
$M$,
  depending only upon $x_0$ and a bound on the norms of the initial data in the
  spaces appearing in \Eqs{cpt1.0}{cpt1}, such
  that $\|\widetilde f(\tau)\|_{L^\infty}$ and $E_\alpha(\tau)$ remain
  bounded by $M$ for $0\leq \tau\leq\tau_+$. Since all the objects above
  were $x_1$-independent, so is $\tau_+$. By the usual continuation
  criterion ({\em cf., e.g.,\/}~\cite[Proposition~1.5, Chapter~16,
  Vol.~III]{Taylor}\footnote{\label{Taylorref}In that reference
  symmetric hyperbolic systems on a torus are considered; however
  simple domain of dependence considerations show that the results
  there apply to the setup here.})  the solution exists on
  $\Omega_{x_1,x_0,\tau_+}$ for all $x_1$; it thus follows that the
  maximally extended solution of the initial value problem considered
  here exists on a set which includes $\Omega_{x_0,\tau_+}$.

To establish point 2, suppose that  a global a-priori $L^\infty$
bound on $\tf$ is known. Then \eq{N1.4} and \eq{N1.5.1} give a
linear integral inequality on $E_\alpha$, and Gronwall's Lemma
gives a global bound on $E_\alpha$. Arguments of the last part of
the proof of point 1 yield the result. \qed

For the purpose of estimating time derivatives of the solutions we
will need a generalization of Theorem~\ref{T2}, which covers the
equations contained by time-differentiating
\Eqs{SE.12.1}{SE.17.1}. There are lots of ways to relax those
hypotheses; for simplicity we shall only make those
generalizations which are strictly necessary for the arguments in
the next section to go through. First, the fact that $f$ is scalar
valued plays no role in our considerations above; henceforth we
assume that $f$ has values in $\R^N$ for some $N\geq 1$.  Next,
the definitions \eq{SE.9} of $e_{\pm}$ and $e_A$ will be kept.
We will consider systems of the form
\begin{deqarr}
&P\left(\begin{array}{c} \varphi \\ \psi \end{array}\right) +
\left(\begin{array}{cc} B_{11} & B_{12} \\ B_{21} & B_{22}
\end{array}\right) \left(\begin{array}{c} \varphi \\ \psi
\end{array}\right) = \left(\begin{array}{c} a \\ b
\end{array}\right) + G\;, & \label{td1a}
\\
&\varphi=\left(\begin{array}{c} \phi_+ \\ \phi_A
\end{array}\right)
\;,\quad \psi=\left(\begin{array}{c} \phi_- \\ \psi_A
\end{array}\right)
&\label{td1b} \arrlabel{td1}
\end{deqarr}
together with
\begin{deqarr}&\phi_{A}
=\frac{1}{(x+\tau+1/2)}h_A(\tilde{f}) + B_{A,\phi} \tilde f\;,&\label{SE.11a} \\
&\psi_A =\frac{1}{(x+\tau+1/2)}h_A(\tilde{f})+ B_{A,\psi} \tilde
f\;,&\label{SE.11b}
\\ &e_-(\tilde{f}) = B_{0}\phi_- + B_{1}\tf\;, & \label{td2}\\
&e_+ (\tilde{f}) = \phi_+
 \;, & \label{td3}
\arrlabel{td4}
\end{deqarr}
for some matrix valued functions $B_{A,\phi}$, $B_{A,\psi}$,
$B_0$, $B_1$, with $B_0$ --- invertible. Here \be
\label{td4P} P = \left(\begin{array}{cc} e_- & \ell^A D_A \\
(\ell^A)^tD_A & e_+ \end{array}\right) \ee is the (geometric)
principal part of \Eqs{SE.12.1}{SE.14.1}. The nonlinear term
$G=G(x^\mu,\tf)$ will be labeled as \be
\label{td4G}G=\left(G_{e_+(\phi_-)},
G_{e_+(\psi_A)},G_{e_-(\phi_A)}, G_{e_+(\phi_-)}\right)\;,\ee with
the order of the components following that of
\Eqs{SE.12.1}{SE.14.1}. The $B_{ab}$'s will be labeled as
$B_{\phi_-,\phi_+}$, $B_{\phi_-,\phi_A}$, {\em etc.}; for example,
in this notation,
the second of \Eqsone{SE.14.1} takes the form
\begin{eqnarray}
\lefteqn{ e_+(\phi_-)
  =  \dea \phi_A  -B_{\phi_-,\phi_-}\phi_- }&& \nn\\
&&
-B_{\phi_-,\phi_+}\phi_+-B_{\phi_-,\phi_A}\phi_A-B_{\phi_-,\psi_A}\psi_A+
b_-+G_{e_+{(\phi_-)}}\;, \label{td8.0}
\end{eqnarray}
with actually $B_{\phi_-,\phi_A}=B_{\phi_-,\psi_A}=0$.

Some effort will be needed to prove the information of point 3 of
the theorem that follows; this is needed to be able to iteratively
apply that theorem in the next section:

\begin{Theorem}\label{T2g} Consider the system \eq{td1}-\eq{td4} with
\begin{eqnarray}
\!\!\!\!\!\!\!\!\!\!\lefteqn{\!\!\!\!\!\!\!\!\!\! \|a(\tau)\|_{\HH
^\alpha_k} +\|b(\tau)\|_{\HH ^\alpha_k}+
\sup_{a,b=1,2}\|B_{ab}(\tau)\|_{\mcC^0_k} +
\sup_{A=1,2;\;\lambda=\phi,\psi}\|
B_{A,\lambda}(\tau)\|_{\mcC^0_{k}}} && \nn \\
&&\phantom{xxxxxxxxxx} +\|B_{0}(\tau)\|_{\mcC^0_k}
+\|B_{0}^{-1}(\tau)\|_{L^\infty} +\|B_{1}(\tau)\|_{\mcC^0_k}
 \leq \tilde C\;,
\label{td5}\end{eqnarray} for some constant $\tilde C$, and
suppose that
\begin{equation}
G(x^\mu,\tf)=\Omega^{-(n+3)/2}H(x^\mu,\Omega^{(n-1)/2}\tf)\;,
\label{td6}
\end{equation}
with $G_{e_-(\phi_A)}=0$, and with $H$ having a uniform zero of
order $\ell$ in the sense of (\ref{S2.71}), with $\ell$ satisfying
\eq{condH}. If the initial data satisfy \eq{cpt1.0}-\eq{cpt1} with
some $k > {n\over2} +1$, $ -1< \alpha <-1/2$, then:
\begin{enumerate}
\item The
conclusions of point 1. of Theorem~\ref{T2} hold with a time
$\tau_+$ depending only upon the constant $\tilde C $ in \eq{td5}
and a bound on the norms of the initial data in the spaces
appearing in \Eqs{cpt1.0}{cpt1}.
\item The
conclusions of point 2. of Theorem~\ref{T2} hold.
\item Under the
hypotheses of point 2. of  Theorem~\ref{T2} 
we also have
\begin{equation}
\|(x+2\tau)\partial_\tau\widetilde
f\|_{L^\infty(\Omega_{x_0,\tau_*})}<\infty\;. \label{td7}
\end{equation}

\end{enumerate}
\end{Theorem}

\remarks 1. The condition $G_{e_-(\phi_A)}=0$ can be weakened to
\begin{equation}
G_{e_-(\phi_A)}(x^\mu,\tf)=\Omega^{-(n+2)/2}H_{e_-(\phi_A)}
(x^\mu,\Omega^{(n-1)/2}\tf)\;, \label{td6.1}
\end{equation}
for some function $H_{e_-(\phi_A)}$ with a uniform zero of order
$\ell$. Similarly it suffices to assume that
\begin{equation}
G_{e_+(\psi_A)}(x^\mu,\tf)=\Omega^{-(n+4)/2}H_{e_+(\psi_A)}
(x^\mu,\Omega^{(n-1)/2}\tf)\;, \label{td6.2}
\end{equation}for some function $H_{e_+(\psi_A)}$ with a uniform zero of order
$\ell$.

2. If the inequality in \eq{condH} is not an equality for $n=2,3$
(no further restrictions for $n\ge 4$), then the result remains
true with $\alpha=-1/2$, see Remark 3. after Theorem~\ref{T2}.

\medskip

\proof Let us start by remarking that, because $\psi_A=\phi_A$, in
equations such as \eq{td8.0} we can replace $B_{\phi_-,\phi_A}$ by
$B_{\phi_-,\phi_A}+B_{\phi_-,\psi_A}$ obtaining a system in which
$B_{\phi_-,\psi_A}=0$. Proceeding similarly with the other
equations we may thus without loss of generality assume that
\begin{equation}
B_{*,\psi_A}=0\;.
\end{equation}
The proof of points 1 and 2 is then identical to that of
Theorem~\ref{T2}, with the following minor changes:
\Eq{equatauphi} is replaced by the equation
\begin{eqnarray}
\lefteqn{ e_-(\phi_+) +B_{\phi_+,\phi_+}\phi_+
  =
}&& \nn\\
&&  \dea \phi_A -B_{\phi+-,\phi_-}\phi_- -B_{\phi_+,\phi_A}\phi_A
+a_++G_{e_-{(\phi_+)}} 
\label{td8.01}
\end{eqnarray}
to which Proposition~\ref{Plemme} still applies, recovering
\eq{eq:gd1}. Further, the equation $\partial_x \tf=
(\phi_--\phi_+)/2$ has to be replaced by
$$\partial_x \tf+\frac {B_1} 2\tf= \frac{B_0\phi_--\phi_+}{2}\;,$$
and the desired conclusion is obtained by
Proposition~\ref{propositionx1}. The remaining arguments do not
require any modifications.

To prove point 3, from \eq{td8.0} we obtain
\begin{eqnarray}
\lefteqn{ e_+[(x+2\tau)\phi_{-}] =}&& \nonumber\\ &&
(x+2\tau)\left(\dea \phi_A  - B_{\phi_-,\phi_-}\phi_- \right.
\nn\\ &&\left.- B_{\phi_-,\phi_A}\phi_A-B_{\phi_-,\phi_+}\phi_++
b_-+ G_{e_+{(\phi_-)}}\right)\;, \label{td8}
\end{eqnarray}
{}From \Eqsone{eq:gd}, \eq{eq:gd1}, and \eq{td2} together with $$
\phi_-,\phi_A\in \HH ^\alpha_k\subset \mcC^\alpha_0\;, \qquad \dea
\phi_A\in \HH ^\alpha_{k-1}\subset \mcC^\alpha_0\;,
$$ we obtain
\begin{eqnarray*}
 e_+[(x+2\tau)\phi_{-}]
&\leq& \hat C x^{-\alpha}\;,
\end{eqnarray*}
for some constant $\hat C$ depending only upon the initial data
and $\|\tf\|_{L^\infty(\Omega_{x_0,\tau_*})}$. Integrating as in
the identity \eq{s6}  we arrive at
\begin{eqnarray*}
\lefteqn{|B_0^{-1}\left\{(x+2\tau)(\partial_\tau\tf-B_1\tf)(x,v,\tau)\right\}|
}& & \\ &&\leq
|B_0^{-1}\left\{(x+2\tau)\partial_\tau\tf(x+2\tau,v,0)\right\}|
+ C\left( \|\tf\|_{L^\infty(\Omega_{x_0,\tau_*})} + \hat C\right)
\\ && \leq C\left(
\|\partial_\tau\tf\|_{\mcC^{-1}_0} +
\|\tf(0)\|_{L^\infty(\Omega_{x_0,\tau_*})} + \hat C\right)
\;,
\end{eqnarray*}
and \Eq{td7} follows. \qed

\subsection{Estimates on the time derivatives of the solutions}
\label{Stdphi} So far we have established existence of solutions
with initial data in weighted Sobolev spaces, as well as weighted
estimates on the space-derivatives of the solutions. The next step
in proving polyhomogeneity is to establish estimates on
time-derivatives. Similarly to the linear case, the question of
corner conditions arises. In order to handle that, we introduce an
index $m$, which corresponds to the number --- perhaps zero --- of
corner conditions which are satisfied by the initial data.  Next,
the definition (\ref{S2.71}) of a uniform zero of order $l$ has to
be strengthened by adding conditions on time-derivatives: we shall
require that for all $M\in\R$,  $|\prho|\le M$, $0\leq i\leq
\min(k,l)$ and for all $0\leq j\leq m$ there exists a constant
$\hat{C}=\hat C(M,m,k)$ such that we have \be\label{S2.71t}
\left\|\frac{\partial^{i+j} F(\tau,\cdot,\prho)}{\partial
    \prho^i\partial \tau^j}\right\|_{{\mcC}_{k+m-i-j}^0} \leq \hat{C}
|\prho|^{l-i}\;.\ee We start with the following:
\begin{Theorem}\label{T2t} Let $\N\ni m\geq 0$, consider  a solution
$f:\Omega_{x_0,\tau_*}\to \R$ of Equation (\ref{SE.1}) satisfying
$$\|\widetilde f\|_{L^\infty(\Omega_{x_0,\tau_*})}<\infty\;,$$ and
suppose that $ $
\beqa\label{cpt1.0t} 0\leq i \leq m+1\qquad
\partial^i_\tau\widetilde f|_{\Sigma_{x_0,0}} & \in&
\HH_{k+m+1-i}^{\alpha} (\Sigma_{x_0,0})\;,\\ 0\leq i \leq m \qquad
\partial_x \partial^i_\tau\widetilde f |_{\Sigma_{x_0,0}} &\in &
\mcC_0^{\alpha}( \Sigma_{x_0,0})\cap
\HH_{k+m-i}^{\alpha-1/2}(\Sigma_{x_0,0})\label{cpt1.1t} \;,  \eeqa
  with  some $k > {n\over2} +1$
and $ -1< \alpha <-1/2$. Suppose, further, that $H$ is smooth in
$f$ and has a uniform zero of order $\ell$ at $f=0,$ in the sense
of (\ref{S2.71t}), with $\ell$ as in \Eq{condH}.  Then for $0 \leq
\tau < \tau_*$ and for $0\leq i \leq m$, $0\leq j+i < k+m-n/2$  we
have
\begin{deqarr}
[(\tau+2x)\partial_\tau]^j\partial_\tau^i\widetilde f
|_{\Sigma_{x_0,\tau}} \in L^\infty
  (\Sigma_{x_0,\tau})\cap \HH_{k+m+1-i-j}^\alpha(\Sigma_{x_0,\tau})\;
,
\\
 \partial_x
[(\tau+2x)\partial_\tau]^j\partial_\tau^{i}\widetilde
f|_{\Sigma_{x_0,\tau}} \in
  \HH^{\alpha-{1\over 2}}_{k+m-i-j} (\Sigma_{x_0,\tau}) \cap
  \mcC_0^\alpha(\Sigma_{x_0,\tau})\;,
\arrlabel{td11}\end{deqarr} and
\begin{equation}0\leq p < k-n/2 \qquad
[(\tau+2x)\partial_\tau]^p\partial_\tau^{m+1} \widetilde
f|_{\Sigma_{x_0,\tau}} \in \HH{^\alpha_{k-p}}
  (\Sigma_{x_0,\tau}) 
\;, \label{td12}\end{equation} with $\tau$-independent bounds on
the norms.
\end{Theorem}

\remark As before, in dimensions $n\ge 4$ the result remains valid
for $\alpha=-1/2$; in dimensions $n=2,3$ the value $-1/2$ for
$\alpha$ is allowed if the inequality in \eq{condH} is not an
equality.

The proof below actually proves the analogous result for the
systems considered in Theorem~\ref{T2g}, provided that obvious
time-derivative conditions on the coefficients are added to
\eq{td5}, the simplest possibility being
\begin{eqnarray}
\lefteqn{
\|\partial_\tau^i a(\tau)\|_{\HH ^\alpha_{k+m-i}}
+\|\partial_\tau^i b(\tau)\|_{\HH ^\alpha_{k+m-i}}+
\sup_{a,b=1,2}\|\partial_\tau^i B_{ab}(\tau)\|_{\mcC^0_{k+m-i}} }&& \nn \\
&&
+ \sup_{A=1,2;\;\lambda=\phi,\psi}\|
\partial_\tau^i B_{A,\lambda}(\tau)\|_{\mcC^0_{k+m-i}}+\|\partial_\tau^i
B_{0}(\tau)\|_{\mcC^0_{k+m-i}}
+\|\partial_\tau^i B_{1}(\tau)\|_{\mcC^0_{k+m-i}}
 \leq \tilde C\;,\nonumber \\ &&
\label{td5t}\end{eqnarray} with $0\le i \le m+k$; the same remark
applies to Corollary~\ref{C1t} below. Before passing to that
proof, we note that an important consequence of Theorem~\ref{T2t}
is that corner conditions will hold at any time $\tau>0$,
regardless of whether or not they hold at $\tau=0$:

\begin{Corollary}
\label{C1t} Under the conditions of point 2 of Theorem~\ref{T2},
for any $0<\tau<\tau_*$ and for $0\leq i< k-1-n/2 $ we have
$$ \qquad \partial_\tau^i\widetilde f |_{\Sigma_{x_0,\tau}} \in L^\infty
  (\Sigma_{x_0,\tau})\cap \HH_{k+1-i}^\alpha(\Sigma_{x_0,\tau})\; ,$$
$$\partial_\tau^{i+1} \widetilde
f|_{\Sigma_{x_0,\tau}} \in \HH{^\alpha_{k-i}}
  (\Sigma_{x_0,\tau}) 
\;,$$
$$ \partial_x \partial_\tau^{i}\widetilde
f|_{\Sigma_{x_0,\tau}} \in
  \HH^{\alpha-{1\over 2}}_{k-i} (\Sigma_{x_0,\tau}) \cap
  \mcC_0^\alpha(\Sigma_{x_0,\tau})\;.$$
\end{Corollary}


We shall need the following simple Lemma:
\begin{Lemma}
\label{L1t} Let  $F(x^\mu,p)$ be a function which is smooth in $p$
at fixed $x^\mu$ and suppose that it has a uniform zero of order
$\ell\geq 1$ in $p$. Then
\begin{enumerate}
\item For all $i\in \N$ the function $\partial_\tau^i
\left(F(x^\mu,u(x^\mu)\right)$ has a uniform zero of order $\ell$,
when viewed as a function of $(u,\partial_\tau u,\ldots,
\partial^i_\tau u)$.
\item
Let $H=\partial_p F$, then $H$ has a uniform zero of order
$\ell-1$.
\end{enumerate}
\end{Lemma}

\proof Let $u=(u^i)$; smoothness of $F$ in $p$ allows us to write
\be\label{td9}F(\vec x,\tau,u)= A_{i_1\ldots i_\ell}u^{i_1} \ldots
u^{i_\ell}\;, \ee with some coefficients $A_{i_1\ldots
i_\ell}=A_{i_1\ldots i_\ell}(\vec x,\tau,u)$ which are smooth in
$u$, and totally symmetric in $i_1,\ldots,i_\ell$; recall that the
summation convention is used throughout. Point 2 immediately
follows from \eq{td9}. From that equation we also obtain
\beqan\partial_\tau F(\tau,\vec x,u) &=& \left(\partial_\tau
A_{i_1\ldots i_\ell}+ \partial_{u^i}A_{i_1\ldots
i_\ell}\partial_\tau u^i\right) u^{i_1} \ldots u^{i_\ell}
\\ && + \ell A_{i_1\ldots i_\ell}u^{i_1} \ldots u^{i_{\ell-1}}
\partial_\tau u^{i_\ell}\;,
\eeqan which proves point 1 for $i=1$. The result then follows by
straightforward induction. \qed

We can pass now to the proof of Theorem~\ref{T2t}:

\medskip

\proof We assume that \Eqs{td1}{td4} are satisfied;
Theorem~\ref{T2g} shows that \eq{td11}-\eq{td12} hold with
$i=j=p=0$. Consider the vector-valued function
$$ (\tf,(x+2\tau)\partial_\tau\tf,\varphi,
(x+2\tau)\partial_\tau\varphi,\psi,
(x+2\tau)\partial_\tau\psi)\;,$$ so that the new function $\tf$ in
\eq{td4} is $(\tf,(x+2\tau)\partial_\tau\tf)$, while the new
functions $\varphi$, respectively $\psi$, in \eq{td1b} are $
(\varphi, (x+2\tau)\partial_\tau\varphi)$, respectively $ (\psi,
(x+2\tau)\partial_\tau\psi)$.
We claim that   a set of equations of the form \eq{td1}-\eq{td4}
holds for those new functions. Consider, for instance, \Eq{td2};
set
$$\hf:= (x+2\tau)\partial_\tau\tf\;,\quad
\hphi_-:=(x+2\tau)\partial_\tau\phi_-\;,$$ {\em etc.}, we have
\begin{eqnarray*}
e_-(\hf) & = &
\partial_\tau\left((x+2\tau)(B_0\phi_-+B_1\tf)\right)
\\ &= & B_0\hphi_- + (2B_0+(x+2\tau)\partial_\tau B_0)\phi_-
\\ & & +B_1\hf + (2B_1+(x+2\tau)\partial_\tau B_1)\tf\;,
\end{eqnarray*}
which is linear in $(\tf,\hf,\phi_-,\hphi_-)$. In fact
\begin{eqnarray*}
e_-\left(\begin{array}{c} \tf \\ \hf \end{array}\right) &=&
\left(\begin{array}{cc} B_0 & 0 \\ 2B_0+(x+2\tau)\partial_\tau B_0
& B_0\end{array}\right) \left(\begin{array}{c} \phi_-\\ \hphi_-
\end{array}\right)
\\ && +
\left(\begin{array}{cc} B_1 & 0 \\ 2B_1+(x+2\tau)\partial_\tau B_1
& B_1\end{array}\right) \left(\begin{array}{c} \tf\\ \hf
\end{array}\right) \;,
\end{eqnarray*}
and the new matrix $B_0$ is again invertible, as desired. Next,
\begin{eqnarray*}
e_-(\hphi_+) & = &
\partial_\tau\left((x+2\tau)\partial_\tau\phi_+\right)
\\ &= & \partial_\tau\left((x+2\tau)\left(
- D_A\phi_A - B_{\phi_+\phi_-}\phi_- \right.\right.\\ && \left.
\left.- B_{\phi_+\phi_A}\phi_A - B_{\phi_+\phi_+}\phi_+ +a_+ +
G_{e_-(\phi_+)} \right)\right)
\\ &= &
-D_A\hphi_A - B_{\phi_+\phi_-}\hphi_-- B_{\phi_+\phi_A}\hphi_A -
B_{\phi_+\phi_+}\hphi_+
\\ && 
+\; \mathrm{linear}(\varphi,\psi) +\widehat a_+ +
G_{e_-(\hphi_+)}
\;,
\\
\widehat a_+ & = & -2D_A\phi_A + \partial_\tau a_+ \ \in\HH
^{\alpha}_{k+m-1}\;,
\\
 G_{e_-(\hphi_+)} & = & \partial_\tau \left(G_{e_-(\phi_+)}\right)(x+2\tau) \;,
\end{eqnarray*}
where ``$\mathrm{linear}$'' denotes terms which are linear in the
relevant variables. The equation for $e_-(\hphi_A)$ is handled in
a similar way. The equations involving only $e_+$ or $\partial_A$
are straightforward, since those operators commute with
multiplication by $(x+2\tau)$.  By Lemma~\ref{L1t} the new
non-linearity has again a zero of order $\ell$, when considered as
a function of $(\tf,{(x+2\tau)\partial_\tau}\tf)$. In order to
apply Theorem~\ref{T2g} we need to check that the initial data are
in the right spaces. Clearly
$$((x+2\tau)\partial_\tau\tf)(0) = x\partial_\tau\tf(0) \in
\HH ^{\alpha+1}_{k+m}\subset \HH ^{\alpha}_{k+m}\cap L^\infty\;,$$
$$\left(\partial_x((x+2\tau)\partial_\tau\tf)\right)(0) =
\left(\partial_\tau \tf+ x\partial_x\partial_\tau\tf\right)(0) \in
\HH ^\alpha_{k+m-1}\subset \mcC^\alpha_{0}\cap\HH
^{\alpha-1/2}_{k+m-1}\;.$$ Condition \eq{cpt1} requires some more
work: \beqan
\left(\partial_\tau((x+2\tau)\partial_\tau\tf)\right)(0) &=&
\left(2 \partial_\tau \tf+ x\partial_\tau^2\tf\right)(0)
\\ &=&
\left(2 \partial_\tau \tf+
x(2\partial_x+e_+)\partial_\tau\tf\right)(0)
\\ &=&
\left(2 \partial_\tau \tf+2x\partial_x\partial_\tau\tf+
xe_+\left(B_0\phi_-+B_1\tf \right)\right)(0)\;. \eeqan The first
two terms are obviously in $\HH ^\alpha_{k+m-1}$, and so is
$xe_+(B_1\tf)= x(\partial_\tau-2\partial_x)(B_-\tf)$. \Eq{td8.0}
gives
\begin{eqnarray*}
 \left(xe_+(\phi_-)\right)(0)
&=&  x\left(\dea \phi_A  -B_{\phi_-,\phi_-}\phi_--B_{\phi_-,\phi_+}\phi_+\right.  \nn\\
&&\left. -B_{\phi_-,\phi_A}\phi_A-B_{\phi_-,\psi_A}\psi_A+
b_-+G_{e_+{(\phi_-)}}
\right)(0)\;.
\end{eqnarray*} The desired
property $\left(xe_+(B_0\phi_-)\right)(0)\in \HH ^\alpha_{k+m-1}$
follows immediately; the only non-trivial term is
$xG_{e_+{(\phi_-)}}$, the $\HH^{\alpha}_{k+m+1}$ norm of which can
be estimated by a function of $\|\tf(0)\|_{L^\infty}$ and
$\|\widetilde f(0)\|_{\HH^{\alpha}_{k+m+1}\Mxonexzeros}$, {\em
cf.\/} \Eq{N1.5.1}.
Now, ${(x+2\tau)\partial_\tau}\tf $ 
is uniformly bounded on $\Omega_{x_0,\tau_*}$ by point 3 of
Theorem~\ref{T2g}, so that we can apply point 2 of
Theorem~\ref{T2g} to conclude that \Eqs{td11}{td12} hold with
$j=p=1$ and $m=0$; straightforward induction establishes
Theorem~\ref{T2t} for the remaining $j$'s and $p$'s.

Consider, now, $m=1$; the result already established with $m=0$
shows that $\partial_\tau\tf(\tau)$ exists and satisfies \eq{td11}
with $i=1$ for any $\tau >0$; similarly \eq{td12} holds with $m=1$
for any $\tau>0$. Now, a calculation similar (but simpler) to the
one done above shows that $(\tf,\partial_\tau\tf)$ satisfies a
system of equations of the form \eq{td1}-\eq{td4} with initial
data satisfying the conditions of Theorem~\ref{T2g} by hypothesis;
the uniform bounds on some interval $[0,\tau_+)$ follow by point 1
of that theorem. We therefore have
$$
\|(\tf,\partial_\tau\tf)\|_{L^\infty(\Omega_{x_0,\tau_*})}<\infty\;.
$$ We can then apply the result already established for $m=0$ to
the system of equations satisfied by $(\tf,\partial_\tau\tf)$ to
obtain the conclusion of Theorem~\ref{T2t} with $m=1$. An
induction upon $m$ finishes the proof. \qed

\subsection{Polyhomogeneous solutions}
\label{Spsnl}

The aim of this section is to establish polyhomogeneity of
solutions of a large class of semi-linear systems of the form
\begin{deqarr}
  \partial_\tau\varphi + B_{11}\varphi + B_{12}\psi&=& L_{11}\varphi +
  L_{12}\psi + a +G_\varphi
  \label{nequa1}\;,\\ \partial_x\psi+ B_{21}\varphi + B_{22}\psi &=&
  L_{21}\varphi +
  L_{22}\psi + b +G_\psi\label{nequa2}\;,\arrlabel{nequa}
\end{deqarr}
with a nonlinearity $$G= (G_\varphi,G_\psi)$$ of the form
\begin{equation}
  \label{m0}
  G= x^{-p\delta} H(x^\mu,x^{q\delta} \psi_1, x^{q\delta+1} \psi_2, x^{q\delta+1} \varphi)\;.
\end{equation}
Here we have decomposed $\psi$ as
\be\label{psidecomp}\psi=\left(\begin{array}{c} \psi_1\\ \psi_2
\end{array}\right)\;; \ee this is motivated by different {\em a
priori\/} estimates we have at our disposal for the appropriately
defined components $\psi_1$ and $\psi_2$ of $\psi$ in the
applications we have in mind. Polyhomogeneity of solutions of
\eq{SE.1} will follow as a special case, see \thm{T2phg} below. We
will need to impose various restrictions on the function $H$, in
order to do that some terminology will be needed. We shall say
that a function $H(x^\mu,u)$ is {\em $\delta$-polyhomogeneous in
$x$ with a uniform zero of order $l$ in $u$} if $H$ is smooth in
$u\in\R^N$ at fixed $x^\mu$, if $H$ satisfies (\ref{S2.71}) for
any $0\le i\leq \min\{ l,k\} $ and any $k\in \N$, if
\begin{equation}
  \label{eq:m1}
\forall i\in\N \quad \partial^i_uH(\cdot, u)\in\mca
\end{equation}
at fixed constant $u$, and if we have the uniform estimate for
constant $u$'s
\begin{equation}
  \label{m1}
  \forall \epsilon >0,  M \geq 0 , i,k\in \N \ \; \exists\;
C(\epsilon,M,i,k) \ \forall |u|\leq M \quad \| \partial^i_u
H(\cdot, u)\|_{\mcC^{-\epsilon}_k} \leq C(\epsilon,M,i,k)\;.
\end{equation}
The qualification ``in $u$'' in ``uniform zero of order $l$ in
$u$'' will often be omitted. The small parameter $\epsilon$ has
been introduced above to take into account the possible
logarithmic blow-up of functions in $\mca$ at $x=0$; for the
applications to the nonlinear scalar wave equation or to the wave
map equation on Minkowski space-time, the alternative simpler
requirement would actually suffice:
\begin{equation}
  \label{m2}
  \forall  \; M \geq 0 , i,k\in \N \ \; \exists\;
C(M,i,k) \ \forall |u|\leq M \quad \| \partial^i_u H(\cdot,
u)\|_{\mcC^{0}_k} \leq C(M,i,k)\;,
\end{equation}
again for constant $u$'s. Clearly functions which are jointly
smooth in $u$ and in $x^\mu$ satisfy the above conditions;
Lemma~\ref{lemmapol} below provides another class of such
functions. The following simple facts about functions in the above
class will be useful:

\begin{Lemma} \label{lemmapol} Let $m_1,m_2,k\in\N$, $m_1\leq m_2$, and
  let  $P(x^\mu,u)$ be a
  polynomial in $u= (u^1,\ldots u^N)$ of the form
$$ P(x^\mu,u) = \sum_{m_1\leq j \leq m_2}P_{i_1\ldots i_j}(x^\mu)u^{i_1}
\ldots u^{i_j} \;, $$ with coefficients  $P_{i_1\ldots i_j}(x^\mu)
\in\mcA_\infty^\delta$. Then:
\begin{enumerate}
\item $P$ is $\delta$-polyhomogeneous in $x$ with a
uniform zero of order $m_1$.
\item If
 $$f\in \mcA_k^\delta 
+ \mcC^\lambda_\infty $$ for some
 $\lambda>0$, then for any $\epsilon >0$ we have
$$P(\cdot,x^{q\delta} f)\in x^{m_1q\delta}(
 \mcA_{k}^\delta
+ \mcC^{\lambda-\epsilon }_\infty)\;.$$
\end{enumerate}
\end{Lemma}

The proof of Lemma~\ref{lemmapol} is elementary and will be left
to the reader.

\begin{Lemma}\label{lemmaGpol} Let $k,q\in\N$ and let $H(x^\mu,u)$ be
  $\delta$-polyhomogeneous with respect to $x$ with a zero of order
  $m$ in $u$. If
 $$f\in \cases{
 \mcA_k^\delta\cap L^\infty
+ \mcC^\lambda_\infty\;, & $q=0\;,$ \cr
 \mcA_k^\delta
+ \mcC^\lambda_\infty\;, & otherwise, 
} $$ for some
 $\lambda>0$, then  for any $\epsilon>0$
$$H(\cdot,x^{q\delta} f)\in x^{mq\delta}(
 \mcA_{k}^\delta
+ \mcC^{\lambda-\epsilon }_\infty)\;.$$
\end{Lemma} \proof We Taylor-expand $H$ in $u$ to order $r$, where $r$
is any number satisfying
$$rq\delta > m q \delta + \lambda\;.$$
We then have
$$H(x^\mu,x^{q\delta}f) = P(x^\mu,x^{q\delta}f) +R\;,$$
where $P$ is a polynomial and $R$ is a remainder. We note that the
coefficients of the expansion of $P$ can be obtained by
differentiating with respect to $u$ and setting $u=0$, and are
therefore in $\mca$ by \eq{eq:m1}. Further, the usual integral
formula for the remainder in a Taylor expansion together with
\eq{m1} shows that $R$ has a uniform zero of order $r$, in the
sense of \Eq{S2.71}. The result follows from Lemma~\ref{lemmapol}
and from \lem{lemmaG1}. \qed

We are ready now to pass to the proof of the non-linear analogue
of Theorem~\ref{Tlemme1}:

\begin{Theorem}
  \label{Tlemme1n}
  Let $p\in \Z, q, 1/\delta\in \N, -1<\beta'\in \R$,
$k\in \N\cup\{\infty\}$, and let $$(\varphi,\psi)\in
\mcC_{\infty}^{\beta'}(\Omega_{x_0,T}) \times
  \mcC_{\infty}^{\beta'}(\Omega_{x_0,T})\;, \quad \psi_1\in L^\infty(\Omega_{x_0,T})
  $$
  ($\psi_1$ as in \Eq{psidecomp}), be  a solution of
  (\ref{nequa}) with $G$ of the form
  \eq{m0}, where $H$ is $\delta$-polyhomogeneous in $x$ with a
  uniform zero of order
  \begin{equation}
    \label{eq:ord}
    m > \frac{p- \frac 1 \delta}{q}\;.
  \end{equation}
 Suppose that Equations~\eq{equa1+}-\eq{HLdelta} hold, and that
  \begin{deqarr}&
    B_{11}\in \left(\mcA_k^\delta \cap L^\infty\right)
    (\Omega_{x_0,T})\;,\qquad
B_{12}, B_{22}, B_{21}\in \mcA_k^\delta (\Omega_{x_0,T}) \;,&
\label{H3second}
\\ &a,b \in \mcA_k^\delta (\Omega_{x_0,T})\;,
\qquad\varphi(0)\in \mcA_k^\delta(M_{x_0})\;.& \label{H4second}
  \end{deqarr}
 Then $$\varphi\in
\left(x^{(mq-p)\delta}\mcA_k^\delta+\mcA_k^\delta\right)
(\Omega_{x_0,T}) =
x^{\min((mp-q)\delta, 0)}{\mycal A}_k^{\delta}(\Omega_{x_0,T})\;, $$ 
$$
\psi\in x^{\min\{(mq-p)\delta+1,1\}} \mcA_{k
}^\delta(\Omega_{x_0,T})+
C_\infty(\overline{\Omega_{x_0,T}})\subset \left(\mcA_k^\delta\cap
  L^\infty\right)(\Omega_{x_0,T})\;.$$
If one further assumes
$$L_{12}^\mu,B_{12},a,\varphi(0), G_\varphi(\cdot,0)\in L^\infty(\Omega_{x_0,T})\;,$$
then it also holds that
$$\varphi\in
\left(x^{(mq-p)\delta}\mcA_k^\delta+\mcA_k^\delta\cap
L^\infty\right)(\Omega_{x_0,T})\;.$$
\end{Theorem}

 \remark Obviously the theorem remains true if we replace $G$ by
a finite sum of nonlinearities satisfying the above hypotheses,
with different $p$'s and $q$'s for each term of the sum.

\medskip

\proof The result is established by a repetition of the proof of
Theorem~\ref{Tlemme1}, using Lemma~\ref{lemmaG1} and
\lem{lemmaGpol} to obtain the necessary estimates on the
non-linear terms. We simply note that the condition on the order
$m$ of the non-linearity guarantees, using \lem{lemmaG1}, that
$$\partial_x\psi =c_2\in \mcC^{\lambda-\epsilon}_\infty\;,$$ with $$\lambda=
\min\{\beta',mq\delta-p\delta\} >-1\;,$$ hence $\psi\in L^\infty$
by integration.  Decreasing $\beta'$  if necessary we may without
loss of generality assume that $\beta'=\lambda$. When applying
Lemma~\ref{lemmaGpol} it is convenient to view the function $H$ as
a function of the variable $f:=(\psi_1,x\psi_2,x\varphi)\in
L^\infty$. The remaining details are left to the reader. \qed

As a straightforward corollary of Theorem~\ref{Tlemme1n} one
obtains:

\begin{Theorem}\label{T2phg}  Let $\delta =1$ in odd space dimensions,
  and let $\delta =1/2$ in even space dimensions. Consider Equation
  \eq{SE.1} on
  ${\R}^{n,1}$, $n\geq 2$, with 
  initial data $$\widetilde f|_{\{\tau={0}\}}\;,\quad \partial
  \widetilde f/\partial \tau|_{\{\tau={0}\}}\in \left(\mca\cap
  L^\infty\right)(M_{x_0}) \;.$$ Suppose further that
  $H(x^\mu,f)$ is smooth in $f$ at fixed $x^\mu$, bounded and
  $\delta$-polyhomogeneous in $x^\mu$ at constant $f$, and has a zero of order
  $\ell$ at $f=0$, with $\ell$ as in \eq{condH}.
Then:
\begin{enumerate}
\item There exists $\tau_+ > {0}$ such that $ f$ exists
  $\Omega_{x_0,\tau_+}$, with
  \begin{equation}
    \label{eq:w17}
    \|\tf\|_{L^\infty(\Omega_{x_0,\tau_+})}\;.
  \end{equation}
\item If the initial data are \emph{compatible} polyhomogeneous in the
  sense that there exists $\lambda < 1$ such that
$$ \forall i \in \N \qquad 
\partial_x \partial^i_\tau \tf (0)\in
\mcC^{-\lambda}_{\infty}(M_{x_0})\;, $$
  then the solution is polyhomogeneous on each neighbourhood
  $\Omega_{x_0,\tau_*}$ of
  $\scrip$ on which $f$ exists and satisfies \eq{eq:w17} with $\tau_+$
  replaced by $\tau_*$.
\end{enumerate}
\end{Theorem}

\proof  Point 1 is  Theorem~\ref{T2} specialized to
polyhomogeneous initial data. To prove point 2 we  set \be \psi =
\left(\begin{array}{c}
  \psi_1=\widetilde{f}\\
  \psi_2=\left(\begin{array}{c}
  \phi_-\\
  \phi_A
  \end{array}\right)
  \end{array}\right)\;,
\ee and \be \varphi= \phi_+\;. \ee Then Equation~(\ref{SE.3})
takes the form (\ref{nequa}) with \beqa &G= -\Omega^{-{n+3\over
2}} H(x^\mu, \Omega^{n-1\over 2}\widetilde{f})
\equiv  -\Omega^{-{n+3\over 2}} H(x^\mu, \Omega^{n-1\over 2}\psi_1)\;,&\\
&G_\varphi= -G\;,
\quad G_{\psi_1}=0\;, \quad G_{\psi_2}= \left(\begin{array}{c}
    -G\\
    0\end{array} \right)\;.&
\eeqa For $n$ even we take
  $\delta=1/2$,
  $p=n+3$, $q= n-1$; the condition  (\ref{eq:ord}) then reads $m> {n+1\over
    n-1}$, which coincides with (\ref{condH}).
For $n$ odd we take $\delta=1$,
   $p = {n+3 \over 2}$, $q= {n-1\over 2}$, and
  (\ref{condH}) guarantees again that \eq{eq:ord} holds.
\qed


 \section{Wave maps}
\label{Swave} Let $(\cN,h)$ be a smooth Riemannian manifold, and
let $f:(\cM,\stsg)\to(\cN,h)$ solve the wave map equation. We will
be interested in maps $f$ which have the property that $f$
approaches a constant map $f_0$ as $r$ tends to infinity along
lightlike directions, $f_0(x)=p_0\in\cN$ for all $x\in\cM$.
Introducing normal coordinates around $p_0$ we can write
$f=(f^a)$, $a=1,\ldots,N=\dim \cN$, with the functions $f^a$
satisfying the set of equations
\begin{equation}
  \label{W.1}
\Box_\stsg f^a + \stsg^{\mu\nu}\Gamma^{a}_{bc}(f) \frac{\partial
  f^b}{\partial x^\mu} \frac{\partial f^c}{\partial x^\nu} = 0\;,
\end{equation}
where the $\Gamma^{a}_{bc}$'s are the Christoffel symbols of the
metric $h$. Setting as before $\tf^a=\Omega^{-{n-1 \over 2}}f^a$,
$\tilde\stsg=\Omega^2\stsg$, we then have from \eq{C.3},
\begin{equation}
  \label{W.2}
\Box_{\tilde\stsg} \tf^a =- \Omega^{-{n-1 \over 2}}
\tilde\stsg^{\mu\nu}\Gamma^{a}_{bc}(\Omega^{n-1 \over 2}\tf)
\frac{\partial   (\Omega^{n-1 \over 2}\tf^b)}{\partial x^\mu}
\frac{\partial (\Omega^{n-1 \over 2}\tf^c)}{\partial x^\nu} + {n-1
\over 4n}(\tilde{R} -R\Omega^{-2})\tf^a\;.
\end{equation}
In particular if $(\cM,\stsg)$ is the Minkowski space-time (and if
we use the same conformal transformation as in Section~\ref{ss1})
we obtain a system of Equations \eq{SE.12.1}-\eq{SE.18} with
$a_A=b_A=0$, with the obvious replacements associated with
$\tf\to\tf^a$, and with $G$ in \eq{SE.18} replaced by
 \begin{eqnarray}
\lefteqn{G^a:= - \Gamma^{a}_{bc} (\Omega^{n-1 \over 2}\tf)
\left\{\Omega^{{n-1 \over 2}} (-\phi^b_+\phi^c_- +
\phi_A^b\phi_A^c) \right.} & & \nonumber \\&& -\left.
  (n-1)\Omega^{{n-3 \over 2}}\tf^c
\left[\left(x\phi^b_+ - (1+x+2\tau)\phi^b_-\right)  -
  (n-1) 
\tf^b\right]\right\}\;.
  \label{W.3}\end{eqnarray}

\subsection{Existence of solutions, space derivatives estimates}
As before, for even space-dimensions $n$ the occurrence of
non-integer powers of $\Omega$ above does not allow the use of the
standard conformal method except for special target manifolds
$(\cN,h)$, \emph{cf.\/}~\cite{ChBGu}. This can be handled  in our
approach, and  we show:

\begin{Theorem}\label{T2w}
Consider Equation
  (\ref{W.1}) on ${\R}^{n,1}$ with initial data given on a hyperboloid
  $\hyp\supset \Sigma_{x_0,0}$ in Minkowski space-time, and satisfying
  \beqa\label{cpt1.0w} \tf ^a|_{\Sigma_{x_0,0}}
  \equiv\Omega^{-{n-1\over 2}}f^a|_{\Sigma_{x_0,0}} & \in& \cases{
  \left(\HH_{k+1}^{\alpha}\cap L^\infty\right)(\Sigma_{x_0,0}) \;, &
  $n\geq 3\;,$\cr \left(\HH_{k+1}^{\alpha}\cap
  {{\mcC}}^0_1\right)(\Sigma_{x_0,0}) \;, & $n=2\;,$} \\ \partial_x
  (\Omega^{-{n-1\over 2}}f^a) |_{\Sigma_{x_0,0}} &\in &
  \HH_k^{\alpha}(\Sigma_{x_0,0})\;, \label{cpt1.1w}\\
  \partial_\tau(\Omega^{-{n-1\over 2}}f^a) |_{\Sigma_{x_0,0}} &\in&
  \cases{ \HH_{k}^{\alpha}(\Sigma_{x_0,0}) \;, & $n\geq 3\;,$\cr
  \left(\HH_{k}^{\alpha}\cap L^\infty\right)(\Sigma_{x_0,0}) \;, &
  $n=2\;.$}
\label{cpt1w}
 \eeqa  for some $k > {n\over2} +1$, $ -1< \alpha
\le -1/2$. Then:
\begin{enumerate}
\item There exists $\tau_+>0$ and a solution $ f^a$
  of Equation~\eq{W.1}, defined on a set containing $\Omega_{x_{0},
    \tau_+}$,  satisfying the given initial conditions, such that
\begin{deqarr}\label{wmca}
\lefteqn{\|\widetilde f^a\|_{\mcC^0_1(\Omega_{x_0,\tau_+})}
<\infty \;, \qquad n=2 }&&
\\
\lefteqn{\|xe_+(\widetilde
  f^a)\|_{L^\infty(\Omega_{x_0,\tau_+})}
+\sum_{i=1}^r\|x X_i\widetilde
  f^a\|_{L^\infty(\Omega_{x_0,\tau_+})}
} && \nn\\ && +\|\widetilde f^a\|_{L^\infty(\Omega_{x_0,\tau_+})}
+\|x \partial_\tau \widetilde
  f^a\|_{L^\infty(\Omega_{x_0,\tau_+})}
<\infty \;, \quad n\geq 3 \;. \arrlabel{wmc}
\end{deqarr}
Here the $X_i$'s are the vector fields defined in
Section~\ref{S2}, {\em cf.\/}~\Eq{champscoord}.
\item Further,  if  
$\tau_*$ is  such that $f^a$
  exists on $\Omega_{x_0,\tau_*}$ with \eq{wmc} holding with $\tau_+=\tau_*$,
then for all $0 \leq
  \tau < \tau_*$  we
  have $$\widetilde f^a |_{\Sigma_{x_0,\tau}} \in L^\infty
  (\Sigma_{x_0,\tau})\cap \HH_{k+1}^\alpha(\Sigma_{x_0,\tau})\; ,$$
  $$\partial_\tau \widetilde f^a|_{\Sigma_{x_0,\tau}} \in
  \HH{^\alpha_k} (\Sigma_{x_0,\tau})
\;,\qquad \partial_x \widetilde f^a|_{\Sigma_{x_0,\tau}} \in
  \HH^{\alpha}_k (\Sigma_{x_0,\tau}) \;,$$
uniformly in $\tau$. If $n=2$ we also have uniform bounds in the
following spaces
$$ \widetilde f^a |_{\Sigma_{x_0,\tau}} \in \left({\mcC}^0_1
  \cap \HH_{k+1}^\alpha\right)(\Sigma_{x_0,\tau})\;,
  \qquad \partial_\tau
  \widetilde f^a|_{\Sigma_{x_0,\tau}} \in \left(\HH{^\alpha_k}
   \cap L^{\infty}\right)(\Sigma_{x_0,\tau})\;.$$
\end{enumerate}\end{Theorem}

\remark Integration of condition \eq{cpt1.1w} implies of course
that $\tf  \in L^\infty(\Sigma_{x_0,0})$.

\medskip

\proof The proof is similar to that of Theorem~\ref{T2}, but
simpler, because we do not need to gain a $1/2$ in the decay rate,
as done in Lemma~\ref{LT2}. We write Equation~\eq{W.1} in the form
\eq{SE.11}-\eq{SE.17.1}, with $a_A=b_A=0$ and with $G$ in
\eq{SE.18} replaced by $G^a$ defined in \eq{W.3}. We write $G^a$
as
\begin{equation}
  \label{eq:W4}
  G^a= A^a + B^a + C^a + D^a + E^a\;,
\end{equation}
with the order of terms in \eq{eq:W4} corresponding to that in
\eq{W.3}. Since we are working in normal coordinates,
$\Gamma^a_{bc}$ has a uniform zero of order one in the sense of
\eq{S2.71} at $f^a=0$. We want to use Equation~\eq{L.10} to get an
a-priori estimate for the solutions of \eq{W.1}; for this we shall
need to estimate the $\HH^{\alpha}_k$ norms of all the terms which
occur in \eq{eq:W4}. The simplest such term is $E^a$:
\begin{eqnarray*}
  \|E^a\|_{\HH^{\alpha}_k}& \equiv & (n-1)^2\| \Gamma^{a}_{bc} (\Omega^{n-1 \over 2}\tf)
(\Omega^{{n-1 \over 2}}\tf^c)(\Omega^{n-1 \over 2}\tf^b)
\Omega^{-1-{n-1 \over 2}}\|_{\HH^{\alpha}_k}
\\ & \approx & (n-1)^2 \| \Gamma^{a}_{bc} (\Omega^{n-1 \over 2}\tf)
(\Omega^{{n-1 \over 2}}\tf^c)(\Omega^{n-1 \over 2}\tf^b)
\|_{\HH^{\alpha+(n+1)/2}_k}\;, \end{eqnarray*} where we have used
the fact that $\Omega/x$ is a smooth, and therefore bounded,
function. The function $\Gamma^{a}_{bc} (\Omega^{n-1 \over 2}\tf)
(\Omega^{{n-1 \over 2}}\tf^c)(\Omega^{n-1 \over 2}\tf^b)$ can be
viewed as a smooth function $F$ of $x^{n-1 \over 2}\tf^a$ with a
uniform zero of order three. We can thus apply \eq{S2.7.2} with
$l=3$ to obtain
\begin{eqnarray}
  \|E(s)\|_{\HH^{\alpha}_k}& \leq  & C(\| \tf(s)\|_{L^\infty}) \|
  \tf\|_{\HH^{\alpha+2-n}_k} \nonumber
\\ & \leq  & C(\| \tf(s)\|_{L^\infty}) \|
  \tf\|_{\HH^{\alpha}_k}\;, \label{W5} \end{eqnarray}
since $n\geq 2$. We note that in dimensions larger than or equal
to three we have at least one  power of $x$ ``left unused'' above,
which will be made use of in estimating the remaining
contributions to $G^a$. We proceed in a similar way with the other
terms; in space dimension $n=2$ we view $\Omega^{(n+1)/2}D^a\equiv
\Omega^{(n+1)/2}(n-1)(1+x+2\tau)\Omega^{{n-3
    \over 2}}\Gamma^{a}_{bc} (\Omega^{n-1 \over 2}\tf)\tf^c\phi^b_- $
as a smooth function $F$ with a uniform zero of order three of
$(x^{n-1 \over 2}\tf^a,x^{n-1 \over 2}\phi^a_- )$, which  leads to
the estimate
\begin{eqnarray}
  \|D(s)\|_{\HH^{\alpha}_k}& \leq  & C(\| \tf(s)\|_{L^\infty},
  \|\phi_- (s)\|_{L^\infty}) \left(\|
  \tf\|_{\HH^{\alpha}_k} + \|\phi_- (s)\|_{\HH^{\alpha}_k}\right)\;.
  \label{W6} \end{eqnarray}
On the other hand, in dimension $3$ or higher  we can view
$\Omega^{(n+1)/2}D^a$ as a function $F$ with a uniform  zero of
order three of $(x^{n-1 \over 2}\tf^a,x^{n-1 \over 2}x\phi^a_- )$,
which implies
\begin{eqnarray}
  \|D(s)\|_{\HH^{\alpha}_k}& \leq  & C(\| \tf(s)\|_{L^\infty},
  \|x\phi_- (s)\|_{L^\infty}) \left(\|
  \tf\|_{\HH^{\alpha}_k} + \|x\phi_- (s)\|_{\HH^{\alpha}_k}\right)\;.
  \label{W6.3} \end{eqnarray}
Regardless of dimension we view
$\Omega^{(n+1)/2}C^a\equiv\Omega^{(n+1)/2} (n-1)x\Omega^{{n-3
\over 2}}\Gamma^{a}_{bc} (\Omega^{n-1 \over
  2}\tf)\tf^c\phi^b_+ $ as a smooth function  with a uniform zero
  of
order three of $(x^{n-1 \over 2}\tf^a,x^{n-1 \over 2}x\phi^a_+ )$,
obtaining thus
\begin{eqnarray}
  \|C(s)\|_{\HH^{\alpha}_k}& \leq  & C(\| \tf(s)\|_{L^\infty},
  \|x \phi_+ (s)\|_{L^\infty}) \left(\|
  \tf\|_{\HH^{\alpha}_k} + \|x\phi_+ (s)\|_{\HH^{\alpha}_k}\right)\;.
  \label{W7} \end{eqnarray}
Viewing $B^a$ as a function of $(x^{n-1 \over 2}\tf^a,x^{n-1 \over
  2}x \phi^a_A )$, and viewing $A^a$ as a function of $(x^{n-1 \over
  2}\tf^a,x^{n-1 \over 2}x\phi^a_- ,x^{n-1 \over 2}x\phi^a_+ )$, one
similarly obtains for $n\geq 3$
\begin{eqnarray}
  \nonumber \|A(s)\|_{\HH^{\alpha}_k}& \leq  & C(\| \tf(s)\|_{L^\infty},
  \| x\phi_- (s)\|_{L^\infty}, \|x \phi_+ (s)\|_{L^\infty}) \times
\\ && \left(\|
  \tf\|_{\HH^{\alpha}_k} + \|x\phi_- (s)\|_{\HH^{\alpha}_k}+ \|x\phi_+
  (s)\|_{\HH^{\alpha}_k}\right)
\;,
  \label{W9.3}
\\
  \|B(s)\|_{\HH^{\alpha}_k}& \leq  & C(\| \tf(s)\|_{L^\infty},
  \| x\phi_A (s)\|_{L^\infty}) \left(\|
  \tf\|_{\HH^{\alpha}_k} + \|x\phi_A (s)\|_{\HH^{\alpha}_k}\right)\;,
\label{W8.3}
                             \end{eqnarray}
while in dimension $2$ it holds that
\begin{eqnarray} \nonumber \|A(s)\|_{\HH^{\alpha}_k}& \leq  & C(\| \tf(s)\|_{L^\infty},
  \| \phi_- (s)\|_{L^\infty}, \|x \phi_+ (s)\|_{L^\infty}) \times
\\ && \left(\|
  \tf\|_{\HH^{\alpha}_k} + \|\phi_- (s)\|_{\HH^{\alpha}_k}+ \|x\phi_+
  (s)\|_{\HH^{\alpha}_k}\right)
\;.
  \label{W9}
\\
  \|B(s)\|_{\HH^{\alpha}_k}& \leq  & C(\| \tf(s)\|_{L^\infty},
  \| \phi_A (s)\|_{L^\infty}) \left(\|
  \tf\|_{\HH^{\alpha}_k} + \|\phi_A (s)\|_{\HH^{\alpha}_k}\right)\;,
\label{W8}
\end{eqnarray}
Summarizing, in space dimension two we have obtained
\begin{eqnarray}
  \nonumber \|G(s)\|_{\HH^{\alpha}_k}& \leq  & C(\| \tf(s)\|_{L^\infty},
  \| \phi_- (s)\|_{L^\infty},
  \| \phi_A (s)\|_{L^\infty}, \|x \phi_+ (s)\|_{L^\infty}) \times
\\ && \left(\|
  \tf\|_{\HH^{\alpha}_k} + \|\phi_- (s)\|_{\HH^{\alpha}_k}+ \|x\phi_+
  (s)\|_{\HH^{\alpha}_k}+ \|\phi_A (s)\|_{\HH^{\alpha}_k}\right)
 \nonumber\\ &\leq  & C(\| \tf(s)\|_{L^\infty},
  \| \phi_- (s)\|_{L^\infty},
  \| \phi_A (s)\|_{L^\infty},\|x \phi_+ (s)\|_{L^\infty}) \times
\nonumber \\ &&\phantom{xxx}\sqrt{E_\alpha(s)}\;,
  \label{W10} \end{eqnarray}
where
\begin{eqnarray}\nonumber
  E_\alpha(t)&=&\|\widetilde f(t)\|^2_{\HH^{\alpha}_{k}}+ \|\phi_-(t)\|^2
_{\HH^{\alpha}_{k}\Mxonexzerot} %
\\&&
+ \|\phi_+(t)\|^2_{\HH^{\alpha}_{k}\Mxonexzerot} +\sum_{A}\|\phi_A
(t)\|^2_{\HH^{\alpha}_{k}\Mxonexzerot}\;.
  \label{eq:LT2w}
\end{eqnarray}
On the other hand in higher dimensions we can write
\begin{eqnarray}
   \|G(s)\|_{\HH^{\alpha}_k} &\leq  & C(\| \tf(s)\|_{L^\infty},
   \| x\phi_A (s)\|_{L^\infty},\|x\phi_- (s)\|_{L^\infty},\|x \phi_+
(s)\|_{L^\infty})
\times \nn \\
&& \phantom{xxx}\sqrt{E_\alpha(s)}\;.
  \label{W10.3} \end{eqnarray}
To obtain a closed inequality from Equations~\eq{L.10} and
\eq{W10} or \eq{W10.3}, we need to control all the $L_\infty$
norms occurring there. Since $k>n/2+1$, from Equation~\eq{W10} and
the weighted Sobolev embeddings we obtain
\begin{eqnarray}
\|G(s)\|_{{\mcC}^{\alpha}_1}& \leq  & C(\| \tf(s)\|_{L^\infty},
  \| \phi_- (s)\|_{L^\infty},
  \| \phi_A (s)\|_{L^\infty}, E_\alpha(s)) 
\;,
  \label{W11}
\end{eqnarray}
in $n=2$, or --- from \eq{W10.3} ---
\begin{eqnarray}
\|G(s)\|_{{\mcC}^{\alpha}_1}& \leq  & C(\| \tf(s)\|_{L^\infty},
   E_\alpha(s)) 
\;,
  \label{W11.3}
\end{eqnarray}
for $n\geq 3$. The identity
\begin{equation}
\label{xint} \tf^a(\tau,x)=\tf^a(\tau,x_0-2\tau) - \frac 12
\int_{x}^{x_0-2\tau} (\phi^a_--\phi^a_+)(\tau,s)\;ds\end{equation}
yields
\begin{eqnarray}
\nonumber \|\tf(s)\|_{L^\infty}  &\leq & C\left(\sqrt{E_\alpha(0)}
+ \|\phi_-(s)\|_{{\mcC}^\alpha_0} +
\|\phi_+(s)\|_{{\mcC}^\alpha_0}\right)
\\ & \leq &C\left(\sqrt{E_\alpha(0)} +  \sqrt{E_\alpha(s)}\right)
\;.
  \label{W13.3}
\end{eqnarray}
for $n\geq 3$, while if $n= 2$ we use the estimate
\begin{eqnarray}
\nonumber \|\tf(s)\|_{L^\infty} + \|\phi_A(s)\|_{L^\infty} &\leq &
C\left(\sqrt{E_\alpha(0)} + \|\phi_-(s)\|_{{\mcC}^\alpha_1} +
\|\phi_+(s)\|_{{\mcC}^\alpha_1}\right)
\\ & \leq &C\left(\sqrt{E_\alpha(0)} +  \sqrt{E_\alpha(s)}\right)
\;.
  \label{W13}
\end{eqnarray}
In \Eqs{W13.3}{W13}, for notational simplicity  we have estimated
$\tf^a(\tau,x_0-2\tau)$ and its angular derivatives by a multiple
of the initial energy $E_\alpha(0)$; strictly speaking, this
should be some functional of $(E_\alpha(0),\tau^*)$ for $\tau^*$
small enough; then such an estimate holds by standard methods  for
$0\le \tau\leq \tau_*<x_0/2$. Further, such an inequality is
correct if we already have a weighted $L^\infty$ bound as assumed
in point 2. of the theorem. If $n\geq 3$ Equations \eq{L.10}  for
$\alpha<-1/2$ or \eq{L.10n} if $\alpha=-1/2$, \eq{W11.3} and
\eq{W13.3} give
\begin{equation}
  \label{eq:glbd3}
  E_\alpha(\tau) \leq
  CE_\alpha({0}) + \int_{0}^{\tau}
\Phi\left( E_\alpha(s)\right) \;ds\;,
\end{equation}
for some constant $C$, and for a  function $\Phi$ which is bounded
on bounded sets, and we conclude as in the proof of
Theorem~\ref{T2}.

If $n=2$, we note the identity
\begin{equation}
  \label{W14}
  \phi_-(\tau,x)=\phi_-(0,x+2\tau) + \int_0^\tau
  e_+(\phi_-)(\sigma,2(\tau-\sigma)+x) \;d\sigma\;.
\end{equation}
{}From the second of Equations~\eq{SE.14.1} we obtain
$$ |e_+(\phi_-)(s,x)|\leq C\left(\|\phi_-(s)\|_{{\mcC}^\alpha_0} +
\|\phi_A(s)\|_{{\mcC}^\alpha_1} + \|\phi_+(s)\|_{{\mcC}^\alpha_0}
+\|G(s)\|_{{\mcC}^\alpha_0}\right) x^\alpha\;,$$ so that
\begin{eqnarray}\nonumber
  |\phi_-(\tau,x)|&\leq& \|\phi_-(0)\|_{L^\infty} + C \int_0^\tau
  \left(\|\phi_-(\sigma)\|_{{\mcC}^\alpha_0} +
\|\phi_A(\sigma)\|_{{\mcC}^\alpha_1} +
\|\phi_+(\sigma)\|_{{\mcC}^\alpha_0} \right.
\\
&&\left.+\|G(\sigma)\|_{{\mcC}^\alpha_0}\right)\left(2(\tau-\sigma)+x\right)^\alpha
\;d\sigma\;.
  \label{W15}\end{eqnarray}
It follows that
\begin{eqnarray}
  \|\phi_-(\tau)\|_{L^\infty}&\leq& \|\phi_-(0)\|_{L^\infty} + C \int_0^\tau
  \left(\sqrt{E_\alpha(\sigma)}
+\|G(\sigma)\|_{{\mcC}^\alpha_0}\right)(\tau-\sigma)^\alpha
\;d\sigma\;.
  \label{W16}\end{eqnarray}
Let
\begin{eqnarray}F(s)\equiv
\| \tf(s)\|_{L^\infty}+
  \| \phi_- (s)\|_{L^\infty}+
  \| \phi_A (s)\|_{L^\infty} + \sqrt{E_\alpha(s)}\;.
 \label{W17}\end{eqnarray}
It follows from \eq{L.10}, \eq{W13} and \eq{W16} that we have
\begin{equation}
  \label{Fineq}
  F(\tau)\leq CF(0)+\int_0^\tau \Phi(F(\sigma))\left(1+
    (\tau-\sigma)^\alpha\right) \;d\sigma\;,
\end{equation}
where $\Phi$ is a function bounded on bounded sets. We have the
following:
\begin{Lemma}\label{Ll}
  There exists a time $\tau_*$, depending only upon $C$, $F(0)$, and
  the function $\Phi$, such that any positive continuous function $
  F:[0,\tau_+)\to\R$ satisfying the inequality \eq{Fineq} with $\alpha
  > -1$ is bounded from above by $CF(0) + 1$ on
  $[0,\max(\tau_+,\tau_*))$.
\end{Lemma}
\proof Let $$M=\sup_{0\leq x\leq CF(0) + 1}|\Phi(x)|\;; $$ if
$M=0$ the result is obviously true, so assume that $M\neq 0$.
{}From Equation~\eq{Fineq} we obtain that on any interval
$[0,\tau)$ on which $F\leq CF(0) + 1$ we have
$$ F(\tau)\leq CF(0)+\int_0^\tau M \left(1+
    (\tau-\sigma)^\alpha\right) \;d\sigma = CF(0)+ M
\left(\tau + {
    \tau^{\alpha+1}\over \alpha + 1}\right)\;.$$
(Equation~\eq{Fineq} with $\tau=0$ shows that $CF(0)\geq F(0)$,
and continuity of $F$ implies that the set of such intervals is
non-empty.) The result is established by choosing
$$\tau_*=\min\left( {
  1\over 2M },\left[{\alpha+1\over 2M }\right]^{1/(\alpha+1)}\right)\;.$$
\qed

Because the existence time $\tau_*$ in Theorem~\ref{T2w} does not
depend upon $x_1$, Theorem~\ref{T2w} with $n=2$ follows again by
an argument identical to the one given at the end of
Theorem~\ref{T2}. \qed

As in the case of the nonlinear wave equation~\eq{SE.1}, in order
to obtain time derivative estimates  we shall need a more general
version of Theorem~\ref{T2w}. Thus, we consider systems of the
form \eq{td1}-\eq{td4P} with a rather more general form of the
non-linearity $G$ appearing there. It should be clear from the
proof of Theorem~\ref{T2w} that it is convenient to treat the case
$n=2$ separately, this will be considered in Section~\ref{S2dtd}
below.  We thus start with a result which holds in dimensions
$n\geq 3$; the same proof gives similar results in dimension $n=2$
for equations with a nonlinearity of higher order:

\begin{Theorem}\label{T2wg} Let  $n\geq 3$ and consider
the system \eq{td1}-\eq{td4} with
\begin{eqnarray}
\lefteqn{ \|a(\tau)\|_{\HH ^\alpha_k} +\|b(\tau)\|_{\HH
^\alpha_k}+ \sup_{a,b=1,2}\|B_{ab}(\tau)\|_{\mcC^0_k} } && \nn \\
&& +\|B_{0}(\tau)\|_{\mcC^0_k} +\|B_{0}^{-1}(\tau)\|_{L^\infty}
+\|B_{1}(\tau)\|_{\mcC^0_k}
 \leq \tilde C\;,
\label{gtd5}\end{eqnarray} for some constant $\tilde C$, with the
nonlinearity $G$ in \Eq{td1a}  of the form
\begin{equation}
\label{gwtd1} G = x^{-(n+3)/2} H(x^\mu,x^{(n-1)/2}\tf,
x^{(n-1)/2}x\phi_A,x^{(n-1)/2}x\phi_+,x^{(n-1)/2}x\phi_-)\;,
\end{equation}
with $G_{e_-(\phi_A)}=0$ ({\em cf.\/} \Eq{td4G}), and with $H$
having a uniform zero of order $\ell\geq 3$ in the sense of
(\ref{S2.71}).  Suppose that the initial data satisfy
\beqa\label{gcpt1.0wg} \tf ^a|_{\Sigma_{x_0,0}}
  \equiv\Omega^{-{n-1\over 2}}f^a|_{\Sigma_{x_0,0}} & \in&
  \left(\HH_{k+1}^{\alpha}\cap L^\infty\right)(\Sigma_{x_0,0}) \;,  \\
\partial_x \tf^a |_{\Sigma_{x_0,0}} &\in &
  \HH_k^{\alpha}(\Sigma_{x_0,0})\;, \label{gcpt1.1w}\\
  \partial_\tau\tf^a |_{\Sigma_{x_0,0}} &\in& \HH_{k}^{\alpha}(\Sigma_{x_0,0}) \;,
\label{gcpt1w}
 \eeqa with
some $k > {n\over2} +1$, $ -1< \alpha \le -1/2$, then:
\begin{enumerate}
\item There exists $\tau_+>0$, depending only upon the constant
$\tilde C $ in \eq{gtd5} and a bound on the norms of the initial
data in the spaces appearing in \Eqs{gcpt1.0wg}{gcpt1w}, and a
solution $ f^a$
  of Equations~\eq{td1}-\eq{td4}, defined on a set containing $\Omega_{x_{0},
    \tau_+}$,  satisfying the given initial conditions, such that
\begin{eqnarray}
\lefteqn{\|xe_+(\widetilde
  f^a)\|_{L^\infty(\Omega_{x_0,\tau_+})}
+\sum_{i=1}^r\|x X_i\widetilde
  f^a\|_{L^\infty(\Omega_{x_0,\tau_+})}
} && \nn\\ && +\|\widetilde f^a\|_{L^\infty(\Omega_{x_0,\tau_+})}
+\|x \partial_\tau \widetilde
  f^a\|_{L^\infty(\Omega_{x_0,\tau_+})}
<\infty \;.
\label{gwmc}
\end{eqnarray}
\item Further,  if  
$\tau_*$ is  such that $f^a$
  exists on $\Omega_{x_0,\tau_*}$ with \eq{gwmc} holding with $\tau_+=\tau_*$,
then for all $0 \leq
  \tau < \tau_*$  we
  have \begin{deqarr}
&\widetilde f^a |_{\Sigma_{x_0,\tau}} \in L^\infty
  (\Sigma_{x_0,\tau})\cap \HH_{k+1}^\alpha(\Sigma_{x_0,\tau})\; ,
&\\ &
\partial_\tau \widetilde f^a|_{\Sigma_{x_0,\tau}} \in
  \HH{^\alpha_k} (\Sigma_{x_0,\tau})
\;, & \\ &
\partial_x \widetilde f^a|_{\Sigma_{x_0,\tau}} \in
  \HH^{\alpha}_k (\Sigma_{x_0,\tau}) \;,
\arrlabel{7.37}
       \end{deqarr}
with uniform bounds in $\tau$; this implies
\begin{eqnarray}
{\|x\partial_\tau \phi_+\|_{L^\infty(\Omega_{x_0,\tau_*})}
+
  \|x\partial_\tau \phi_A\|_{L^\infty(\Omega_{x_0,\tau_*})}
+\| (x+2\tau)\partial_\tau \widetilde
f^a\|_{L^\infty(\Omega_{x_0,\tau_*})} <\infty }
\;.\nn\\
\label{gwmc2}
\end{eqnarray}
If $k>n/2+2$ then we also have
\begin{eqnarray}
\|x (x+2\tau)\partial_{\tau}
\phi_-\|_{L^\infty(\Omega_{x_0,\tau_*})} <\infty \;.
\label{gwmc2a}
\end{eqnarray}
\end{enumerate}
\end{Theorem}

\proof The transition from Theorem~\ref{T2w} to Theorem~\ref{T2wg}
is rather similar to that from Theorem~\ref{T2} to
Theorem~\ref{T2g}. We note that the estimates done in the course
of the proof of Theorem~\ref{T2w}, with $n\geq 3$ there, can be
summed up in the inequality \be\| x^{-(n+1)/2}
H(x^\mu,x^{(n-1)/2}\hat f)\|_{\HH ^\alpha_k} \leq C(\|\hat
f\|_{L^\infty}) \|\hat f\|_{\HH ^\alpha_k}\;, \label{gw4}\ee where
$$ \hat f:= (\tf, x\phi_A,x\phi_+,x\phi_-)\;.$$ The minor
modifications of the proof of Theorem~\ref{T2w} needed to obtain
\eq{7.37} and the estimate \eq{gwmc2} on $(x+2\tau)\partial_\tau
\widetilde f$ are identical to the ones described in the proof of
Theorem~\ref{T2g}.  The estimate on $\|x\partial_\tau
\phi_+\|_{L^\infty(\Omega_{x_0,\tau_*})}$ is obtained directly
from \Eq{td8.01} and from \eq{gw4}. The estimate on
$\|x\partial_\tau \phi_A\|_{L^\infty(\Omega_{x_0,\tau_*})}$ is
obtained from the \eq{td1a}--equivalent of the first of
Equations~\eq{SE.15.1}. Next, for $k> n/2 +2$ \Eqsone{td8.0} and
\eq{gw4} give \be e_+(\phi_-) \in \HH ^{\alpha}_{k-1} \subset
\mcC^\alpha_1\;,\label{tdw6}\ee Differentiating \Eq{W14} with
respect to $x$ gives
\begin{equation}
  \label{W14m}
  \partial_x\phi_-(\tau,x)=\partial_x\phi_-(0,x+2\tau) + \int_0^\tau
  \left(\partial_xe_+(\phi_-)\right)(\sigma,2(\tau-\sigma)+x) \;d\sigma\;,
\end{equation}
which together with \eq{tdw6} implies, by straightforward
integration,
\begin{equation}
   \label{W14n}
  x(x+2\tau)|\partial_x\phi_-(\tau,x)| \leq C
\end{equation}
This,  \eq{tdw6}, and the identity
$$\partial_\tau \phi_- = (\partial_\tau-2\partial_x+2\partial_x)
\phi_-= e_+(\phi_-) + 2\partial_x\phi_-$$ establish \eq{gwmc2a}.
\qed

\subsection{Estimates on the time derivatives of the solutions, $n \ge 3$}
\label{Stdwave}

To control the time derivatives of the solutions, as in
Section~\ref{Stdphi} we introduce an index $m$ which counts the
number of corner conditions which are eventually satisfied by the
initial data at the ``corner'' $\tau=x=0$. As before we make a
formal statement only for solutions of the wave-map equation
\eq{W.1}, it should be clear from the proof that an analogous
statement holds for solutions of \eq{td1}-\eq{td4} under
appropriate conditions on the coefficients there.

\begin{Theorem}\label{T2wt} In dimension $n\ge 3$ let $\N\ni m\geq 0$.  Consider  a solution
$f:\Omega_{x_0,\tau_*}\to \R$ of Equation (\ref{W.1}) satisfying
\begin{eqnarray}
\lefteqn{\|xe_+(\widetilde
  f^a)\|_{L^\infty(\Omega_{x_0,\tau_*})}
+\sum_{i=1}^r\|x X_i\widetilde
  f^a\|_{L^\infty(\Omega_{x_0,\tau_*})}
} && \nn\\ && +\|\widetilde f^a\|_{L^\infty(\Omega_{x_0,\tau_*})}
+\|x \partial_\tau \widetilde
  f^a\|_{L^\infty(\Omega_{x_0,\tau_*})}
<\infty \;, \label{gwtmc}
\end{eqnarray}
and suppose that \beqa\label{gcpt1.0wgt} 0\leq i \leq m+1\qquad
\partial^i_\tau\tf ^a|_{\Sigma_{x_0,0}}
   & \in& \HH_{k+m+1-i}^{\alpha}(\Sigma_{x_0,0}) \;,  \\
0\leq i \leq m \qquad
\partial_x \partial^i_\tau\tf^a |_{\Sigma_{x_0,0}} &\in &
  \HH_{k+m-i}^{\alpha}(\Sigma_{x_0,0})\;, \label{gcpt1.1wt}
 \eeqa with
some $k > {n\over2} +2$, $ -1< \alpha \le-1/2$.   Then for $0 \leq
\tau < \tau_*$ and for $0\leq i \leq m$, we have
\begin{deqarr}
\lefteqn{ 0\leq j+i < k+m-n/2} &&\nn\\&&
[(\tau+2x)\partial_\tau]^j\partial_\tau^i\widetilde f^a
|_{\Sigma_{x_0,\tau}} \in L^\infty
  (\Sigma_{x_0,\tau})\cap \HH_{k+m+1-i-j}^\alpha(\Sigma_{x_0,\tau})\;
,
\\
\lefteqn{ 0\leq j+i < k+m-n/2-1}  &&\nn\\&& \partial_x
[(\tau+2x)\partial_\tau]^j\partial_\tau^{i}\widetilde
f^a|_{\Sigma_{x_0,\tau}} \in
  \HH^{\alpha}_{k+m-i-j} (\Sigma_{x_0,\tau}) \;,
\arrlabel{td11wt}\end{deqarr} and
\begin{equation}0\leq p < k-n/2 \qquad
[(\tau+2x)\partial_\tau]^p\partial_\tau^{m+1} \widetilde
f^a|_{\Sigma_{x_0,\tau}} \in \HH{^\alpha_{k-p}}
  (\Sigma_{x_0,\tau}) 
\;, \label{td12wt}\end{equation} with $\tau$-independent bounds on
the norms.
\end{Theorem}

\proof The proof is an inductive application of
Theorem~\ref{T2wg}, as in the proof of Theorem~\ref{T2t}, and will
be omitted. \qed

\subsection{Estimates on the time derivatives, $n=2$}
\label{S2dtd}

In space-dimension two the following equivalent of
Theorem~\ref{T2wg} holds:
\begin{Theorem}\label{2dT2wg} Let  $n=2 $, consider
the system \eq{td1}-\eq{td4}, suppose that \eq{gtd5} holds for
some constant $\tilde C$, with the nonlinearity $G$ in \Eq{td1a}
of the form
\begin{equation}
\label{2dgwtd1} G = x^{-3/2} H(x^\mu,x^{1/2}\tf,
x^{1/2}\phi_A,x^{1/2}\phi_-,x^{3/2}\phi_+)\;,
\end{equation}
with $G_{e_-(\phi_A)}=0$ ({\em cf.\/} \Eq{td4G}), and with $H$
having a uniform zero of order $\ell\geq 3$ in the sense of
(\ref{S2.71}).  Suppose that the initial data satisfy

  \beqa\label{2dcpt1.0w} \tf ^a|_{\Sigma_{x_0,0}}
  \equiv\Omega^{-{1\over 2}}f^a|_{\Sigma_{x_0,0}} & \in& \left(\HH_{k+1}^{\alpha}\cap
  {{\mcC}}^0_1\right)(\Sigma_{x_0,0}) \;, \\ \partial_x
  (\Omega^{-{1\over 2}}f^a) |_{\Sigma_{x_0,0}} &\in &
  \HH_k^{\alpha}(\Sigma_{x_0,0})\;, \label{2dcpt1.1w}\\
  \partial_\tau(\Omega^{-{1\over 2}}f^a) |_{\Sigma_{x_0,0}} &\in&
  \left(\HH_{k}^{\alpha}\cap L^\infty\right)(\Sigma_{x_0,0}) \;.
\label{2dcpt1w}
 \eeqa  for some $k > 2$, $ -1< \alpha
\le -1/2$.
%
Then:
\begin{enumerate}
\item There exists $\tau_+>0$, depending only upon the constant
$\tilde C $ in \eq{gtd5} and a bound on the norms of the initial
data in the spaces appearing in \Eqs{2dcpt1.0w}{2dcpt1w}, and a
solution $ f^a$
  of Equations~\eq{td1}-\eq{td4}, defined on a set containing $\Omega_{x_{0},
    \tau_+}$,  satisfying the given initial conditions, such that
\begin{eqnarray}
\|\widetilde f^a\|_{\mcC^0_1(\Omega_{x_0,\tau_+})}
<\infty \;.
\label{2dgwmc}
\end{eqnarray}
\item Further,  for any
$\tau_*$  such that $f^a$
  exists on $\Omega_{x_0,\tau_*}$ with \eq{gwmc} holding with $\tau_+=\tau_*$,
 we
  have for all $0 \leq
  \tau < \tau_*$
$$ \widetilde f^a |_{\Sigma_{x_0,\tau}} \in \left({\mcC}^0_1
  \cap \HH_{k+1}^\alpha\right)(\Sigma_{x_0,\tau})\;,
  \qquad \partial_\tau
  \widetilde f^a|_{\Sigma_{x_0,\tau}} \in \left(\HH{^\alpha_k}
   \cap L^{\infty}\right)(\Sigma_{x_0,\tau})\;,$$
   $$\partial_x \widetilde f^a|_{\Sigma_{x_0,\tau}} \in
  \HH^{\alpha}_k (\Sigma_{x_0,\tau}) \;,$$
with bounds  uniform in $\tau$. This implies
\begin{eqnarray}
{\|x\partial_\tau \phi_+\|_{L^\infty(\Omega_{x_0,\tau_*})} +\|
\partial_\tau \widetilde f^a\|_{L^\infty(\Omega_{x_0,\tau_*})}
<\infty } \;. \label{2dgwmc2}
\end{eqnarray}
 If $k>4$ then we also have
\begin{eqnarray}
  \|(x+2\tau)\partial_\tau \phi_A\|_{L^\infty(\Omega_{x_0,\tau_*})} <\infty \;.
\label{2dgwmc2a1}
\end{eqnarray}
If $k>4$ and if $\partial_\tau^2\tf|_{\Sigma_{x_0,\tau}}\in
\HH^{-1}_{k-1}$ then it further holds that
\begin{eqnarray}
\|(x+2\tau)\partial_{\tau}
\phi_-\|_{L^\infty(\Omega_{x_0,\tau_*})}
<\infty \;. \label{2dgwmc2a2}
\end{eqnarray}
\end{enumerate}
\end{Theorem}

\proof The proof of point 1. is essentially the same as that of
Theorem~\ref{T2w}, with the modifications discussed in the proof
of  Theorem~\ref{T2g}. We note that the key estimates \eq{W10} and
\eq{W11} hold in exactly the same form here, similarly for
\Eqs{W17}{Fineq}. The estimate on $\partial_\tau \tf$ in
\eq{2dgwmc2} follows from the definition of the norm in
\eq{2dgwmc}. The estimate on $\|x\partial_\tau
\phi_+\|_{L^\infty(\Omega_{x_0,\tau_*})}$ is obtained directly
from \Eq{td8.01} and from \eq{W10}.  To obtain \eq{2dgwmc2a1} one
needs to prove a bound on $\partial_A\phi_-$. This is obtained by
differentiating \eq{W14} with respect to $v^A$ and using the
already known uniform bound for $G$ in $\cH^{\alpha}_{k}$, so that
$\partial_A G\in
\cH^{\alpha}_{k-1}\subset \mcC^\alpha_0$. 
Finally,
$$e_+((x+2\tau)\partial_\tau \phi_-) = (x+2\tau)\partial_\tau \left(e_+(\phi_-)\right)\;,$$
and integrating as in \eq{W14} one finds
\begin{eqnarray}
 \lefteqn{
(x+2\tau)\partial_\tau\phi_-(\tau,x)=(x+2\tau)\partial_\tau\phi_-(0,x+2\tau)}
&& \nn \\ && + \int_0^\tau \left\{
  (2(\tau-\sigma)+x)\partial_\tau\left(e_+(\phi_-)\right)(\sigma,2(\tau-\sigma)+x)  \right\}\;d\sigma\;.
 \label{W14nn}\end{eqnarray}
 The term at the right-hand-side of the first line of \eq{W14nn}
 is bounded because of the hypothesis on $\partial_\tau^2\tf$.
 Expressing $e_+(\phi_-)$ by the right-hand-side of \eq{td8.0},
 one immediately finds that all the linear terms that arise after
 differentiation with respect to $\tau$ are in
 $\HH^{\alpha-1}_{k-1}$ or better, and therefore
give a finite contribution when integrated upon.
The contribution from the non--linearity $G$ can be rewritten as
$$x^{-1}\left\{H_{\tf}\partial_\tau \tf+
{ H}_{\phi^A}\partial_\tau \phi^A+ { H}_{\phi_-}\partial_\tau
\phi_- +{ H}_{\phi_+}x\partial_\tau \phi_+\right\}\;,$$ with
appropriate functions $H_{*}$ which, by Lemma~\ref{L1t}, all have
a uniform zero of order $l-1\ge 2$ in their arguments. This easily
implies that the coefficients (including the $x^{-1}$ factor) in
front of the $\tau$ derivatives are in $L^\infty$, and since each
of the $\tau$-derivative terms is in  $ \cH^{\alpha-1}_{k-1}$ or
better, the whole  term is in $ \cH^{\alpha-1}_{k-1} \subset
C^{\alpha-1}_0$. This is sufficient to lead to a finite
contribution in \eq{W14nn}, and \eq{2dgwmc2a2} follows. \qed

We finally arrive at the two-dimensional  equivalent of
Theorem~\ref{T2wt}; comments identical  to those made in
Section~\ref{Stdwave} apply here. The main difference  is that in
dimension $2$ we need the $L^\infty$ bound on $\partial _\tau\tf$
to obtain existence, which leads to the compatibility condition
\eq{2dgcpt1.3wgt} on the second $\tau$ derivatives of $\tf$ when
one attempts to iteratively apply Theorem~\ref{2dT2wg}. The proof
is again identical to that of Theorem~\ref{T2t} and will be
omitted. Let us just mention that one easily checks that the
conditions spelled out below guarantee that the initial data for
the inductive system of equations are in the right spaces for the
iterative application of Theorem~\ref{2dT2wg}. Further,
\Eqs{2dgwmc2}{2dgwmc2a2} provide the \emph{a-priori} bounds which
guarantee that the existence time of the solution will not shrink
at each iteration step.

\begin{Theorem}\label{2dT2wt} In space-dimension two let $\N\ni m\geq 0$. Consider  a solution
$f:\Omega_{x_0,\tau_*}\to \R$ of Equation (\ref{W.1}) satisfying
\begin{eqnarray}
\|\widetilde f^a\|_{\mcC^0_1(\Omega_{x_0,\tau_+})}
<\infty \;.
\label{2dgwmct}
\end{eqnarray}
and suppose that \beqa\label{2dgcpt1.0wgt} 0\leq i \leq m\qquad
\partial^i_\tau\tf ^a|_{\Sigma_{x_0,0}}
   & \in& (\HH_{k+m+1-i}^{\alpha}\cap \mcC^0_1)(\Sigma_{x_0,0}) \;,  \\
0\leq i \leq m \qquad
\partial_x \partial^i_\tau\tf^a |_{\Sigma_{x_0,0}} &\in &
  \HH_{k+m-i}^{\alpha}(\Sigma_{x_0,0})\;,
  \\
\partial^{m+1}_\tau\tf ^a|_{\Sigma_{x_0,0}}
   & \in& (\HH_{k}^{\alpha}\cap L^\infty)(\Sigma_{x_0,0}) \;, \label{2dgcpt1.1wt2}
   \\ \partial_\tau^{m+2}\widetilde
  f^a|_{\Sigma_{x_0,0}}&\in& \HH^{-1}_{k-1}(\Sigma_{x_0,0})\;,\label{2dgcpt1.3wgt}
 \eeqa with
some $k > 4$, $ -1< \alpha \le-1/2$.   Then for $0 \leq \tau <
\tau_*$ and for $0\leq i \leq m$, we have
\begin{deqarr}
\lefteqn{ 0\leq j+i < k+m-3} &&\nn\\&&
[(\tau+2x)\partial_\tau]^j\partial_\tau^i\widetilde f^a
|_{\Sigma_{x_0,\tau}} \in
(\HH_{k+m+1-i-j}^\alpha\cap \mcC^0_1)(\Sigma_{x_0,\tau})\; ,
\\
\lefteqn{ 0\leq j+i < k+m-3}  &&\nn\\&& \partial_x
[(\tau+2x)\partial_\tau]^j\partial_\tau^{i}\widetilde
f^a|_{\Sigma_{x_0,\tau}} \in
  \HH^{\alpha}_{k+m-i-j} (\Sigma_{x_0,\tau}) \;,
\arrlabel{td11wtn}\end{deqarr} and
\begin{equation}0\leq p < k-3 \qquad
[(\tau+2x)\partial_\tau]^p\partial_\tau^{m+1} \widetilde
f^a|_{\Sigma_{x_0,\tau}} \in \HH{^\alpha_{k-p}}
  (\Sigma_{x_0,\tau}) 
\;, \label{2dtd12wt}\end{equation} with $\tau$-independent bounds
on the norms.
\end{Theorem}
\qed

\subsection{Polyhomogeneous solutions}
\label{Sphgwave}

We are finally ready to prove polyhomogeneity at $\scri$ of
solutions of the wave map equation:

\begin{Theorem}\label{Twavemap}  Let $\delta =1$ in odd space dimensions,
  and let $\delta =1/2$ in even space dimensions. Consider Equation
  \eq{W.1} on
  ${\R}^{n,1}$, $n\geq 2$, with 
  initial data \begin{eqnarray} \partial^i_\tau \widetilde f^a|_{\{\tau={0}\}}
  \in (\mca\cap
L^\infty)(M_{x_0})\;,\quad i=0,1, & n=
  2\;,
  \\
  \widetilde f^a|_{\{\tau={0}\}}\in (\mca\cap L^\infty)(M_{x_0})\;,\quad
  \partial_\tau
  \widetilde f^a|_{\{\tau={0}\}}\in \mca(M_{x_0})\;, & n\ge
  3\;.
  \end{eqnarray}
Then:
\begin{enumerate}
\item There exists $\tau_+ > {0}$ such that $ f^a$ exists
 on $\Omega_{x_0,\tau_+}$, with
 \begin{deqarr}\label{wmcap}
\lefteqn{\|\widetilde f^a\|_{\mcC^0_1(\Omega_{x_0,\tau_+})}
<\infty \;, \qquad n=2\;, }&&
\\
\lefteqn{\|xe_+(\widetilde
  f^a)\|_{L^\infty(\Omega_{x_0,\tau_+})}
+\sum_{i=1}^r\|x X_i\widetilde
  f^a\|_{L^\infty(\Omega_{x_0,\tau_+})}
} && \nn\\ && +\|\widetilde f^a\|_{L^\infty(\Omega_{x_0,\tau_+})}
+\|x \partial_\tau \widetilde
  f^a\|_{L^\infty(\Omega_{x_0,\tau_+})}
<\infty \;, \quad n\geq 3 \;. \arrlabel{wmcp}
\end{deqarr}
\item If the initial data are \emph{compatible} polyhomogeneous in the
  sense that 
$$ \forall i \in \N \qquad \partial^i_\tau \tf^a (0)\in
L^{\infty}(M_{x_0})\;,
$$
  then the solution is polyhomogeneous on each neighbourhood
  $\Omega_{x_0,\tau_*}$ of
  $\scrip$ on which $f$ exists and satisfies \eq{wmcp} with $\tau_+$
  replaced with $\tau\ast $.
\end{enumerate}
\end{Theorem}



\proof Existence of solutions follows from Theorem~\ref{T2w}.
Theorems~\ref{T2wt} and \ref{2dT2wt} give the time-derivative
estimates which are necessary in Theorem~\ref{Tlemme1n}. In order
to apply that last theorem, we set \be \varphi =
\left(\begin{array}{c}
  \phi_+^c\\
  \phi_A^c
  \end{array}\right)\;,
\ee and \be \psi_1=(\widetilde f^c)\;, \quad \psi_2= (\phi_-^c)\;.
\ee Equation~(\ref{W.2}) takes then the form (\ref{nequa}). As in
Theorem~\ref{T2phg}, for $n\ge 4$ even we take
  $\delta=1/2$,
  $p=n+3$, $q= n-1$;  while for $n\ge 3$ odd we take $\delta=1$,
   $p = {n+3 \over 2}$, $q= {n-1\over 2}$. For $n=2$ we set
  $\delta=1/2$,
  $p=3$, $q= 1$. The non-linearity here has a uniform zero of order $3$,
   which is compatible with the hypotheses of Theorem~\ref{Tlemme1n}, and the result follows
   by that last theorem. \qed


{\bf Acknowledgements:} We are grateful to Helmut Friedrich for
many useful comments on a previous version of this paper.

\appendix

\section {Function spaces, embeddings, inequalities} \label{S2}

Throughout this paper the letter $C$ denotes a constant the exact
value of which is irrelevant for the problem at hand, and which
may vary from line to line.

Let $M$ be a smooth manifold such that
$$\overline{\! M} \equiv M \cup \partial M$$
is a compact manifold\footnote{We use the convention in which
manifolds are always open sets. Thus, a manifold with boundary
does \emph{not} contain its boundary as a point set.} with smooth
boundary $\partial M$.
Throughout this work the symbol $x$ stands for a smooth defining
function for $\partial M$, \emph{i.e.}, a smooth function on
$\overline{\! M}$ such that $\{x=0\}=
\partial M$, with $dx$ nowhere vanishing on $\partial M$. It follows
that there exists $x_0 > 0$ and a compact neighborhood $K$ of
$\partial M$ on which $x$ can be used as a coordinate, with $K$
being diffeomorphic to $[0,x_0] \times \partial M$. For $0 \leq
x_1 < x_2 \leq x_0$ we set
\begin{deqarr}
&M_{x_1}=\{p \in M \;|\;0< x(p) < x_1\}\;,& \label{tpM1}\\
&M_{x_1,x_2}=\{p \in M \;|\; x_1 < x(p) < x_2\}\;,&\\
&\tilde{\! \partial} M_{x_1}=\{p \in M \;|\; x(p)=x_1\} \approx
\partial M\;.\arrlabel{tpM}
\end{deqarr}
In  what follows the symbol $\Omega$ will generally denote one of
the sets $M, M_{x_1}$, or $M_{x_1,x_2}$. Any subset of $
\bM_{x_0}$ can be locally coordinatized by coordinates
$y^i=(x,v^A)$, where the $v^A$'s can be thought of as local
coordinates on $\partial M$. We cover $\pM$ by a finite number of
coordinate charts ${\cO}_i$ so that the sets
$\overline{\Omega}{}_i$, where
 $$\Omega_i:=(0,x_0)\times{\cO}_i\;,$$ cover
$M_{x_0}$. We use the usual multi-index notation for partial
derivatives: for $\beta = (\beta_1,\ldots,\beta_n) \in \N^n$ we
set $\partial^\beta =
\partial^{\beta_1}_1\ldots \partial^{\beta_n}_n$. We will write
$\partial_v^\beta$ for derivatives of the form $
\partial^{\beta_2}_2\ldots \partial^{\beta_n}_n$, which do not involve
the $x^1\equiv x$ variable.

If $\locO $ is an open set, for $k\in \N\cup\infty$ we let
$C_k(\locO )$ denote the usual space of $k$-times differentiable
functions on $\locO $; the symbol $C_k(\ovlocO )$ is used to
denote the set of those functions in $C_k(\locO )$ the derivatives
of which, up to order $k$, extend by continuity to $\ovlocO$. We
emphasize that no uniformity is assumed in $C_k(\locO )$, so that
functions there could grow without bound when approaching the
boundary of $\locO$. Nevertheless, the symbol $\|\cdot\|_{C_k}$
will denote the usual supremum norm of $f$ and its derivatives up
to order $k$.
 The symbol $C_{k+\lambda}(\locO)$ denotes the space of
$k$-times continuously differentiable functions on $\locO$, with
$\lambda$-H{\"o}lder continuous $k$'th derivatives.


For $\alpha \in \R$, $k \in \N$ and $\lambda \in (0, 1]$, we
define ${\mcC}_0 ^ \alpha (\Omega_i)$ (respectively $
{\mcC}^\alpha_{0+\lambda} (\Omega_i)$ , ${\mcC}^\alpha_k
(\Omega_i)$ , ${\mcC}^\alpha_{k+\lambda} (\Omega_i)$) as the
spaces of appropriately differentiable functions such that the
respective norms \beqa \|f\|_{{\mcC}_0 ^ \alpha (\Omega_i)}
&\equiv&
\sup_{p\in \Omega_i} |x^{-\alpha} f(p)|
\,,
\nonumber\\
\|f\|_{{\mcC}^\alpha_{0+\lambda} (\Omega_i)} &\equiv&
\|f\|_{{\mcC}^\alpha_0 (\Omega_i)} + \displaystyle \sup_{y \in
\Omega_i} \,\sup_{y\neq y' \in B(y,{x(y)\over 2}) \cap
  \Omega_i}\,{x(y)^{-\alpha-\lambda}|f(y)-f(y')|\over
  |y-y'|^\lambda}\,,\nonumber\\
\|f\|_{{\mcC}^\alpha_k (\Omega_i)} &\equiv& \displaystyle \sum_{0
  \leq|\beta|\leq k} \ \ \ \|x^{\beta_1}\partial^\beta f\| _{{\mcC}_{0} ^
  {\alpha}(\Omega_i)}\,,
\nonumber\\
\|f\|_{{\mcC}^\alpha_{k+\lambda} (\Omega_i)} &\equiv&
\|f\|_{{\mcC}^\alpha_{k-1}(\Omega_i)} + \displaystyle
\sum_{|\beta|=k} \|x^{\beta_1}\partial^\beta
f\|_{{\mcC}^\alpha_{0+\lambda} (\Omega_i)}\,, \eeqa are finite.
Let $\locO$ be an open subset of $M$, or a submanifold with
boundary in $M$; for such sets we define: \beqa \|f\|_{{\mcC}_k ^
  \alpha (\locO)} & \equiv & \sup_i \| f\|_{{\mcC}_k ^\alpha
  (\Omega_i\cap \locO)} + \|f\|_{C_k
  (\complement M_{x_0/2} \cap\locO)}\;,\nonumber \\
\|f\|_{ {\mcC}_{k+\lambda} ^ \alpha (\locO)} & \equiv & \sup_i
\|f\|_{{\mcC}_{k+\lambda}^\alpha (\Omega_i\cap \locO)} +
\|f\|_{C_{k+\lambda} (\complement M_{x_0/2}\cap \locO)}\;.  \eeqa
 %
We note that $f \in {\mcC}^{\alpha+\sigma}_{k+\lambda} (\Omega)$
if and only if $x^{-\sigma} f \in {\mcC}^\alpha_{k+\lambda}
(\Omega)$.

We define the spaces ${\HH }^\alpha_k (\ooi)$ as the spaces of
those functions in $H_k^{\loc}(\ooi)$ for which the norms
 $\| \cdot \|_{{\HH }^\alpha_k(\ooi)}$ are finite, where
\begin{equation}
\|f\|^2_{{\HH }^\alpha_k (\ooi)} = \displaystyle \sum_{0\leq
|\beta|\leq k} \displaystyle \int_{\ooi} (x^{-\alpha
  +\beta_1}
 \partial^{\beta} f)^2 {dx\over x}\volu\;. \label{S2.0x}
\end{equation}
Here $\dnu$ is a measure on $\partial M$ arising from some smooth
Riemannian metric on $\partial M$.
 This is equivalent to
\begin{equation} \displaystyle
\sum_{0\leq \beta_1+|\beta|\leq k} \displaystyle \int_{\ooi}
(x^{-\alpha} (x\partial_x)^{\beta_1}\partial_v^{\beta} f)^2
{dx\over x}\volu \;,\label{S2.0x.1}\end{equation} and it will
sometimes be convenient to use \eq{S2.0x.1} as the definition of
$\|f\|^2_{{\HH }^\alpha_k (\ooi)}$.  For $\locO $'s such that
$\ooi\subset \locO $ the spaces ${{\HH }^\alpha_k(\locO )}$ are
defined as the spaces of those functions in $H_k^{\loc}(\locO )$
for which the norm squared \be \|f\|^2_{{\HH }^\alpha_k (\locO )}
= \sum_i \|f\|^2 _{\Hak (\ooi)}
\; + \|f\|^2 _{H_k (\locO \cap \complementaire M_{x_0/2})}
\ee is finite.  We note the equivalence of norms,
$$\|f\|_{H_0(\locO )} \approx \|f\|_{{\HH }^{-1/2}_0(\locO )}
\;,$$ and that ${\HH }^{\alpha}_k(M_{x_1,x_2})=H _k(M_{x_1,x_2})$
for all $\alpha$ and $k$ whenever $x_1>0$, the norms being
equivalent, with the constants involved depending upon $x_1$ and
$x_2$, and degenerating in general when $x_1$ tends to zero.

It is often awkward to work with coordinate charts, in order to
avoid that one can proceed as follows: Choose a fixed smooth
complete Riemannian metric $\backg$ on $\bM$. Let $x$ be any
smooth defining function for $\partial M$, we let $X_1$ be the
gradient of $x$ with respect to the metric $\backg$; rescaling
$\backg$ by a smooth function if necessary we may without loss of
generality assume that $X_1$ has length one in the metric $\backg$
in a neighbourhood of $\partial M$.  As before we cover $\pM$ by a
finite number of coordinate charts ${\cO}_i$ with associated
coordinates $v^A$; the $v^A$'s are then propagated to a
neighbourhood of $\partial M$ by requiring $$X_1(v^A)=0\;.$$ This
leads to a covering of $M_{x_0}$ of the kind already used, and one
easily checks that $$X_1=\partial_x$$ in the resulting local
coordinates. This gives then a globally defined vector
$\partial_x$ on $M_{x_0}$. For $i=2,\ldots,r$ we let $X_i$ be any
smooth vector fields on $\partial M$ satisfying the condition that
at any $p\in\partial M$ the linear combinations of the $X_i$
exhaust the tangent space $T_p\partial M$. (If $\partial M$ is a
sphere $S^{n-1}$, a convenient choice is the collection of all
Killing vectors of $(S^{n-1},\mathring h)$, where $\mathring h$ is
the unit round metric on $S^{n-1}$.) Over the domain of a chart
$(v^A)$ of $\partial M$, one thus has
\begin{deqarr} \partial_A &=& \sum_{i=2}^r f_A^i(v^B) X_i\;,\\ X_i&=&
\sum_{A=2}^{n} X_i^A(v^B)\partial_A \arrlabel{champscoord}\;,
\end{deqarr} for some locally defined smooth functions $f_A^i,
X_i^A$; clearly things can be arranged so that those functions are
bounded, together with all their partial derivatives. We propagate
the $X_i$'s to $M_{x_0 }$ by requiring $$[X_1,X_i]=0\;,$$
equivalently \be
\partial_x X_i^A = 0\;. \label{champscoord2}\ee It follows that
\eq{champscoord} still holds with $x$-independent functions. For
any multi-index $\beta=(\beta_1,\beta_2,\ldots,\beta_{r})\in \N^r$
we set,  on $M_{x_0}$, \be\label{decaldef} \decal^\beta f =
X_1^{\beta_1} X_2^{\beta_2}\cdots X_r^{\beta_r}f =
\partial_x^{\beta_1} X_2^{\beta_2}\cdots X_r^{\beta_r}f\;.\ee
It follows that we have \beqan& \|f\|_{{\mcC}^\alpha_k (M_{x_0})}
\approx \displaystyle \sum_{0
  \leq|\beta|\leq k}  \|x^{\beta_1}\decal^\beta f\| _{{\mcC}_{0} ^
  {\alpha}(M_{x_0})}\,,&
\\ &\|f\|^2_{{\HH }^\alpha_k (M_{x_0})} \approx \displaystyle
\sum_{0\leq |\beta|\leq k} \displaystyle \int_{M_{x_0}}
(x^{-\alpha
  +\beta_1}
 \decal^{\beta} f)^2 {dx\over x}\volu&
\eeqan (where $\approx$ denotes the fact that the norms are
equivalent), {\em etc.} Here, $|\beta|=\beta_1+\ldots+\beta_r$.

There is a useful way of rewriting $\|\cdot\|_{{\HH }^\alpha_k
  (M_{x_0})}$ which proceeds as follows: for $f\in{\HH }^\alpha_k
(M_{x_0})$, $s \in (1,2)$, and $n \in \N$ we set
\begin{equation}f_n(s,v) = f(x=x_0{s\over 2^n},
v)\;;\label{S2.1}\end{equation} letting $\approx$ denote again
equivalence of norms one then has, after a change of variables,
\begin{eqnarray}\|f\|^2_{{\HH }^\alpha_k(M_{x_0})} &=&\displaystyle
\sum_{n\geq 1} \displaystyle \sum_{0 \leq |\beta| \leq k}
\int_{[2^{-n}{x_0},2^{1-n}{x_0}]\times \partial M}
|x^{-\alpha+\beta_1} \decal^\beta f (x,v)|^2 {dx\over x}  \dnu
\nonumber \\ & \approx &x{^{-2\alpha}_0}\displaystyle \sum_{n\geq
1} \displaystyle \sum_{0\leq |\beta| \leq k}2^{2 n \alpha}
\int_{[1,2]\times \partial M}|\decal^\beta f_n(s,v)|^2 ds\,
\dnu\nonumber \\ &= &x{^{-2\alpha}_0} \displaystyle \sum_{n\geq 1}
2^{2n\alpha} \|f_n\|^2_{H{_k}([1,2]\times \partial
M)}\;.\label{S2.2}\end{eqnarray} More precisely, we write
$A\approx B$ if there exist constants $C_1,C_2 > 0$ such that $C_1
A \leq B \leq C_2 A$. In $(\ref{S2.2})$ the relevant constants
depend only upon $\alpha$ and $k$.

 It turns out to be useful to
have a formula similar to \eq{S2.2} for functions in ${\HH
}^\alpha_k
  (M_{x_2,x_1})$; this can be done for any $x_1$ and $x_2$, but in
order to obtain uniform control of certain constants it is
convenient to require $2 x_2\leq x_1 $. For such values of $x_1$
and $x_2$ we let $n_0(x_1,x_2)\in \N$ be such that ${x_1 \over
  2^{n_0+1}}\leq x_2 \leq {x_1\over 2^{n_0}}$. For $n\in\N$, $n\geq 1$, and for any $f:
M_{x_2,x_1} \to \R^N$ we then define $f_n: (1,2)\times\partial M
\to \R^N $ by \begin{eqnarray} \nonumber
  n\leq n_0\,, &&f_n(s,v)=f(x_1{s\over 2^n},v)\;,\\
  n= n_0+1\,, && f_n(s,v)= f(x_2\, s , v)\;,\nonumber\\
  n>n_0+1\,, && f_n = 0 \;.
\label{coco25}\end{eqnarray} (This coincides with the definition
already given for $M_{x_1}$, when this set is thought of as being
an ``$M_{x_2,x_1}$ with $ x_2=0$'', if we set $n_0 = +\infty$.)  A
calculation as in \eq{S2.2} shows that for any $2 x_2\leq x_1\leq
x_0$, there exist constants $C_1$ and $c_1$, independent of $x_0$,
$x_1$ and $x_2$, such that for all $ f \in \Hak (M_{x_2,x_1}) ,$
\beqa\nonumber
\lefteqn{c_1 x_1^{-2\alpha}\sum_n \{2^{n\alpha} \|f_n\| _{H_k
([1,2]\times \partial
  M)}\}^2 } & \\ &\leq \|f\|^2_{\Hak(M_{x_2,x_1})}\leq C_1
  x_1^{-2\alpha}\sum_n \{2^{n\alpha} \|f_n\| _{H_k ([1,2]\times
  \partial M)}\}^2\;. \label{gequiv} \eeqa Equation~\eq{S2.2} leads
  one to introduce\footnote{The symbol $\cB$ might suggest to the reader
  that we specifically have Besov spaces in mind; this is not the case,
  and we hope that the notation will not lead to confusion.} spaces ${\cB}^\alpha_{k+\lambda}$, that arise
  naturally from weighted Sobolev embeddings, \emph{cf.\/} Equation
  (\ref{S2.5.1}) below: we define
\begin{eqnarray}\|f\|^2_{{\cB}^\alpha_{k+\lambda}(M_{x_0})}  &= &x{^{-2\alpha}_0} \displaystyle
\sum_{n\geq 1} 2^{2n\alpha} \|f_n\|^2_{C{_{k+\lambda}}([1,2]\times
\partial M)}\;,\label{S2.2.1}\end{eqnarray} $f_n$ as in \eq{S2.1},
and we set $${\cB}^\alpha_{k+\lambda} (M_{x_0}) = \{f \in {
C}_{k+\lambda}(M_{x_0}) \;| \; \|f\|_{{\mycal
B}^\alpha_{k+\lambda}(M_{x_0})} < \infty\}\;.$$ Clearly $$ {\mycal
B}^\alpha_{k+\lambda} (M_{x_0})\subset {\mycal
C}^\alpha_{k+\lambda}(M_{x_0})\;.$$ Since the general term $f_N$,
as well as sums of the form $\Sigma_{n\geq N}f_n$, of a convergent
series tend to zero as $N$ tends to infinity, for $f\in {\cal
B}^\alpha_{k+\lambda} (M_{x_0})$ we actually have \begin{equation}
\label{S2.2.2} \displaystyle \lim_{x_1 \to 0}
\|f\|_{{\mcC}^\alpha_{k+\lambda} ( M_{x_1})}=0\;.
  \end{equation}
  We have the trivial inclusion,
\begin{eqnarray}
\label{S2.3} &\alpha' > \alpha \quad \Longrightarrow \quad
{\mcC}^{\alpha'}_{k+\lambda}(M_{x_1}) \subset {\HH
}^\alpha_k(M_{x_1}) \;. &
\end{eqnarray}
 The fact that the inequality $\alpha^\prime > \alpha$
in \eq{S2.3} is strict has various annoying consequences, which
are best avoided by introducing yet another space --- the space
$\Gak$ of functions in $H^k_{\loc}(M_{x_0})$ for which the norm
squared \be \|f\|^2_{\Gak (M_{x_0})} = \sup_{n\geq 1} \;\left\{
\sum_{0\leq
    \beta\leq k} \int_{[2^{-n}x_0,2^{1-n}x_0]\times\partial M}
  |x^{-\alpha+\beta_1}\decal^{\beta} f(x,v)|^2
{\dx\over x}\volu\,\right\} 
\ee is finite. We note that $\|f\|_{\Gak (M_{x_0})}$ is equivalent
to
\be x_0^{-\alpha}\sup_{n\geq 1}\;\left\{ 2^{n\alpha}\| f_n\| _{H_k
    ([1,2]\times \pM)}\right\} \;, \ee with $f_n(s,v) = f({x_0s\over
  2^n},v)$, as in \eq{S2.1}.  To define the $\Gak(M_{x_2,x_1})$'s,
assuming again that $x_2\leq x_1/2$, we let $I_n(x_1,x_2)$ be
defined as \beqa \nonumber n\leq n_0\;,&& I_n = (2^{-n}x_1,
2^{1-n}x_1)\;,\\\nonumber n=n_0+1\;,&& I_{n_0+1} = (x_2, 2x_2)
\;,\\\label{coco26} n>n_0+1 \;,&& I_n =\emptyset\;, \eeqa where
$n_0$ is as in \eq{coco25}.  For all $f\in H_k^{\loc}
(M_{x_2,x_1})$ we set \be \|f\|^2_{\Gak(M_{x_2,x_1})} = \sup_n
\{\sum_i \sum_{0\leq |\beta|\leq
  k}\int_{\ooi \cap\{ I_n\times \partial M\}} (x^{-\alpha +\beta_1}
\decal^{\beta} f)^2 {dx\over x}\volu\;\} \ee (we identify
$(a,b)\times\partial M$ and $M_{a,b}$). Similarly to \eq{gequiv},
there exist constants $c_2$ and $C_2$, which do \emph{not} depend
upon $x_0$, $x_1,$ and $x_2$, such that for all $2 x_2 \leq
x_1\leq x_0$, \be \label{Gequiv} c_2 x_1^{-\alpha}\sup_n
\|f_n\|_{H_k ([1,2]\times\partial M)} \leq
\|f\|_{\Gak(M_{x_2,x_1})} \leq C_2 x_1^{-\alpha}\sup_n
\|f_n\|_{H_k
  ([1,2]\times\partial M)} \;.\ee
We have the obvious inequality \be \|f\|_{\Gak (\Omega)} \leq
\|f\|_{\Hak (\Omega)}\;, \ee together with the modified version of
\eq{S2.3},
\begin{eqnarray}
\label{S2.3.1} &\alpha' \geq \alpha \quad \Longrightarrow \quad
{\mcC}^{\alpha'}_{k+\lambda} \subset {\cG}^\alpha_k \;; &
\end{eqnarray}
in particular the function $(x,v)\to x^\alpha$ is in
$\Gak(M_{x_0})$.

If $S_k$ denotes a space of functions, where $k \in \N$ is a
differentiability index, we set
$$S_\infty \equiv \cap _{k\in\N}S_k\;,$$
\emph{e.g.}, ${\cG}^\alpha_\infty\equiv \cap_{k\in\N}\Gak$,
\emph{etc}.

 We note the following:
\begin{Proposition} \label{PS1.1} Let \,$\Omega = M$, or \,$\Omega =
M_{x_1}$, $0<x_1\leq x_0$, or \,$\Omega =M_{ x_2,x_1}$, $2x_2 <
x_1 \leq x_0$, and let $ {\HH }^\alpha_k= {\HH
}^\alpha_k(\Omega)$, \emph{etc}. For $k'\in\N$, $\lambda\in
[0,1]$, $0 \le k'+\lambda\le k - n/2 \not\in \N $ or $0 \le
k'+\lambda< k  - n/2 \in \N $ we have the continuous embeddings
\begin{eqnarray}
  {\HH }^\alpha_{k}
\subset {\cB}^\alpha_{k'+\lambda} \subset
{\mcC}^\alpha_{k'+\lambda} \;,
\label{S2.4}
\qquad
  {\HH }^\alpha_{k}
\subset{\cG}^\alpha_k \subset {\mcC}^\alpha_{k'+\lambda } \;,
\end{eqnarray}
and there exists an $x_2$-independent constant $C$ such that we
have
\begin{eqnarray}
\label{S2.5.1} &\forall f \in {\HH }^\alpha_k \qquad
\|f\|_{{\cB}^\alpha_{k'+\lambda}(\Omega)}\leq C \|f\|_{{\HH
}^\alpha_k(\Omega)}\;,& \\ \label{S2.5.2.1} &\forall f \in
{\cG}^\alpha_k \qquad
\|f\|_{{\mcC}^\alpha_{k'+\lambda}(\Omega)}\leq C \|f\|_{{\cal
G}^\alpha_k(\Omega)}\;.&
\end{eqnarray}
\end{Proposition}

\proof (\ref{S2.5.1})-(\ref{S2.5.2.1}) follow immediately from
(\ref{S2.2}) and \eq{gequiv}, together with the standard Sobolev
embedding; the remaining inclusions in \eq{S2.4} are trivial. \qed

\medskip

All other inequalities involving Sobolev spaces have their
counterpart in the weighted setting; we shall in particular need
various  weighted versions of the Moser inequalities. The reader
should note the different weights for the members of \Eq{S2.7.2}
below --- this shift of weights in this inequality is the key to
our handling of nonlinear equations.

\begin{Proposition} \label{PS1.2} Let \,$\Omega = M$, or \,$\Omega =
M_{x_1}$, $0 < x_1 \leq x_0$, or \,$\Omega = M_{x_2,x_1}$, $2 x_2
< x_1 \leq x_0$, and let ${\HH }^\alpha_k = {\HH }^\alpha_k
(\Omega)$, \emph{etc.}
\begin{enumerate}
\item There exists a constant $C = C (\alpha,\alpha',\beta,k,x_1)$ such that,
for all $f \in {\HH }^{\alpha'}_k \cap {\mcC}^\alpha_0$ and $g \in
{\HH }^\beta_k \cap {\mcC}^{\alpha+\beta-\alpha'}_0$, we have
\be\label{S2.5.2}\|fg\|_{{\HH }^{\alpha+\beta}_k} \leq
C\left(\|f\|_{{\mcC}^{\alpha\phantom{'}}_0} \|g\|_{{\HH }^\beta_k}
+ \|f\|_{{\HH }^{\alpha'}_k} \|g\|_{{\mcC}^{\alpha +
\beta-{\alpha'}}_0}\right)\;.\ee
Further, $\forall \ |\gamma| \leq k$,
\begin{eqnarray}\label{S2.6}\hspace{-0.5cm}
\lefteqn{\|x^{\gamma_1}\decal^\gamma (fg) -
(x^{\gamma_1}\decal^\gamma f) g\|_{{\HH }^{\alpha+\beta}_0} \leq C
\left(\|f\|_{{\mcC}^\alpha_0}
\|g\|_{{\HH }^\beta_k}+\phantom{\sum_{i=2}^r}\right. } \nn \\
&&\left.\|f\|_{{\HH }^{\alpha'}_{k-1}} \left(\|
x\partial_xg\|_{{\mcC}^{\alpha+\beta-{\alpha'}}_0}+\sum_{i=2}^r\|
X_ig\|_{{\mcC}^{\alpha+\beta-{\alpha'}}_0}\right)\right)\;,\end{eqnarray}
where the vector fields $X$ are defined in \Eq{champscoord}.
\item Let $F \in C_k (M \times \R^N)$ be a function such that
for all $B\in \R^+$ there exists a constant $C_1=C_1(B)$ so that,
for all $p \in \R^N$, $|p| \leq B$, we have
$$\left\|F(\cdot,p)\right\|_{{\mcC}^0_k(M_{x_0})} \leq C_1\;.$$ Then
for all $\alpha < 0$, $\beta\in\R$, and $ B \in \R^+$ there exists
a constant $C_2 (B,k,\alpha,\beta,x_1)$ such that for all
$\R^N$-valued functions $f \in {\HH }^{\alpha-\beta}_k(\Omega)$
with $\|x^\beta f\|_{L^ \infty (\Omega)} \leq B$ we have
\be\label{S2.7}\left\|F(\cdot,x^\beta f)\right\|_{{\HH
}^\alpha_k}\leq C_2 (1+\|f\| _{{\HH }^{\alpha-\beta}_k})\;.\ee
\end{enumerate}
Further, if $F$ has a {\em uniform} zero of order $l>0$ at $p=0$,
in the sense that  for all $B\in \R$ there exists a constant
$\hat{C}(B)$ such that for all $|\prho|\le B$ and $0\leq i\leq
\min(k,l)$, \be\label{S2.71}
\left\|\frac{\partial^iF(\cdot,\prho)}{\partial
    \prho^i}\right\|_{{\mcC}_{k-i}^0(M_{x_0})} \leq \hat{C}(B)
|\prho|^{l-i}\;,\ee then for all $\alpha \in \R$, $\beta \geq0$,
there exists a constant $C_3 (\hat{C},l,k,\alpha,\beta,B)$ such
that, for all $f\in {\HH }^{\alpha-l\beta}_k (\Omega)$ with
$\|f\|_{L^\infty(\Omega)} \leq B$, we have \be\label{S2.7.2}
\left\|F(\cdot,x^\beta f)\right\|_{{\HH }^\alpha_k} \leq C_3
\|f\|_{{\HH }^{\alpha-l\beta}_k} \;. \ee
\end{Proposition}

\noindent\textbf{Remark:} The hypothesis \eq{S2.71} will hold if
$F$ is \emph{e.g.\/} a polynomial in $p$ with coefficients of
$p^j$ vanishing for $j<l$, and being functions belonging to
${\mcC}^{0}_k$ for $j\geq l$.

\medskip

\proof We shall give a detailed proof of (\ref{S2.7}) and
\eq{S2.7.2}, the inequalities (\ref{S2.5.2})-(\ref{S2.6}) follow
by an analogous argument using \cite[Volume III, p.~10,
Equations~(3.21)-(3.22)]{Taylor}, \emph{cf.\/} the calculation of
Proposition~\ref{PS1.2.1} below. Let, similarly to (\ref{S2.1}),
$$F_n(s,v)=F \left((x={x_0s\over 2^n},v); ({x_0s\over 2^n})^\beta
f(x={x_0s\over 2^n},v)\right)\;;$$ from Equation (\ref{S2.2}) we
have \be\label{S2.8}\|F(\cdot,x^\beta f)\|^2_{{\HH }^\alpha_k
(M_{x_0})} \approx x^{2\alpha}_0 \displaystyle\sum_{n\geq1} 2^{2 n
\alpha}\|F_n\|^2_{H_{k}([1,2]\times \partial M)}\;.\ee We have the
obvious bound $$\displaystyle \sup_{[1,2]\times \partial M}
\left|\left({x{_0}s\over 2^n}\right)^\beta f\left({x_0s\over
2^n},v\right)\right| \leq \|x^\beta f\|_{L^\infty (M_{x_0})}\leq
M\;.$$ Further the partial derivatives of $(s,v) \to F_n (s,v,p)$
with respect to $s$ and $v$ at $p \in \R^N$ fixed, $|p| \leq M$,
can be bounded by a constant depending only upon $$\displaystyle
\sup_{|p|\leq M} \|F(\cdot,p)\|_{{\mcC}^0_k(M_{x_0})}\;.$$ The
usual Moser inequalities \cite{Taylor}[Volume III, p.~11,
Equation~(3.30)] give
$$\|F_n\|^2_{H{_k}([1,2]\times \partial M)} \leq
C\left(1+2^{-2n\beta}\|f_n\|^2 _{H{_k}([1,2]\times \partial
M)}\right)\;,$$ with $f_n$ as in (\ref{S2.1}), and with a constant
$C$ depending upon $k$ and $M$. Inserting this in (\ref{S2.8}) one
obtains (recall that $\alpha<0$)
\begin{eqnarray}
\|F(\cdot,x^\beta f)\|^2_{{\HH }^\alpha_k (M_{x_0})}&\leq& C
\displaystyle \sum_{n\geq 1} 2 ^{2 n \alpha} (1+2^{-2 n \beta}
\|f_n\|^2_{H{_k}([1,2]\times \partial M)})\nonumber \\
&\leq& C \left(1+ \|f\|_{{\HH }^{\alpha-\beta}_k
(M_{x_0})}\right)\;.\label{S2.9}\end{eqnarray} This establishes
(\ref{S2.7}) for \,$\Omega =M_{x_0}$, and (\ref{S2.7}) with
\,$\Omega =M$ readily follows. The remaining \,$\Omega$'s are
handled in a similar way.

To establish (\ref{S2.7.2}), we note the inequality
$$
\left|{\partial^{|\gamma|+i} F_n (\cdot,\prho) \over \partial
    y^\gamma \partial\prho^i} \right| \leq C|\prho|^{\max
  (l-i,0)}\;,$$
which follows from \eq{S2.71} when ${|\gamma|+i}\leq k$.  Letting
$y$ stand for $(s,v)\in [1,2]\times \partial M$, it then follows
that for $|\sigma|\leq k$ we have \beqan |\partial^\sigma F_n|
&=&\left|\sum_{|\gamma|+|\sigma_1|+\cdots +|\sigma_i|=|\sigma|}
  C(\sigma_1,\ldots, \sigma_i,\beta)\left({x_0\over
      2^n}\right)^{\beta(|\sigma_1|+\cdots +|\sigma_i|)} \right.
\\
&&\left. \qquad\qquad\times{\partial^{|\gamma|+i} F_n \over
\partial
    y^\gamma \partial\prho^i} \partial^{\sigma_1}(s^\beta f_n)
  \cdots\partial^{\sigma_i}(s^\beta f_n)\right|\\
&\leq& 2^{-l\beta n}C \sum_{|\sigma_1|+\cdots+|\sigma_i|\leq
|\sigma|} | \partial^{\sigma_1} (s^\beta
f_n)|\cdots|\partial^{\sigma_i}(s^\beta f_n)|\; .  \eeqan The
usual inequalities \cite[Volume~III, Chapter~13,
Section~3]{Taylor} give
$$
\|F_n\|_{ H_k ([1,2]\times\partial M)}\leq C(k,M)2^{-l\beta n}
\|f_n\|_{H_k([1,2]\times\partial M)}\;,$$ for some constant
$C(k,M)$, and one concludes from \eq{S2.8}, as in (\ref{S2.9}).
\qed

We have the following sharper version of \eq{S2.5.2}-\eq{S2.6}:
\begin{Proposition}\label{PS1.2.1}
  Let \,$\Omega = M$, or \,$\Omega = M_{x_1}$, $0<x_1\leq x_0$, or
  \,$\Omega =M_{ x_2,x_1}$, $2x_2 \leq x_1 \leq x_0$, and let $\Hak = \Hak
  (\Omega)$, {\rm etc}. There exists a constant $C_s =
  C_s(\alpha,\beta,k)$ such that, for all $f\in \Hak \cap
  \boa$ and $g\in \Gbk\cap \cob$ 
we  have
\be \|fg\|_{\HH_k^{\alpha+\beta}}\leq
  C_s(\|f\|_{\boa}\|g\|_{\Gbk} + \|f\|_{\Hak}\|g\|_{\cob})\,,
\label{Mo1}\ee
Moreover it also holds that
\begin{eqnarray}\label{Mo2}\hspace{-0.5cm}
\lefteqn{\forall |\gamma|\leq k\ , \ \|x^{\gamma_1}\decal^\gamma
(fg) - (x^{\gamma_1}\decal^\gamma f) g\|_{{\HH }^{\alpha+\beta}_0}
} \nn \\ &&\leq C \left(\|f\|_{\boa}
\|g\|_{{\GG}^\beta_k}+\|f\|_{{\HH}^{\alpha}_{k-1}} \left(\|
x\partial_xg\|_{{\mcC}^{\beta}_0}+\sum_{i=2}^r\|
X_ig\|_{{\mcC}^{\beta}_0}\right)\right)\;,\end{eqnarray} where the
vector fields $X$ are defined in \Eq{champscoord}.
\end{Proposition}
\remark A useful, though less elegant, inequality related to
\eq{Mo1} is
\be \forall \ |\gamma +\sigma| \leq k \qquad \|x^{\gamma_1}
(\decal^\gamma f) x^{\sigma_1}(\decal^\sigma
g)\|_{\HH_0^{\alpha+\beta}}\leq
  C_s(\|f\|_{\boa}\|g\|_{\Gbk} + \|f\|_{\Hak}\|g\|_{\cob})\,.
\label{Mo1.1}\ee

\proof We will prove (\ref{Mo2}), the proof of (\ref{Mo1}) is
essentially identical.  When $\Omega=M_{x_0}$ we do the rescaling
$f_n(s,v) = f({x_0s\over 2^n},v)$, $g_n(s,v)=g({x_0 s\over
2^n},v)$, we then have, for all $|\gamma|\leq k$,
\beqa\nonumber
\lefteqn{ \|x^{\gamma_1}\decal ^{\gamma}(fg)
  -(x^{\gamma_1}\decal ^{\gamma}f)g\|^2_{\HH_0^{\alpha+\beta}}}&&
\\ & \approx & x_0^{-2(\alpha+\beta)}\sum_{n} 2^{2n(\alpha+\beta)}
\|\decal ^{\gamma}(f_ng_n)
  -(\decal ^{\gamma}f_n)g_n\|^2_{H_0([1,2]\times\partial M)}\nonumber
\\
&\leq&
 Cx_0^{-2(\alpha+\beta)}\sum_n 2^{2n(\alpha+\beta)}
\left(\| f_n\|^2_{L^\infty}\|g_n\|^2_{H_{k}} +
\|f_n\|^2_{H_{k-1}}\|\decal  g_n\|^2_{L^\infty}\right)\nonumber
\\
&\leq& C x_0^{-2(\alpha+\beta)} \left( \left(\sum_n 2^{2n\alpha}\|
    f_n\|^2_{L^\infty} \right) \sup_n  \left(2^{2n\beta}\|
    g_n\|^2_{H_{k}}\right) \right.
\nonumber \\ & & +\left.
 \left(\sum_n 2^{2n\alpha}\|f_n\|^2_{H_{k-1}} \right) \sup_n  \left(
   2^{2n\beta}\|\decal  g_n\|^2_{L^\infty}\right)\right)\nonumber
\\
&\approx& C\left(\|f\|^2_{\boa} \|g\|^2_{\GG_{k}^{\beta}} +
  \|f\|^2_{\HH^\alpha_{k-1}} \|g\|^2_{\mcC^\beta_1}\right)
\nonumber\\
& \leq & C_s\left(\|f\|_{\boa} \|g\|_{\GG_{k}^{\beta}} +
  \|f\|_{\HH^\alpha_{k-1}} \|g\|_{\mcC^\beta_1}\right)^2\,.
\label{Moser}\eeqa (In the third line above we have used the
inequality \cite[Volume III, p.~10, Equation~(3.22)]{Taylor}.) The
case $\Omega=M$ follows immediately from the above; the case
$\Omega=M_{x_2 x_1}$ is treated similarly using
\eq{coco25}-\eq{gequiv} and \eq{coco26}-\eq{Gequiv}.
\qed

Similar results can be proved in weighted H{\"o}lder spaces:

\begin{Lemma}\label{interpolation}  Let \,$\Omega = M$, or \,$\Omega =
  M_{x_1}$, $0<x_1\leq
  x_0$, or \,$\Omega =M_{ x_2,x_1}$, $2x_2 \leq x_1 \leq x_0$, and let
  $\mcCak = \mcCak (\Omega)$. Let $f\in \mcC_k^\alpha\cap \mcC_0^\beta$
  and $g\in \mcC_k^\gamma\cap \mcC_0^\delta$ with
  $\alpha+\delta=\gamma+\beta=\sigma$. Then we have $fg\in
  \mcC_k^\sigma$ and \be \|fg\|_{\mcC_k^\sigma}\leq C_i (
  \|f\|_{\mcC_0^\beta} \|g\|_{\mcC_k^\gamma} + \|g\|_{\mcC_0^\delta}
  \|f\|_{\mcC_k^\alpha})\;,\label{ineqinterpol} \ee
\end{Lemma}
\proof The proof is very similar to that of Propositions
\ref{PS1.2} and \ref{PS1.2.1}.  We use the same conventions as in
(\ref{coco25}), (\ref{coco26}).  We have $ \|fg\|_{\mcC_k^\sigma}
\approx\sup_n 2^{n\sigma}\|f_ng_n\|_{C_k(\omega)}$, where
\begin{equation}\label{defomegan}\omega\equiv[1,2]\times\partial M\;, \ee
similarly for $f$ and $g$.  The interpolation inequality
\cite[Appendix A]{Hormander}  $$\|f_n g_n\|_{C_k(\omega)} \leq
C(\|f_n\|_\infty \|g_n\|_{C_k(\omega)} +
\|g_n\|_\infty\|f_n\|_{C_k(\omega)})$$ leads to the conclusion.
\proofend

We have the following $\mcC_k^\beta$ equivalent of the second part
of Proposition~\ref{PS1.2}, with a similar proof, based on
Lemma~\ref{interpolation}:

\begin{Lemma}\label{lemmaG1}
  Let $F$ be a function satisfying the hypotheses of point 2 of
  Proposition~\ref{PS1.2}, with a uniform zero of order $l$ in $p$ in the
  sense of \Eq{S2.71}.  Then, for any $\epsilon>0$, $\beta\in\R$
and $f\in {\mycal C}_k^{\beta}   \cap L^\infty $ we have
$F(.,x^\epsilon f)\in {\mcC}_k^{\beta+l     \epsilon} $, and there
exists a constant $C$ depending upon $\|f\|_{L^\infty}$ such that
\be
  \|F(.,x^\epsilon f)\|_{\mcC_k^{\beta+l\epsilon}} \leq
  C(\|f\|_{\infty}) \|f\|_{\mcC_k^\beta}\;.  \ee
\end{Lemma}

The space of polyhomogeneous functions $\cAp=\cApM$ is defined as
the set of smooth functions on $\bM$ which have an asymptotic
expansion of the form
\begin{equation}
  \label{eq:3new}
  f \sim \sum_{i=0}^\infty \sum_{j=0}^{N_i} f_{ij} {x^{n_i}} \ln
  ^j x\;,
\end{equation}
for some sequences $n_i,N_i$, with $n_i\nearrow\infty$. The
polyhomogeneous expansions of the introduction are of this form if
$r$ there is replaced by $1/x$; this corresponds to the conformal
transformation of Section~\ref{ss1}, which brings ``null infinity"
to a finite distance. We emphasize that we allow non-integer
values of the $n_i$'s; however, we shall mostly be interested in
rational ones, as those arise naturally in the problem at hand.
Here the symbol $\sim$ stands for ``being asymptotic to'': if the
right-hand-side is truncated at some finite $i$, the remainder
term falls off appropriately faster. Further, the functions
$f_{ij}$ are supposed to be smooth on $\bM$, and the asymptotic
expansions should be preserved under differentiation. It is easily
checked that the space $\cAp$ is independent of the choice of the
function $x$, within the class of defining functions for $\pM$.

\section{ODE's in  weighted spaces}\label{SODEsws}
In our handling of PDE's below we will need ODE estimates to
obtain information about solutions, we thus begin with some
\emph{a priori\/} estimates in weighted spaces for ODE's. While
the results are well-known in principle, and easy to prove, we
present them in detail here because their precise form is
necessary for our arguments later in this work. For a vector $w$
we denote by $\|w\|$ or by $|w|$ the usual Euclidean norm, while
for a matrix $b$ the symbol $\|b\|$ denotes its matrix norm.

\subsection{Solutions of $\partial_\tau\varphi +b \varphi = c$ in
weighted spaces} \label{SODEwc}

Let $\cO$ be an open subset of $\pM$, which might be the whole of
$\pM$, or a coordinate patch of $\pM$ with coordinates $v^A$,
whichever appropriate in the context; we set
\begin{equation}\label{defcu}\cUxx\equiv(x_2,x_1)\times
\cO\times [0,T]\;, \ee
\begin{equation}\label{defSigma}
\Sigx\equiv(x_2,x_1)\times\cO\;,\ee with  $0\leq x_2< x_1$. The
time variable $\tau$ will usually be the last variable, so $\tau$
will run from $[0,T]$ whenever $\cUxx$ is involved. Strictly
speaking, $\cUxx$ should carry an extra $T$ index, but we have not
done that in order not to overburden notation. To avoid
ambiguities we emphasize that the spaces $\mcC_{k}^0
  (\cUxx)$ in the Proposition below are defined as in the previous
  section, with the $v^A$ variables there corresponding here to some
  local coordinates on $\cO$ \emph{together with} the time variable
  $\tau$; the time derivative $\partial_\tau$ should be
  understood as a one-sided one at $\tau=0$ and at $\tau=T$.
\begin{Proposition}\label{Plemme}Let $\alpha\in\R$,  $b\in \mcC_{k}^0
  (\cUxx,\textrm{End}(\R^N))$,
  $c\in \mcC_k^\alpha (\cUxx,\R^N)$, then the unique solution $\varphi$
  of the equation \be \partial_\tau\varphi + b \varphi  = c \;,
\label{ode1} \ee with initial data $\tilde{\varphi} \equiv
\varphi|_{\tau=0}$ $\in \mcC_k^\alpha (\Sigx,\R^N)$ is in $
\mcC_k^\alpha (\cUxx,\R^N)$ with \be
\|\varphi\|_{\mcC_k^\alpha(\cUxx)} \leq C\left(n,N,k,T,x_1,
  \|b\|_{\mcC_k^0(\cUxx)}\right)\left(
  \|\tilde{\varphi}\|_{\mcC_k^\alpha(\Sigx)}+
  \|c\|_{\mcC_k^\alpha(\cUxx)}\right)
.\label{ineqtauCk} \ee We also have the estimates \be \|\varphi
(\tau)\|_{\mcC_0^\alpha(\Sigx)} \leq C e^{\|b\|_\infty \tau }
\left(\|\varphi (0)\|_{\mcC_0^\alpha(\Sigx)} + \int_0^\tau e^{-
    \|b\|_{\infty} s} \|c(s)\|_{\mcC_0^\alpha(\Sigx)} \; ds\right)\;.
\label{ineqtauC0}\ee
\beqa
\nonumber
 \|\varphi(\tau)\|_{{\mcC}^\alpha_k(\Sigx)} &\leq& Ce^{C \|b\|_\infty \tau} \times
\left( \|\varphi(0)\|_{{\mcC}^\alpha_k(\Sigx)} +  \int_0^\tau
e^{-C
|b|_\infty s}\| c(s)\|_{{\mcC}^\alpha_k(\Sigx)} \; ds \right.\\
&& +\left. \int_0^\tau e^{(1-C)|b|_\infty
\;s}\|b(s)\|_{{\mcC}^\alpha_k(\Sigx)}\left(
\|\varphi(0)\|_{{\mcC}^\alpha_0(\Sigx)}
\phantom{\int_0^\tau}\right. \right.\nonumber \\ & &\left.\left.
\phantom{xxxxxxxxx}+ \int_0^se^{-|b|_{\infty}\;t}
\|c(t)\|_{{\mcC}^\alpha_0(\Sigx)} \;dt\;\right) \; ds \;\right)\;.
\label{estimateC1} \eeqa
\end{Proposition}

\remarks

1. Analogous results in $\BB_k^\alpha$ spaces can be proved by
similar arguments.

2. An {\em a-priori\/} estimate in weighted Sobolev spaces for
\eq{ode1} follows from Proposition~\ref{PL.1} below by setting
$E^\mu_-\partial_\mu=\partial_\tau\otimes \textrm{id}$  and
$L\equiv\psi\equiv b\equiv 0$ there.

\medskip

\proof Let $k\in \N^\ast $, and let
$\beta=(\beta_1,\beta_2,\ldots\beta_n)$ be a multi-index with
$|\beta|\leq k$; $\partial^\beta \varphi$ verifies the equation
\beqa
 \partial_\tau \partial^\beta \varphi  &=& - \partial^\beta (b\varphi)
+ \partial^\beta c\;. \label{commutation} \eeqa Let $\epsilon >0$
and set
$$e(.,t,\epsilon) = \left( \epsilon+
\sum_{|\beta|\leq k} x^{2(\beta_1-\alpha)}\langle \partial^\beta
\varphi ,\; \partial^\beta \varphi\rangle \right)^{1/2}\;,$$
$$E(t,\epsilon) = \|e(.,t,\epsilon)\|_{L^\infty(\Sigx)}\;.$$ When
$k=0$ one easily finds $$\partial_\tau  e \leq \|b\|e + |c|\;,$$
and \eq{ineqtauC0} readily follows. For $k>0$ we  have \beqan
\partial_\tau  e &=& {1\over e}  \sum_{|\beta|\leq k}
x^{2(\beta_1-\alpha)}\langle \partial_\tau \partial^\beta
\varphi,\;\partial^\beta \varphi \rangle \;,\\
 & \leq & {1\over e}  \sum_{|\beta|\leq k}
 x^{2(\beta_1-\alpha)}|\partial^\beta(-b\varphi+c)|
 \; |\partial^\beta \varphi| \;,
\\ & \leq & { C(k,n)\over e} (\|b\varphi\|_{{\mcC}^\alpha_k(\Sigx)}
 + \|c\|_{{\mcC}^\alpha_k(\Sigx)}) e\;,\\ & \leq & C(k,n)
 (\|b\varphi\|_{{\mcC}^\alpha_k(\Sigx)} + \|c\|_{{\mcC}^\alpha_k(\Sigx)})\;,
\eeqan where $C(k,n)$ is a constant depending upon $k$ and the
space dimension $n$, and which arises from the inequality $
\sum_{i=1}^p |a_i|\leq\sqrt{p}\sqrt{\sum_i |a_i|^2}$ for any real
sequence $(a_i)$.
  The weighted interpolation inequalities, Lemma~\ref{interpolation},
  imply $$\|b\varphi\|_{{\mcC}^\alpha_k(\Sigx)} \leq
C(\|b\|_{L^\infty(\Sigx)}\|\varphi\|_{{\mcC}^\alpha_k(\Sigx)} +
\|b\|_{{\mcC}^0_k(\Sigx)}
\|\varphi\|_{{\mcC}^\alpha_0(\Sigx)})\;,$$ where  $C$  is a
constant which depends upon $k$, $N$ and $n$. It follows that
\beqan
\partial_\tau  e & \leq & C \left(\|b\|_{L^\infty(\Sigx)} \|\varphi\|_{\mcCak(\Sigx)} +
 \|b\|_{{\mcC}^0_k(\Sigx)} \|\varphi\|_{{\mcC}^\alpha_0(\Sigx)} +
\|c\|_{\mcCak(\Sigx)}\right)\\
& \leq &C \left(\|b\|_\infty E(\epsilon,t) +
\|b\|_{{\mcC}^0_k(\Sigx)} \|\varphi\|_{{\mcC}^\alpha_0(\Sigx)} +
\|c\|_{\mcCak(\Sigx)}\right)\;, \eeqan with perhaps a different
constant $C$. By integration we obtain \beqan e(t) &\leq& e(0) + C
\int_0^t\left(\|b\|_\infty E(s,\epsilon) +\|b(s)\|_{\mcCzk(\Sigx)}
\|\varphi(s)\|_{\mcCaz(\Sigx)} + \|c(s)\|_{\mcCak(\Sigx)}\right)ds
\;, \eeqan from which we deduce \beqan E( t,\epsilon) \leq
E(0,\epsilon)+  C\int_0^t\left( \|b\|_\infty E(s,\epsilon) +
\|b(s)\|_{{\mcC}^0_k(\Sigx)}
\|\varphi(s)\|_{{\mcC}^\alpha_0(\Sigx)} +
\|c(s)\|_{\mcCak(\Sigx)}\right) ds \;. \eeqan Using Gronwall's
Lemma and  letting $\epsilon\to 0$ one obtains \beqan E(t,0)&\leq&
e^{C\|b\|_{\infty}t}  E(0,0)
\\ & & + C\int_0^t
e^{C\|b\|_{\infty}(t-s)} \Big(\|b(s)\|_{\mcCzk(\Sigx)}
\|\varphi(s)\|_{\mcCaz(\Sigx)} + \|c(s)\|_{\mcCak(\Sigx)}\Big)ds
\;. \eeqan The estimate~\eq{ineqtauC0} for
$\|\varphi\|_{\mcCaz(\Sigx)}$ inserted in the last inequality
leads to Equation~\eq{estimateC1}. The time-derivative estimates
follow immediately from the above and from the equation satisfied
by $\varphi$. \proofend

\subsection{Solutions of $\partial_x \phi + b \phi =c$ in weighted spaces}
\label{Sodex} All the results in this section, as well as in
Section~\ref{Sphgx} below, remain valid if we replace the set
$\cUxx $ defined in \Eq{defcu} with $\Sigx $ defined in
\eq{defSigma}
--- the time dimension does not play a preferred role in the
current problem. We start with the following elementary result;
the point is to ensure that the relevant constants are $x_2$
independent:
\begin{Lemma} \label{integrationx}
Let $g\in \mcC_k^\alpha(\cUxx ,\R^N)$, $0\le x_2< x_1$, then $f$
defined for $\alpha> -1$ by
$$f(x,v^A,\tau)= \int_{x_2}^x g(s,v^A,\tau) \; ds$$ is in
$\mcC_k^{\alpha+1}(\cUxx , \R^N)$, with $$
\|f\|_{\mcC_k^{\alpha+1}(\cUxx )} \leq\max\left\{1 , {1\over
\alpha+1}\right\} \|g\|_{\mcC_k^\alpha(\cUxx )}\;.$$ Similarly
$f_2$ defined by $$f_2(x,v,\tau) =- \int_{x}^{x_1} g(s,v,\tau) \;
ds$$ satisfies $$(1+ (\ln x)^2 )^{-1/2} f_2 \in\mcC_k^{0 }(\cUxx
)\mbox{ for $\alpha=-1$}\;,
$$
$$f_2\in\mcC_k^{\min \{ \alpha+1,0\} }(\cUxx )\mbox{  for $\alpha<0$ and $\alpha\neq -1$}\;,$$ with
$$\|f_2\|_{\mcC_k^{\min\{ \alpha+1,0\}}  (\cUxx )}\leq
\max\left\{1,\left|{1\over
1+\alpha}\right|,\left|{x_1^{\alpha+1}\over
1+\alpha}\right|\right\} \|g\|_{\mcC_k^\alpha(\cUxx )}\;. $$
\end{Lemma}
\proof We have the trivial relations \beqan
 \int_{x_2}^x s^\alpha \;ds & \leq & {1\over \alpha+1} x^{\alpha+1} \ \mbox{ for $\alpha> -1$}\; , \\
\int_{x}^{x_1} s^{-1} \; ds & = & \ln x_1 - \ln x \;,
 \eeqan
as well as the commutation rules: \beqan
 \partial_x \int_a^x g \;dx \: &= & g(x) \;,\\
\partial_{v^A}\int_a^x \;g dx &=& \int_a^x \partial_{v^A}g \; dx \;,\\ \partial_\tau \int_a^x
g \; dx &=& \int_a^x \partial_\tau g \;dx \;. \eeqan Note that \be
 \|f\|_{\mcC_k^{\alpha+1}(\cUxx )} =
\|\partial_xf\|_{\mcC_{k-1}^{\alpha}(\cUxx )} + \sum_{0\le
i+|\delta|\leq k} \|\partial_\tau^i \partial_{v^A}^\delta
f\|_{\mcC_0^{\alpha+1}(\cUxx )}\;, \ee
 with
$ \|\partial_xf\|_{\mcC_{k-1}^{\alpha}(\cUxx )} =
\|g\|_{\mcC_{k-1}^\alpha(\cUxx )}$. To estimate
$\partial_\tau^i\partial_{v^A}^\delta f$ one writes \beqan
|\partial_\tau^i \partial_v^\delta f| &\leq&  \int_{x_2}^{x} |\partial_\tau^i\partial_v^\delta g|\;ds \;,\\
&\leq& \int_{x_2}^x \|\partial_\tau^i\partial_v^\delta g\|_{\mcC_0^\alpha} s^{\alpha} \; ds \;,\\
&\leq& {1\over \alpha+1} x^{\alpha+1}
\|\partial_\tau\partial_v^\delta g\|_{\mcC_0^\alpha}\;. \eeqan
 The results for $f_2$ are established in a similar way.  \proofend

We shall use the following notation \be\label{scridef}
\mcS_{x_2}=\{x=x_2\}\;, \ee with the range of the other variables
being in principle clear from the context; this is the equivalent
of  the set $\,\tilde {\!
\partial} M_{x_2}$ of \Eq{tpM} when the set-up described there is assumed.

\begin{Proposition} \label{propositionx1}
  Let $0\le x_2 <x_1$, suppose that $b\in \mcC_{k}^{-\epsilon} (\cUxx
  ,End(\R^N))$, $0\leq \epsilon <1$, $c\in \mcC_k^\alpha
  (\cUxx ,\R^N)$, and let $\phi $ be a solution in $C_k^{\loc}(\cUxx
  )$ of the
  equation \be \partial_x \phi + b \phi = c \;.
  \label{ode2}
  \ee Then the following hold:
\begin{enumerate}
\item If $\alpha <-1$ , then $\phi \in \mcC_k^{\alpha+1}(\cUxx )$
  and we have, for $\alpha+2-\epsilon \neq 0$ and for $x_2\le x_3\le x_1
  $ small enough so that $C(\|b\|_{\mcC_0^{-\epsilon}},x_3) <1$,
  $x_3\ne 0$,
  \be
  \|\phi\|_{\mcC_0^{\alpha+1}({\mcU}_{x_2,x_3})} \leq {1\over 1-
    C(\|b\|_{\mcC_0^{-\epsilon}},x_3)}
(x^{-\alpha-1}_3\|\phi\|_{C_0(\mcS_{x_3})} + {1\over |1+\alpha|}
  \|c\|_{\mcC_0^\alpha({\mcU}_{x_2,x_3})})
  \;,
  \ee
  where
  \be\label{constant}
  C(\|b\|_{\mcC_0^{-\epsilon}},x_3) = {x_3^{1-\epsilon}\over
    |2+\alpha-\epsilon|}\|b\|_{\mcC_0^{-\epsilon}({\mcU}_{x_2,x_3})}\;.\ee
  Moreover, if $x_2\le x_3\le x_1$ is small enough so that
  $C_i \,C(\|b\|_{\mcC_0^{-\epsilon}},x_3) <1$,
  where $C_i$ is the constant in the interpolation inequality
(\ref{ineqinterpol}),
  then
  \beqa
  \|\phi\|_{\mcC_k^{\alpha+1}({\mcU}_{x_2,x_3})} &\leq&
C_\alpha(\|b\|_{\mcC_0^{-\epsilon}},C_i,x_3)
  \left(\|\phi(x_3)\|_{C_k(\mcS_{x_3})}+
\|c\|_{\mcC_k^\alpha({\mcU}_{x_2,x_3})}
  \right.\nonumber \\
  && \left. + \|b\|_{\mcC_k^{-\epsilon}({\mcU} _{x_2,x_3})}(
\|\phi(x_3)\|_{C_0(\mcS_{x_3})} +
    \|c\|_{\mcC_0^\alpha({\mcU}_{x_2,x_3})} )\right)\;,
\nn \\ &&
  \label{ineqx}\eeqa
  with  $C_\alpha(\|b\|_{\mcC_0^{-\epsilon}},C_i,x_3)$ an increasing function in
  the first and third variable.

\item If $\alpha=1$, then $(1+(\ln x)^2)^{-1/2}\phi \in \mcC_k^0(\cUxx )$.

\item If $\alpha>-1$, then $\phi_{x_2} \equiv \lim_{x\to x_2} \phi$ is in
$C_k(\mcS_{x_2})$, with
 \be
\label{diffest}\phi-\phi_{x_2} \in \mcC_k^{1-\epsilon}(\cUxx
)+\mcC_k^{\alpha+1}(\cUxx )\;,\ee $\phi \in\mcC_k^{\alpha+1}(\cUxx
)$ if $\phi_{x_2}=0$, and
  \be
  \|\phi\|_{L^\infty({\mcU}_ {x_2,x_3})} \leq {1\over 1- C'(\|b\|_{\mcC_0^{-\epsilon}},x_3)} \left(
  \|\phi\|_{L^\infty(\mcS_{x_3})} + {x_3^{1+\alpha}\over 1+\alpha}
\|c\|_{\mcC_{0}^\alpha(\mcU_{x_2,x_3})} \right) \label{ineqx1}
 \ee
   for $x_2\le x_3\le x_1$ small enough so that
   $$C'(\|b\|_{\mcC_0^{-\epsilon}},x_3):= {x^{1-\epsilon}_3 \over
   1-\epsilon} \|b\|_{\mcC_0^{-\epsilon}(\mcU_{x_2,x_3})}<1\;.$$
   Moreover
  for $x_3$ small enough  so that $C_i C'(\|b\|_{\mcC_0^{-\epsilon}},x_3)
  <1$ we also have
  \beqa
  \|\phi\|_{\mcC_k^{0}({\mcU}_{x_2,x_3})} &\leq&
C_\alpha'(\|b\|_{\mcC_0^{-\epsilon}},C_i,x_3)
  \left(\|\phi(x_3)\|_{C_k(\mcS_{x_3})}+
\|c\|_{\mcC_k^\alpha({\mcU}_{x_2,x_3})}
  \right.\nonumber \\
  && \left. + \|b\|_{\mcC_k^{-\epsilon}({\mcU} _{x_2,x_3})}(
\|\phi(x_3)\|_{C_0(\mcS_{x_3})} +
    \|c\|_{\mcC_0^\alpha({\mcU}_{x_2,x_3})} )\right)\;,
\nn\\&&
  \label{ineqx2}\eeqa
  with  $C_\alpha'$ an increasing function in its first and third
  argument.
\end{enumerate}
  \end{Proposition}

\remarks  1. The inequalities above are standard when $x_2>0$ and
when the constants are allowed to depend upon $x_2$, regardless of
whether or not $x_3$ can be made small. As already mentioned, the
point here is to make sure that the constants do not blow up as
$x_2$ gets small.

2. In case 2. log-weighted estimates are easily derived; they
will, however, not be needed in what follows.

 \proof 1. For simplicity, we will write $\mcC_k^\delta$ for
$\mcC_k^\delta({\mcU}_{x_3,x_2})$. Let $\phi$ be a (local)
solution of (\ref{ode2}), corresponding to  initial data at
$\{x=x_1\} $ in
$C_k(\mcS _{x_1})$. For $a>0$ set%
\newcommand{\betaone}{\beta_1}
$$e_a(x,v^A, \tau):= (a+ \sum_{|\beta|\leq k}
x^{2\betaone }\langle\partial^\beta \phi |
\partial^\beta \phi\rangle)^{1/2}\;,$$
 and $e:= e_0$.  Let $x_3 \in ]x_2,x_1[\cap ]0,1]$ be such that
${x_3^{1-\epsilon} \over |2+\alpha-\epsilon|}
\|b\|_{\mcC_0^{-\epsilon}} <1$. We have for all $x_2<x\leq x_3$,
\beqa - \partial_x e_a &=& -{1\over e_a}\sum \betaone x^{2\betaone
-1}\langle\partial^\beta\phi | \partial^\beta \phi\rangle
\ \ \mbox{ I } \nn \\ 
 && - {1\over  e_a}\sum_{|\beta|\leq k}
x^{2\betaone }\langle\partial^\beta\partial_x\phi|\partial^\beta
\phi\rangle \ \ \mbox{ II }, \label{main} \eeqa Since $\betaone $
is non-negative we have $-\partial_xe_a (x, v^A, \tau) \leq \mbox{
II }$; further \beqa \mbox{II}&=&{1\over e_a}\sum_{|\beta|\leq k}
x^{2\betaone }
\langle\partial^\beta(b \phi- c) |\partial^\beta \phi \rangle\nonumber\\
 &\leq&{1\over e_a}\sum_{|\beta|\leq k}(
|x^{\betaone }\partial^\beta c |+ |x^{\betaone }\partial^\beta
(b\phi)|)\, | x^{\betaone }\partial^\beta \phi | \nonumber \\
&\leq& \sum_{|\beta|\leq k} |x^{\betaone }\partial^\beta c |+
|x^{\betaone }\partial^\beta (b\phi)| \;. \label{preineq} \eeqa
Clearly \beqan
\sum |x^{\betaone }\partial^\beta c|&=&x^{\alpha} \sum |x^{-\alpha+\betaone }\partial^\beta c |\\
&\leq&  x^{\alpha} \|c\|_{\mcC_k^\alpha}\;,\\
\sum_{|\beta|\leq k} |x^{\betaone } \partial^\beta(b\phi)| &=&
x^{\alpha+1-\epsilon}\sum_{|\beta|\leq k}
 |x^{-\alpha-1+\epsilon+\betaone } \partial^\beta (b\phi)|\\
&\leq& x^{\alpha+1-\epsilon}
\|b\phi\|_{\mcC_k^{\alpha+1-\epsilon}}\;, \eeqan
 which gives
\beqa -\partial_xe_a &\leq & x^{\alpha}\|c\|_{\mcC_k^{\alpha}}
+x^{\alpha+1- \epsilon} \|b\phi\|_{\mcC_k^{\alpha+1-\epsilon}}
\;.\label{ineqbase}
 \eeqa
Consider, first, the case $k=0$;  in this case \eq{ineqbase} reads
$$-\partial_xe_a \leq x^\alpha \|c\|_{\mcC_0^\alpha} + x^{\alpha+1-\epsilon} \|b\|_{\mcC_0^{-\epsilon}}
\|\phi\|_{\mcC_0^{\alpha+1}}\;,$$ which, after integrating over
$[x_3,x]$ and passing to the limit $a\to 0$, gives (recall that
$\alpha<-1$) \beqa
 e(x,v^A,\tau)&\leq& e(x_3,v^A,\tau) + \left(-{x^{\alpha+1} \over
(1+\alpha)} + {x^{\alpha+1}_3 \over (1+\alpha)}\right)
\|c\|_{\mcC_0^\alpha} \nonumber\\
&&+\left({x^{\alpha+2-\epsilon}_3 \over (2+\alpha-\epsilon)} -
{x^{\alpha+2-\epsilon} \over (2+\alpha-\epsilon)}\right)
\|b\|_{\mcC_0^{-\epsilon}}\|\phi\|_{\mcC_0^{\alpha+1}} \nonumber \\
&\leq& \|\phi\|_{C_0(\mcS_{x_3})} +{x^{\alpha+1}
\over |1+\alpha|} \|c\|_{\mcC_0^\alpha} \nonumber\\
&& + \left({x^{\alpha+2-\epsilon}_3 \over (2+\alpha-\epsilon)} -
{x^{\alpha+2-\epsilon} \over (2+\alpha-\epsilon)}\right)
\|b\|_{\mcC_0^{-\epsilon}}\|\phi\|_{\mcC_0^{\alpha+1}}\;.\label{lasteq}
\eeqa Suppose for the moment that $\alpha+2-\epsilon<0$;
Equation~(\ref{lasteq}) yields \be e(x,v^A,\tau) \leq
\|\phi\|_{C_0(\mcS_{x_3})} +{x^{\alpha+1} \over |1+\alpha|}
\|c\|_{\mcC_0^\alpha} + {x^{\alpha+2-\epsilon} \over
|2+\alpha-\epsilon|}
\|b\|_{\mcC_0^{-\epsilon}}\|\phi\|_{\mcC_0^{\alpha+1} }\;, \ee and
since $x^{-1-\alpha} \leq x_3^{-1-\alpha} \leq 1$ we obtain \beqan
 x^{-\alpha-1} e(x,v^A,\tau) &\leq &
x_3^{-1-\alpha} \|\phi\|_{C_0(\mcS_{x_3})} + {1\over |1+\alpha|}
\|c\|_{\mcC_0^\alpha} + {x^{1-\epsilon}_3 \over
|2+\alpha-\epsilon|} \|b\|_{\mcC_0^{-\epsilon}}
\|\phi\|_{\mcC_0^{\alpha+1}}\;. \eeqan On the other hand, if
$\alpha+2-\epsilon>0$ then \beqan e(x,v^A,\tau)&\leq&
\|\phi\|_{C_0(\mcS_{x_3})} +{x^{\alpha+1} \over |1+\alpha|}
\|c\|_{\mcC_0^\alpha} + {x^{\alpha+2-\epsilon}_3 \over
(2+\alpha-\epsilon)})
\|b\|_{\mcC_0^{-\epsilon}}\|\phi\|_{\mcC_0^{\alpha+1} }\;, \eeqan
which gives \beqan x^{-\alpha-1} e(x,v^A,\tau) &\leq &
x_3^{-1-\alpha} \|\phi\|_{C_0(\mcS_{x_3})} + {1\over |1+\alpha|} \|c\|_{\mcC_0^\alpha}\\
&& + {x^{1-\epsilon}_3 \over |2+\alpha-\epsilon|}
\|b\|_{\mcC_0^{-\epsilon}} \|\phi\|_{\mcC_0^{\alpha+1}}\;. \eeqan
The inequality $\|\phi\|_{\mcC_0^{\alpha+1}({\mcU} _{x_2,x_3})}
\leq \sup_{[x_2,x_3]} x^{-1-\alpha}e$ shows  that in all cases we
have
$$
\|\phi\|_{\mcC_0^{\alpha+1}({\mcU}_{x_2,x_3})} \leq {1\over
1-C(\|b\|_{\mcC_{0}^{-\epsilon}},x_3)} ( x_3^{-1-\alpha}
\|\phi\|_{C_0(\mcS_{x_3})} + {1\over |1+\alpha|}
\|c\|_{\mcC_0^\alpha}) \;, $$
with the constant as in \Eq{constant}. Consider, now, any $0< k\in
\N$; Equation~(\ref{ineqbase}) and the  interpolation inequality
(\ref{ineqinterpol}) give
$$ 
 -\partial_xe_a \leq x^\alpha\|c\|_{\mcC_k^\alpha} +
x^{\alpha+1-\epsilon} C_i ( \|b\|_{\mcC_0^{-\epsilon}}
\|\phi\|_{\mcC_k^{\alpha+1}} + \|b\|_{\mcC_k^{-\epsilon}}
\|\phi\|_{\mcC_0^{\alpha+1}})\;.
$$
An argument identical to the one before, considering separately
the cases $\alpha+2-\epsilon
>0$ or $<0$, leads to \beqan \|\phi\|_{\mcC_k^{\alpha+1}} &\leq&
{1\over 1- C_i\,
C(\|b\|_{\mcC_0^{-\epsilon}},x_3)}\left(x_3^{-1-\alpha}
\|\phi\|_{C_k(\mcS_{x_3})}
+ {1\over |1+\alpha| } \|c\|_{\mcC_k^\alpha}\right)\\
& &+ {C_i\over 1- C_i\, C(\|b\|_{\mcC_0^{-\epsilon}},x_3) }
{x_3^{1-\epsilon}\over |2+ \alpha -\epsilon|}
\|b\|_{\mcC_k^{-\epsilon}}
\|\phi\|_{\mcC_0^{\alpha+1}}\\
&\leq &  {C(x_3)\over 1- C_i\, C(\|b\|_{\mcC_0^{-\epsilon}},x_3)}
\left(\|\phi\|_{C_k(\mcS_{x_3})} + \|c\|_{\mcC_k^\alpha}\right. \\
&&\left. + {C_i\over 1-C\|b\|_{\mcC_{0}^{-\epsilon}},x_3)} \,
\|b\|_{\mcC_k^{-\epsilon}}\left( x^{-\alpha-1}_3
\|\phi\|_{C_0(\mcS_{x_3})} + {1\over
  |1+\alpha|} \|c\|_{\mcC_0^{\alpha}}\right)\right)\;,
\eeqan which gives (\ref{ineqx}). We have thus shown that
$\phi\in\mcC_k^{\alpha+1}({\mcU}_{x_2,x_3})$; the property that
$\phi \in \mcC_k^{\alpha+1}(\cUxx )$ immediately follows.

2. The proof is identical, except for a few obvious modifications
in the calculations.

3.  To obtain the $L^\infty$ estimate, we start from
\eq{main}-(\ref{preineq}) with $k=0$, which upon integration and
passing to the limit $a\to 0$ gives \beqan e(x,v^A,\tau)&\leq&
e(x_3,v^A,\tau) + {x^{\alpha+1}_3 \over 1+\alpha}
\|c\|_{\mcC_0^\alpha} +{x^{1-\epsilon}_3 \over 1-\epsilon}
\|b\|_{\mcC_0^{-\epsilon}}\|\phi\|_{\mcC_0^{0}}\;, \eeqan from
which we deduce
$$
 \|\phi\|_{L^\infty({\mcU}_{x_2,x_3})} \leq
\|\phi\|_{L^\infty(\mcS_{x_3})} + {x_3^{\alpha+1}\over \alpha+1}
\|c\|_{\mcC_0^\alpha} +{x^{1-\epsilon}_3 \over 1-\epsilon}
\|b\|_{\mcC_0^{-\epsilon}}
\|\phi\|_{L^\infty({\mcU}_  {x_2,x_3})}\;,
$$
and (\ref{ineqx1}) follows. The proof of (\ref{ineqx2}) is similar
to that of the analogous statement in point 1. From what has been
said it can be seen that $\phi_{x_2} \equiv \lim_{x\to x_2} \phi$
exists and is in $C_k(\mcS_{x_2})$. It remains to show that
  $\phi-\phi_{x_2}$ satisfies \eq{diffest}. When $b$ is a multiple
  of the identity, we can integrate (\ref{ode2}) to obtain \be \label{intrep} \phi(x,\cdot)
= \phi_{x_2}(\cdot) e^{-\int_{x_2}^x b(s,\cdot)ds} + \int_{x_2}^x
e^{\int_{x}^y b(s,\cdot)ds} c(y,\cdot)dy\;, \ee from which the
result easily follows. The general case can be established by
manipulations similar to the previous ones. \qed

\subsection{Polyhomogeneous solutions of $\partial_\tau\varphi + b
\varphi = c$} \label{Sphgtau}

We pass now to an analysis of ODE's with polyhomogeneous sources.
The results here have an auxiliary character, and several of them
are rather elementary; they will be needed to  handle the real
problem at hand, with partial differential operators.
 Let
$\cO$ be an open subset of $\partial M$, we set
\be\label{defcunotwo} {\mcU}_{x_1} = ]0,x_1]\times \cO\times
[0,T]\;.\ee It will be seen in Sections~\ref{sslwe} and
\ref{Swave} that\footnote{This is due to occurrence of the factor
$\Omega^{(n-1)/2}$ in equations such as \eq{SE.2}.} integer
space-dimensions force us to consider polyhomogeneous expansions
with half-integer power of $x$; in order to account for that, we
introduce an index
$$\delta = {1\over d}\;,$$ where  $d$ is a non-zero integer, $d\in
\N^\ast $. We will mostly be interested in the case $d=1/2$ or
$d=1$, however other values are also possible in the formalism
here. Results analogous to the ones below hold for the general
polyhomogeneous expansions of \Eq{eq:3new}, which can be
established by similar methods. We find it of interest that a
consistent framework can be obtained in the setting considered
below:
\begin{Proposition} \label{proprietetau}
  Let $\beta\in\R$ and consider the system
  \begin{deqarr}& \label{systemtau.0}\partial_\tau\varphi + b \varphi = c \;,&\\ &
\displaystyle  \varphi|_{\{\tau=0\}}(x,v)\equiv
\tilde{\varphi}(x,v) =x^{\beta} \sum_{i=0}^p \sum_{j=0}^{N_i} x^{i
\delta} \ln ^j x \;\tilde{\varphi}_{ij}(x,v) \;+ \;
  \tilde{\varphi}_{p\delta+\beta+\epsilon}(x,v) \;,
 & \nn \\ &\label{systemtau2} &
\\ &\tilde{\varphi}_{ij}\in C_\infty(\overline{
\{\tau=0\}})\;,\qquad
  \tilde{\varphi}_{p\delta+\beta+\epsilon}
\in\mcC_{\infty}^{p\delta+\beta+\epsilon}(\{\tau=0\})\;, &
\label{systemtau}
  \arrlabel{tode}\end{deqarr}
  with
  \begin{deqarr} &
\displaystyle  b (x,v,\tau)= \sum_{i=0}^{p}\sum_{j=0}^{N_i'}
x^{i\delta} \ln^j x \;b_{ij}(x,v,\tau) \: + \:
b_{p\delta+\epsilon}(x,v,\tau)\;, &\\ &
  b_{p\delta+\epsilon} \in \mcC_\infty^{p\delta+\epsilon}(\U_{x_1})\;,\qquad
  b_{ij} \in C_\infty(\overline{\U_{x_1}})
\;, &\\& \displaystyle  c(x,v,\tau) = x^{\beta}\sum_{i=0}^{p}
\sum_{j=0}^{N_i''} x^{i\delta} \ln^j x \;c_{ij}(x,v,\tau) \: + \:
c_{p\delta+\beta+\epsilon}(x,v,\tau)\;, & \\ &
  c_{p\delta+\beta+\epsilon} \in
\mcC_\infty^{p\delta+\beta+\epsilon}(\U_{x_1})\;,\qquad c_{ij} \in
C_\infty(\overline{\U_{x_1}})
\;, & \arrlabel{todep}\end{deqarr}
  where $0<\epsilon<\delta$, and $(N_i),(N'_i),(N''_i)$ are sequences with integer
  values, and with
$$b\in  L^\infty{(\U_{x_1})}\;.$$
  Then the solution $\varphi$ takes the form
  \be \varphi(x,v,\tau) = x^{\beta}\sum_{i=0}^{p}
  \sum_{j=0}^{M_i} x^{i \delta} \ln^j x \;\varphi_{ij}(x,v,\tau) \: + \: \varphi_{p\delta+\beta+\epsilon}(x,v,\tau)\;,
 \ee
  with $\varphi_{ij}\in C_\infty(\overline{\U_{x_1}})$,
  $M_k$ is an integer sequence and  $\varphi_{p\delta+\beta+\epsilon}\in
  \mcC_\infty^{p\delta+\beta+\epsilon}(\U_{x_1})$.
\end{Proposition}

To prove the proposition we shall need the following lemma:
\begin{Lemma}
  Under the hypotheses of Proposition~\ref{proprietetau}, suppose that
  in addition  we have
$$ \tilde{\varphi}_{p\delta+\beta+\epsilon}= b_{p\delta+\epsilon} =
c_{p\delta+\beta+\epsilon}=0\;.$$ Then for any
 $\epsilon\in ]0,\delta[$ we have
 \be
 \varphi =x^{\beta} \sum_{i=0}^{p} \sum_{j=0}^{M_i} x^{ i \delta}
 \ln^j x \;\varphi_{ij} \: + \varphi_{p\delta+
\beta+\epsilon}\;,\label{formephi} \ee  with $\varphi_{ij}\in
C_\infty
(\overline {\U_{x_1}})$, 
 $\varphi_{p\delta + \beta+\epsilon}\in\mcC_{\infty}^{p \delta +\beta+\epsilon}
 (\U_{x_1})$,
for some integer-valued sequence $M_k$.
\end{Lemma}
\proof Inserting (\ref{formephi}) in the equation
(\ref{systemtau.0}) and tracking  the coefficients in front of
$x^{i \delta}\ln^j x$ one finds the following set of equations:
\beqan M_0 =\max \{N_0,N_0''\}\;, & & M_{i+1} = \max \{
\max_{0\leq k\leq i} M_k+
 N'_{i-k} , \; N''_{i+1},\; N_{i+1}\} \;,\\
i \in \lsemantics0,p\rsemantics \;,\; j\in \lsemantics
 0,M_i\rsemantics \;,&&\partial_\tau \varphi_{ij} + \sum_{k=0}^i
 \sum_{l=0}^{\min\{N_k', j\}}
 b_{kl}\varphi_{i-k\;j-l} = c_{ij} \;,\\
\partial_\tau\varphi_{p\delta+\beta+\epsilon} + b
 \varphi_{p\delta+\beta+\epsilon} &=&-\sum_{i=p+1}^{2p}
 \;x^{\beta}\sum_{j=0}^{M_i} x^{i\delta}\ln^j x \;\{ \sum_{k=0}^i
 \sum_{l=0}^{\min\{N_k', j\}}
b_{kl}\varphi_{i-k\;j-l} \} \;. \eeqan Here $\lsemantics a,b
\rsemantics:= [a,b]\cap \N$. This system is easily solved: one
begins with $i=0$ and solves the equations for $j$  running  from
$0$ to $M_0$. This can then be repeated for $i=1$, {\em etc},
until $i=p$ is reached. This provides the functions
$\varphi_{ij}$.  Finally, one solves the last equation for the
remainder term $\varphi_{p\delta + \beta +\epsilon}$, with initial
value zero, noting that the right hand side  of the resulting
equation is in $\mcC_\infty^{p\delta+\beta+\epsilon}(\U_{x_1})$,
and one concludes using Proposition~\ref{Plemme}. \proofend.
\medskip

\noindent{\sc Proof of Proposition \ref{proprietetau}}: With the
notation of the proposition, we set $b_{\phg} = b -
b_{p\delta+\epsilon}$, $c_{\phg}= c-c_{p\delta+\beta+\epsilon}$,
$\tilde \varphi_{\phg} = \tilde \varphi - \tilde
\varphi_{p\delta+\beta+\epsilon} $.  We use the Lemma above to
obtain a solution $\varphi_{\phg}$  of the problem  \beqa
  \partial_\tau\varphi + b_{\phg} \varphi &=&
  c_{\phg}\;,\label{equalemmetau}\\
  \varphi|_{\Sigma}= \tilde{\varphi} &=& x^{\beta}\sum_{i=0}^p
  \sum_{j=0}^{N_i} x^{i\delta} \ln ^j x \;\tilde{\varphi}_{ij}(x,v) \;.
  \eeqa  Then we
solve $$\partial_{\tau}\varphi'+  b\varphi' =
  c_{p\delta+\beta+\epsilon}- b_{p\delta+\epsilon}\varphi_{\phg}$$ with
$\varphi'|_{\tau=0} = \tilde \varphi_{p\delta+\beta+\epsilon}$.
According to
  Proposition~\ref{Plemme} we have $\varphi'\in
  \mcC_{\infty}^{p\delta+\beta+\epsilon}(\U_{x_1})$.
To conclude we set $\varphi = \varphi_{\phg} + \varphi'$ which is
of
  the required form,  and solves (\ref{systemtau.0}). \qed

\subsection{Polyhomogeneous solutions of $\partial_x\varphi + b
\varphi = c$} \label{Sphgx}

\begin{Proposition}\label{propositionx}
Let $\varphi $ be a solution in $C_\infty ^{\loc} (\U_{x_1})$ of
\be
\partial_x\varphi + \frac bx \varphi =
c\;,\label{systemx} \ee and suppose that \eq{todep} holds with
some $\epsilon\in ]0,\delta[$, $\beta\in \R $,
and with some integer-valued sequences $(N'_i),(N''_i)$. If
$$ b=o(x)$$
(equivalently, $b_{0j}(0,v,\tau)=0$), then \be \varphi =
\sum_{i=0}^{p} \sum_{j=0}^{M_i} x^{i\delta} \ln^j x
\;\widehat\varphi_{ij} \: + x^{\beta+1} \sum_{i=0}^{p}
\sum_{j=0}^{M_i} x^{i\delta} \ln^j x \;\varphi_{ij} \: + \:
\varphi_{p\delta+1+\beta+\epsilon}\;, \label{formephix} \ee with
$$\widehat\varphi_{ij} \:,\varphi_{ij}\in
C_\infty(\overline{\U_{x_1}})\;, \qquad
\varphi_{p\delta+1+\beta+\epsilon}\in
\mcC_\infty^{p\delta+1+\beta+\epsilon}({\U_{x_1})}\;,$$ for some
integer sequence $(M_i)$.
\end{Proposition}
\proof Proposition~\ref{propositionx1} shows that for $\beta > -1$
the limit $$\varphi_0(\cdot):=\lim_{x\to 0} \varphi(x,\cdot)$$
exists and is a smooth function on $\cO\times[0,T]$. When $b$ is a
multiple of the identity matrix the result is obtained by a
straightforward analysis of the formula \be \label{intrep1}
\varphi(x,\cdot) = \varphi_0(\cdot) e^{-\int_{0}^x b(s,\cdot)ds} +
\int_{0}^x e^{\int_{x}^y b(s,\cdot)ds} c(y,\cdot)dy\;, \ee using
the estimates of Lemma~\ref{integrationx}. For $\beta < -1 $, and
again for $b$ --- a multiple of the identity matrix --- we use
instead \be \label{intrep1+} \varphi(x,\cdot) =
\varphi(x_1/2,\cdot) e^{-\int_{x_1/2}^x b(s,\cdot)ds} +
\int_{x_1/2}^x e^{\int_{x}^y b(s,\cdot)ds} c(y,\cdot)dy\;. \ee In
the general case, we first note that it follows from
Proposition~\ref{propositionx1} that there exists $\lambda\in\R$
such that $\psi\in \mcC^\lambda_\infty$. We then write
\be\label{insback}\partial_x \psi - c= - \frac bx \psi \in
\mcC^{\lambda+\delta-1}_\infty\;;\ee integrating gives
$$\psi - \int_0^xc  \in \mcC^{\lambda+\delta}_\infty\;.$$
Inserting this equation   in the right-hand-side of \eq{insback}
and integrating again one obtains a similar equation with a
remainder term falling-off one power of $\delta$ faster. The
result is proved by repeating this procedure a finite number of
times.
 \qed

\bibliographystyle{amsplain}

\addcontentsline{toc}{1}{\bigskip \noindent  {\bf References \hfill}}

\bibliography{../../references/hip_bib,%
../../references/reffile,%
../../references/vienna,%
../../references/newbiblio,%
../../references/newbiblio2,%
../../references/netbiblio,%
../../references/bibl,%
../../references/marot}
\end{document}